\documentclass[11pt]{article}

\usepackage[margin=1in]{geometry}
\usepackage{amsmath,amssymb,amsthm,mathtools,bm}
\usepackage{booktabs}
\usepackage{graphicx}
\usepackage{subcaption}
\usepackage{placeins}
\usepackage{enumitem}
\usepackage{microtype}
\usepackage{hyperref}

\graphicspath{{./}}
\newcommand{\meshimage}[2][\linewidth]{\includegraphics[width=#1]{#2}}
\newcommand{\resultpanel}[1]{\meshimage[\linewidth]{#1}}

\hypersetup{
  hidelinks,
  pdfauthor={Xianping Li},
  pdftitle={A Unified Isoparametric Moving Mesh PDE Method for P1 and Q1 Elements with Certified Local Nonsingularity}
}

\newtheorem{theorem}{Theorem}[section]
\newtheorem{lemma}[theorem]{Lemma}
\newtheorem{proposition}[theorem]{Proposition}
\newtheorem{corollary}[theorem]{Corollary}
\theoremstyle{definition}

\newtheorem{remark}[theorem]{Remark}
\theoremstyle{definition}

\newcommand{\R}{\mathbb{R}}
\newcommand{\T}{\mathcal{T}}

\newcommand{\Khat}{\widehat K}
\newcommand{\x}{\bm{x}}
\newcommand{\xiB}{\bm{\xi}}
\newcommand{\M}{\mathbb{M}}
\newcommand{\B}{\bm{B}}
\newcommand{\A}{\bm{A}}
\newcommand{\J}{\bm{J}}
\newcommand{\HH}{\bm{H}}
\newcommand{\I}{\bm{I}}
\newcommand{\Tframe}{\bm{T}}
\newcommand{\PiB}{\bm{\Pi}}
\newcommand{\cof}{\operatorname{cof}}
\newcommand{\tr}{\operatorname{tr}}
\newcommand{\rank}{\operatorname{rank}}
\newcommand{\diag}{\operatorname{diag}}
\newcommand{\dd}{\,\mathrm{d}}
\newcommand{\EA}{\mathrm{EA}}
\newcommand{\val}{\mathrm{val}}

\newcommand{\dist}{\operatorname{dist}}

\title{A Unified Isoparametric Moving Mesh PDE Method for $P_1$ and $Q_1$ Elements with Certified Local Nonsingularity}
\author{Xianping Li%
\thanks{College of Integrative Sciences and Arts, Arizona State University, AZ, 85212, USA (\textit{xianping@asu.edu})}
}
\date{}

\begin{document}
\maketitle

\begin{abstract}
We develop a unified isoparametric moving mesh partial differential equation
(MMPDE) method for $P_1$ simplicial and $Q_1$ tensor-product elements,
including elements on embedded surfaces. Mesh movement follows a prescribed
symmetric positive-definite metric through equidistribution and alignment.
For $Q_1$ elements, we prove that bounded, exactly integrated
adaptation energy prevents local Jacobian rank loss on a fixed finite mesh
at sharp critical exponents. The result also gives local full-rank preservation for
embedded bilinear surface elements. A complete physical-coordinate gradient
accounts for metric variation at moving quadrature points. For the discrete MMPDE,
we use whole-element validity tests and develop sufficient Bernstein certificates
for trilinear hexahedra. An explicit positive-margin
bound on each fixed stage is provided by an optional certificate barrier.
Numerical examples illustrate metric-driven
movement on one-dimensional, planar, volume, and graph-surface meshes. Limitations of
pointwise hexahedral Jacobian sampling are also presented.
\end{abstract}


\section{Introduction}\label{sec:intro}

Moving mesh methods relocate mesh points without changing mesh connectivity,
according to a prescribed measure of size, shape, and direction. In several
dimensions, this information is commonly specified by a symmetric positive-
definite (SPD) metric tensor $\M=\M(\x)$. The metric may be designed for
approximation accuracy, diffusion-aligned elements,
or another geometric or application-specific purpose. The metric determines preferred element size,
shape, and orientation through the equidistribution and alignment conditions
\cite{Huang2001,HuangSun2003,Huang2005quality,HuangRussell2011}. The moving mesh PDE
(MMPDE) approach describes the evolution of mesh points by a system of differential
equations derived from a meshing functional
\cite{HuangRenRussell1994,HuangRussell2001,BuddHuangRussell2009}.

For simplicial meshes, the geometric discretization of variational meshing
functionals and the nonsingularity of the resulting MMPDE trajectories have been
studied in \cite{HuangKamenski2015,HuangKamenski2018}. Extensions to surface meshes
are considered in \cite{KolasinskiHuang2020}, and a unified treatment of curves,
surfaces, and domains using affine simplicial elements is given in
\cite{ZhangHuang2026}. In those works, only simplicial meshes are considered. 
In their formulations, an affine element map has a constant
Jacobian matrix, a property that is especially useful in the analysis of element
volume and nonsingularity.

On the other hand, quadrilateral and hexahedral meshes are widely used in finite element computations,
particularly when mesh lines are to be aligned with boundaries or directional
features of a solution. Their approximation properties and practical advantages
depend on the problem, discretization, and mesh quality, as discussed in
\cite{DAzevedo2000,CifuentesKalbag1992,SchneiderEtAl2019}. Generation and
optimization of these meshes continue to present difficulties
\cite{BommesEtAl2013,PietroniEtAl2022,Knupp2000a,Knupp2000b}.
Unlike an affine simplex, a bilinear quadrilateral or trilinear hexahedron has an
isoparametric Jacobian that generally varies within the element. Consequently,
a positive determinant at a finite set of sample points does not necessarily
ensure a positive determinant everywhere for a trilinear hexahedron. Meanwhile, the use of
Bernstein polynomial bounds for verifying geometric validity have been studied
in \cite{JohnenRemacleGeuzaine2013,JohnenWeillRemacle2017}.
For trilinear hexahedra, a recent boundary-based characterization also permits
whole-element validity to be checked without using a full-volume Bernstein
representation \cite{Zhang2026HexBoundaries}.
Classical quadrilateral and hexahedral smoothing, untangling, and
Jacobian-quality optimization are discussed in
\cite{Knupp2000a,Knupp2000b,RuizGironesEtAl2015}. Variational and
harmonic-coordinate mesh movement provide complementary routes to
$r$-adaptation \cite{CaoHuangRussell1999,HuangRussell2011}. The present
barrier acts on coefficients that bound the entire determinant polynomial,
rather than only on determinants sampled at a finite set of quadrature points.

In this work, we develop a unified isoparametric MMPDE formulation for $P_1$
simplicial and $Q_1$ tensor-product elements. Our objective is to preserve the
equidistribution--alignment structure of the variational MMPDE while controlling
local validity of non-affine elements. The method uses the intrinsic pullback
metric $\B_K^T\M\B_K$, where $\B_K$ is the Jacobian of the element map, and
projects nodal velocities onto the admissible geometric tangent spaces.

A principal theoretical contribution is the nonsingularity analysis for
non-affine $Q_1$ elements. We prove that, on a fixed finite mesh with
bounded nodal coordinates and a uniformly SPD metric, bounded exactly
integrated adaptation energy prevents local Jacobian rank loss. The
critical exponents are $p=3/2$ for bilinear and $p=4/3$ for trilinear
elements, and both thresholds are sharp. This extends the coercivity-based
analysis of \cite{HuangKamenski2018} from affine simplicial meshes to the
present isoparametric $Q_1$ formulation. The result also establishes
local full-rank preservation for embedded bilinear surface elements,
although it does not guarantee prescribed projected orientation or
global injectivity.

For the discrete MMPDE, we use whole-element validity checks to preserve
local orientation at accepted mesh configurations. Planar $Q_1$
elements admit an exact four-corner determinant test, whereas trilinear
hexahedra are verified using sufficient Bernstein coefficient bounds.
An optional certificate barrier provides an explicit positive-margin
estimate along a continuous gradient flow within a fixed certificate
stage. These mechanisms complement the exact-energy nonsingularity
theorem because finite numerical quadrature alone does not provide its
guarantee.

Our proposed method accepts any prescribed SPD metric. We use nodal quadratic least-squares
Hessian recovery for both element families without subdividing $Q_1$ cells into simplices.

The rest of the paper is organized as follows. 
Section~\ref{sec:framework} develops the isoparametric MMPDE, and
Section~\ref{sec:validitymain} presents the validity theory and discrete
integration. Metric construction and representative examples follow in
Sections~\ref{sec:adaptationmain} and~\ref{sec:numerics}.
Conclusions are presented in Section~\ref{sec:conclusion}.

\section{Isoparametric moving mesh formulation}\label{sec:framework}
In this section, we develop the isoparametric moving mesh PDE. 
We first describe the element geometry and the geometric constraints on moving vertices. 
\subsection{Isoparametric geometry}\label{sec:iso}

Let $S$ be a connected $m$-dimensional geometric object embedded in $\R^d$, where
\[ 1\le m\le d. \]
The principal cases are a curve ($m=1$), a surface ($m=2<d$), and a bulk domain ($m=d$).  
Let $\T_h$ be a conforming mesh of $S$ with vertices $\{\x_i\}_{i=1}^{N_v}$ and elements $K$.
Let $N$ denote the number of elements in $\T_h$.

For the tensor-product family considered here, the reference element is
\[ \Khat=[0,1]^m. \]
A $Q_1$ element has $2^m$ local vertices indexed by $\bm{\alpha}\in\{0,1\}^m$.  
The tensor-product basis functions are
\begin{equation}\label{eq:q1basis}
N_{\bm{\alpha}}(\xiB) =\prod_{r=1}^m \xi_r^{\alpha_r}(1-\xi_r)^{1-\alpha_r}.
\end{equation}
The element map is given by
\begin{equation}\label{eq:map}
F_K(\xiB)=\sum_{\bm{\alpha}\in\{0,1\}^m} N_{\bm{\alpha}}(\xiB)\x_{K,\bm{\alpha}}, \qquad \xiB\in\Khat,
\end{equation}
and its Jacobian matrix is as follows:
\begin{equation}\label{eq:Bdef}
\B_K(\xiB)=D_{\xiB}F_K(\xiB) =\sum_{a=1}^{2^m}\x_{K,a}\bigl(\nabla_{\xiB}N_a(\xiB)\bigr)^T \in\R^{d\times m}.
\end{equation}
For $m=d$, regularity means $\det\B_K\ne0$. An oriented bulk element is
valid when $\det\B_K>0$ relative to the chosen ordering of the reference
vertices.  For $m<d$, regularity means $\rank\B_K=m$.  Local full rank does
not, by itself, imply that $F_K$ is injective on the entire reference cell.

The formulas below also apply to the $P_1$ family, with a regular (equilateral) reference $m$-simplex 
of unit edge length and volume $V_{\rm ref}=\sqrt{m+1}/(m!\sqrt{2^m})$. 
Its affine map has a constant Jacobian, and $\A=\alpha\I_m$ characterizes a metric-regular simplex. 
For $m=1$, affine $P_1$ and $Q_1$ coincide. The reference-shape convention must be 
accounted for in alignment comparisons; see \cite{HuangKamenski2015,ZhangHuang2026}.

For $m<d$, assume that the prescribed geometry is a sufficiently smooth embedded
manifold (possibly with boundary and corners treated as separate strata), and
mesh vertices are required to remain on it.  The restriction is imposed on the
vertices only and a polynomial element spanning these vertices need not lie
exactly on the smooth manifold between nodes.  Let
\[ \Tframe_i=[\bm t_{i,1},\ldots,\bm t_{i,m}]\in\R^{d\times m} \]
be an orthonormal basis of the tangent space of $S$ at $\x_i$.  Define the orthogonal tangent projector
\begin{equation}\label{eq:Pi}
\PiB_i=\Tframe_i\Tframe_i^T.
\end{equation}
For a bulk domain, take $\PiB_i=\I_d$.  Then we have
\[ \PiB_i=\begin{cases} \I_d, & m=d,\\ \I_d-\bm n_i\bm n_i^T, & m=d-1,\\ \bm\tau_i\bm\tau_i^T, & m=1, \end{cases} \]
whenever a single normal $\bm n_i$ or tangent $\bm\tau_i$ is available.  
This packages the separate projection cases in \cite{ZhangHuang2026} into one formula.

Boundary restrictions can be incorporated by replacing $\PiB_i$ with the projector onto the allowable boundary tangent space 
or by setting $\PiB_i=0$ at fixed vertices.

\subsection{Metric equidistribution and alignment}\label{sec:metric}
With the isoparametric geometry established, we formulate equidistribution and alignment in the prescribed metric. 
Let $U\subset\R^d$ be a neighborhood containing all element images
$F_K(\Khat)$ in the configurations under consideration.  For a bulk domain,
$U$ may be the physical domain together with a neighborhood of its boundary.
For a smooth embedded surface or curve, $U$ is an ambient tubular
neighborhood.  Assume an ambient extension of the desired metric is defined
on $U$ and is symmetric and uniformly positive definite:
\begin{equation}\label{eq:Mbound}
0<\underline\rho\,\I_d\preceq \M(\x)\preceq \overline\rho\,\I_d,
\qquad \x\in U.
\end{equation}
An ambient extension is necessary because the interior of a flat or bilinear
polynomial surface element generally lies off the exact curved surface,
even when its vertices lie on that surface.

For each element $K$ and each reference point $\xiB\in\Khat$, define the intrinsic pullback metric by
\begin{equation}\label{eq:Adef}
\A_K(\xiB) =\B_K(\xiB)^T\M(F_K(\xiB))\B_K(\xiB) \in\R^{m\times m}.
\end{equation}
If $\B_K$ has rank $m$, then $\A_K$ is symmetric positive definite.

The local metric $m$-volume density is defined as
\begin{equation}\label{eq:metricdensity}
\mu_K(\xiB)=\det(\A_K(\xiB))^{1/2}.
\end{equation}
The metric measure of $K$ is
\begin{equation}\label{eq:metricvolume}
|K|_{\M} =\int_{\Khat}\mu_K(\xiB)\dd\xiB.
\end{equation}
For $m=d$, one has
\[ \mu_K=|\det\B_K|\det(\M)^{1/2}. \]
For a surface element in $\R^3$, $\mu_K$ is the corresponding metric surface-area density.

Let
\begin{equation}\label{eq:sigma}
\sigma_h=\sum_{K\in\T_h}|K|_{\M}.
\end{equation}
For equal-volume reference cells, the elementwise equidistribution target is given by
\begin{equation}\label{eq:eqdist}
|K|_{\M}=\frac{\sigma_h}{N}, \qquad K\in\T_h.
\end{equation}
Equation~\eqref{eq:eqdist} is an elementwise condition since $|K|_{\M}$ is the integral of the metric density over $K$. 
A stronger, pointwise condition would require $\mu_K(\xiB)$ to be approximately uniform throughout the reference elements. 
These conditions are not equivalent for non-affine $Q_1$ maps because the density can vary within a cell.

In addition to equidistribution, the metric prescribes the local element shape. 
The pointwise tensor-product alignment condition takes the form
\begin{equation}\label{eq:alignment}
\frac1m\tr(\A_K) =\det(\A_K)^{1/m}.
\end{equation}
By the arithmetic--geometric mean inequality, equality holds if and only if $\A_K=\alpha\I_m$ for some $\alpha>0$. 
Equivalently, the columns of $\M^{1/2}\B_K$ are mutually orthogonal and have equal length.

A useful pointwise alignment quality measure is
\begin{equation}\label{eq:Qali}
Q_{\rm ali}(K,\xiB) =\frac{\tr(\A_K)}{m\det(\A_K)^{1/m}} \ge 1.
\end{equation}
An ideal metric element has $Q_{\rm ali}=1$. On the tensor-product
reference cell this corresponds locally to a metric square ($m=2$) or
cube ($m=3$). For $P_1$, it describes a regular metric simplex only
when the regular reference simplex specified above is used.

In the bulk case $m=d$, if
\[ \M=R\diag(\lambda_1,\ldots,\lambda_m)R^T \]
and the element coordinate directions are approximately aligned with the eigenvectors,
then $\A_K\approx\alpha I$ implies the local physical scales
\begin{equation}\label{eq:hi}
h_i\propto \lambda_i^{-1/2}.
\end{equation}
Thus, physical element sizes decrease in directions where the metric has larger eigenvalues, 
producing anisotropic elements aligned with the prescribed metric.

\subsection{Integrated equidistribution--alignment energy}\label{sec:energy}
The preceding metric conditions motivate an energy functional whose nodal variation defines the mesh motion.
Following the equidistribution--alignment functional used for simplicial meshes \cite{Huang2001,HuangKamenski2015,ZhangHuang2026}, define
\begin{equation}\label{eq:Jdef}
\J_K(\xiB)=\A_K(\xiB)^{-1}, \qquad D_K(\xiB)=\det\J_K(\xiB).
\end{equation}
For parameters $p>1$ and $0<\theta\le1/2$, let
\begin{equation}\label{eq:G}
G(\J,D) =\theta D^{-1/2}[\tr(\J)]^{mp/2} + (1-2\theta)m^{mp/2}D^{(p-1)/2}.
\end{equation}
The default values used in the underlying simplicial work are
\begin{equation}\label{eq:defaults}
p=\frac32, \qquad \theta=\frac13.
\end{equation}
For $0<\theta\le1/2$ both terms in \eqref{eq:G} are
nonnegative, which is used in the barrier bound. The parameter $p>1$ determines the strength of 
the metric-alignment and equidistribution penalties and is essential in the scaling and nonsingularity analysis. 

The isoparametric adaptation energy is defined by
\begin{equation}\label{eq:IEA}
I_{\EA}(\T_h) =\sum_{K\in\T_h} \int_{\Khat} G(\J_K(\xiB),D_K(\xiB))\dd\xiB.
\end{equation}
For an affine element with a constant stage metric, the integrand is constant
and \eqref{eq:IEA} equals $|\Khat|G$. The factor $|\det\B_K|$ associated
with physical volume is already present implicitly in this density.
Indeed, for an orientation-preserving bulk map ($m=d$), put
$\J_{\xi x}=\B_K^{-1}$ and
$\HH=\J_{\xi x}\M^{-1}\J_{\xi x}^T$,
together with $D_K=\det\HH=(\det\B_K)^{-2}(\det\M)^{-1}$, we have
\begin{equation}\label{eq:pullbackidentity}
G(\HH,\det\HH)
=\det\B_K\,\mathcal G(\J_{\xi x},\M),
\end{equation}
where the classical physical-coordinate density is
\begin{align}\label{eq:physicaldensity}
\mathcal G(\J_{\xi x},\M)
={}&\theta\sqrt{\det\M}
  [\tr(\J_{\xi x}\M^{-1}\J_{\xi x}^T)]^{mp/2}\notag\\
&+(1-2\theta)m^{mp/2}(\det\J_{\xi x})^p
 (\det\M)^{(1-p)/2}.
\end{align}
Therefore, when the coordinate maps and prescribed metric agree, the
reference-cell energy in \eqref{eq:IEA} equals the classical energy integrated
over the physical domain. This equality of energies does not imply that
the corresponding gradient flows in physical and computational coordinates
are identical. For embedded elements, the change of variables uses the
intrinsic volume density and the induced tangent metric. We derive the
variation with respect to physical nodal coordinates below, after
specifying how the metric is evaluated during each relaxation stage.

\begin{remark}[Exact integration versus numerical quadrature]\label{rem:exactvsquadrature}
For non-affine $Q_1$ elements, $\B_K$ varies within the reference cell, so
\eqref{eq:IEA} is defined by an exact integral. Numerical calculations use
finite quadrature, initially with at least two Gaussian nodes per coordinate
and higher-order rules for convergence checks. The exact-energy result below
does not apply automatically to a finite quadrature sum, which can miss a
near-singular region. Moreover, the density depends on $\B_K$ through
$\B_K^T\M\B_K$ and does not distinguish the signs of bulk Jacobian
determinants. A discrete method therefore requires a separate oriented,
whole-element validity check.
\end{remark}

The following lemma shows that bounded adaptation energy alone does not bound the nodal coordinates, 
so compactness in the subsequent theorem must be imposed separately.

\begin{lemma}[Scaling under nodal dilation for a constant metric]
\label{lem:Q1nodalscaling}
For $\M=\I$ and $p>1$, if every element map is replaced by
$F_K^{(s)}=sF_K$ with $s>0$, then the exactly integrated energy
satisfies
\[
 I_{\EA}(sX)=s^{m(1-p)}I_{\EA}(X),
\]
whenever the original configuration has full-rank element maps.
\end{lemma}
\begin{proof}
Nodal dilation sends $\B_K$ to $s\B_K$ and $\A_K$ to
$s^2\A_K$. Consequently $\J_K$ changes to $s^{-2}\J_K$ and
$D_K$ to $s^{-2m}D_K$. Each term of $G$ in \eqref{eq:G}
therefore acquires the factor $s^{m(1-p)}$. Integrating over
the unchanged reference cells proves the identity.
\end{proof}

We next show that, under a compactness assumption, bounded exact adaptation energy prevents 
bilinear and trilinear elements from losing local Jacobian rank.

\begin{theorem}[Exact-energy nonsingularity of $Q_1$ elements]
\label{thm:Q1exactcoercivity}
Let $m\in\{2,3\}$ and $d\ge m$. Consider a fixed finite collection
of bilinear ($m=2$) or trilinear ($m=3$) isoparametric maps
$F_K:[0,1]^m\to\R^d$. Assume that the nodal-coordinate vector $X$
ranges over a compact set $\mathcal K$ and that the metric values
throughout the corresponding element images satisfy the uniform bounds
\eqref{eq:Mbound}. Let $0<\theta\le1/2$, and suppose
\begin{equation}\label{eq:Q1criticalp}
p\ge p_m:=\begin{cases}3/2,&m=2,\\4/3,&m=3.\end{cases}
\end{equation}
The energy in \eqref{eq:IEA} is understood as an extended nonnegative
integral, with integrand $+\infty$ where the element Jacobian loses
column rank. For every $E<\infty$ there exists $\delta>0$ such that
\begin{equation}\label{eq:Q1sigmabound}
X\in\mathcal K,\quad I_{\EA}(X)\le E
\quad\Longrightarrow\quad
\min_K\min_{\xiB\in\Khat}\sigma_{\min}(\B_K(\xiB;X))\ge\delta.
\end{equation}
The conclusion applies to embedded bilinear and trilinear maps as a statement of local full column rank. 
The assumption that $X$ belongs to a compact set is essential and is not implied by a bound on the adaptation energy. 
The constant $\delta$ depends on the fixed finite element collection, $\mathcal K$, $E$, the metric bounds, $p$, and $\theta$. 
The result does not provide an explicit lower bound or one that is uniform under mesh refinement. The thresholds in
\eqref{eq:Q1criticalp} cannot be lowered for the stated general
class of maps.
\end{theorem}
\begin{proof}
For clarity we separate the singular-value estimate, the geometric
nonintegrability argument, and the compactness step.
At a full-column-rank Jacobian $\B$, denote its singular values by
$\sigma_1\ge\cdots\ge\sigma_m>0$ and put
\[
H_p(\B)=\left(\prod_{j=1}^m\sigma_j\right)
\left(\sum_{j=1}^m\sigma_j^{-2}\right)^{mp/2}.
\]
For every $v\in\R^m$, the metric bounds give
$\underline\rho|\B v|^2\le v^T\A v\le\overline\rho|\B v|^2$. Equivalently,
$\underline\rho\,\B^T\B\preceq\A\preceq\overline\rho\,\B^T\B$.
Eigenvalue monotonicity for symmetric matrices then gives
\[
\sqrt{\det\A}\ge\underline\rho^{m/2}\prod_j\sigma_j,
\qquad
\tr(\A^{-1})\ge\overline\rho^{-1}\sum_j\sigma_j^{-2}.
\]
The second term of $G$ is nonnegative for $0<\theta\le1/2$. Keeping only the first term gives
\begin{equation}\label{eq:Q1energyLB}
G(\J_K,D_K)\ge c_0H_p(\B_K),\qquad
c_0=\theta\underline\rho^{m/2}\overline\rho^{-mp/2}>0.
\end{equation}
All constants in this estimate are independent of the configuration
$X\in\mathcal K$ and the reference point.

Lemma~\ref{lem:Q1singularenergy}, proved in
Appendix~\ref{app:Q1coercivity}, establishes the essential geometric
fact: if a bilinear or trilinear map has a rank-deficient Jacobian at
any point of its closed reference cube, then
$\int_{\Khat}H_p(\B)\dd\xiB=+\infty$ for the corresponding
exponent in \eqref{eq:Q1criticalp}. This includes singularities at
faces, edges, and corners. In particular, a rank-one trilinear
Jacobian creates a segment of singular Jacobians, a fact needed to
obtain the sharp three-dimensional exponent.

Suppose that \eqref{eq:Q1sigmabound} is false. Then a sequence
$X_n\in\mathcal K$ satisfies $I_{\EA}(X_n)\le E$ while the minimum
of $\sigma_{\min}(\B_K)$ over elements and reference points tends
to zero. Because the element collection is finite and both
$\mathcal K$ and $\Khat$ are compact, pass to a subsequence on
which one fixed element $K_*$ and reference points $\xiB_n\to\xiB_*$
realize that minimum and $X_n\to X_*\in\mathcal K$.
The $Q_1$ Jacobian depends continuously on the nodal coordinates and
reference points, so
$\rank\B_{K_*}(\xiB_*;X_*)<m$.

Extend $H_p$ by $+\infty$ at rank-deficient matrices. The
arithmetic--geometric mean inequality yields
\begin{equation}\label{eq:Q1lowersemicont}
H_p(\B)\ge m^{mp/2}
\left(\prod_{j=1}^m\sigma_j\right)^{1-p}.
\end{equation}
On the full-rank set, $H_p$ is continuous. When full-rank matrices converge to a rank-deficient matrix, 
their singular-value product tends to zero. Since $p>1$, the right-hand side of \eqref{eq:Q1lowersemicont} diverges, 
proving that the extension of $H_p$ by $+\infty$ is lower semicontinuous. Fatou's lemma, the geometric lemma, and
\eqref{eq:Q1energyLB} now imply
\[
+\infty
=\int_{\Khat}H_p(\B_{K_*}(\xiB;X_*))\dd\xiB
\le\liminf_{n\to\infty}
\int_{\Khat}H_p(\B_{K_*}(\xiB;X_n))\dd\xiB
\le E/c_0,
\]
a contradiction. The sharpness examples are given in
Appendix~\ref{app:Q1coercivity}.
\end{proof}

The sublevel-set result applies to continuous mesh trajectories that retain an upper bound on the same exact energy.

\begin{corollary}[Exact-energy preservation of local $Q_1$ regularity]
\label{cor:Q1exactflow}
Let $X(t)$ be a continuous path in $\mathcal K$ starting from a
full-rank $Q_1$ mesh, and suppose the same exact energy
\eqref{eq:IEA} satisfies $I_{\EA}(X(t))\le E<\infty$ throughout
the path. Under the hypotheses of Theorem~\ref{thm:Q1exactcoercivity},
all element maps remain full rank with the common bound
\eqref{eq:Q1sigmabound}. In the bulk case $m=d$, initially positive
Jacobians remain positive and satisfy
$\det\B_K\ge\delta^m$ everywhere. For embedded bilinear surface
cells ($m=2<d$), only local full rank is implied. Prescribed projected
orientation and global nonoverlap require additional conditions.
\end{corollary}
\begin{proof}
Theorem~\ref{thm:Q1exactcoercivity} applies at every time, with the same $\delta$. 
For each bulk element the function $(t,\xiB)\mapsto\det\B_K(\xiB;X(t))$ is continuous 
and nonzero on the connected set $[0,T]\times\Khat$. Since it is positive at $t=0$, it remains positive. 
Moreover, $|\det\B_K|=\prod_{j=1}^m\sigma_j(\B_K)\ge\delta^m$.
For embedded maps, the same estimate implies full column rank but no
fixed orientation of a projected determinant.
\end{proof}

\begin{remark}[Role of the exact energy]
The conclusion applies to a continuously dissipative MMPDE stage if
the exact $I_{\EA}$ is nonincreasing, or if it is bounded above by a
nonincreasing total energy with a nonnegative barrier. Such dissipation
requires the gradient of that same exact functional. Differentiating
a quadrature sum, changing the prescribed metric between stages, or
using a distinct indirect-coordinate gradient flow does not by itself
establish the hypothesis. The result provides a conditional sublevel-set bound. The nodal
compactness assumption must be verified independently. Lemma~\ref{lem:Q1nodalscaling} shows this 
by uniformly dilating all nodal positions when the metric is constant. 
Such a dilation does not establish that an individual interior node 
can escape while the boundary nodes remain fixed. A sufficient
condition is to keep every node in a prescribed compact region,
with the metric uniformly SPD on all element images. Confinement to
a compact embedded surface also bounds the nodal positions, whereas
fixing only the boundary nodes of a bulk mesh does not in general
bound interior nodes. At the critical exponents, the examples in
Appendix~\ref{app:Q1coercivity} have logarithmic energy divergence,
so the nonconstructive constant $\delta$ is not a practical discrete
validity margin. Unlike the explicit simplicial bounds in
\cite{HuangKamenski2018}, the present estimate is qualitative for a
fixed collection of non-affine elements. With additional compactness
and local-Lipschitz hypotheses it can support an ODE continuation
argument, but it does not alone imply global existence or replace
discrete whole-element acceptance.
\end{remark}

A metric stage begins with a prescribed physical-space SPD tensor field
$\M^{(n)}:U\to\R^{d\times d}$, constructed from the physical solution or
specified analytically.  The field is fixed during this stage and
its definition is not updated in response to the changing mesh.  There are,
however, two different ways to evaluate it during a direct mesh relaxation.
They define different finite-dimensional meshing energies and different
gradient flows.

In the spatial-field realization, the prescribed metric is evaluated at the moving physical points:
\begin{align}
 \M_K^{\rm sp}(\xiB;X)&=\M^{(n)}(F_K(\xiB;X)),
 \label{eq:spatialmetric}\\
 \A_K^{\rm sp}(\xiB;X)&=\B_K(\xiB;X)^T
                 \M_K^{\rm sp}(\xiB;X)\B_K(\xiB;X),
 \label{eq:Aspatial}\\
 I_{\EA,\rm sp}^{(n)}(X)
 &=\sum_K\int_{\Khat}
 G\left((\A_K^{\rm sp})^{-1},
             \det((\A_K^{\rm sp})^{-1})\right)\dd\xiB.
 \label{eq:IEAspatial}
\end{align}
The field $\M^{(n)}$ is independent of pseudo-time, but
$\M_K^{\rm sp}$ is not independent of $X$.  Its derivative must
therefore be included in the discrete variation.  This is the natural
finite-element discretization of the physical-space metric functional
\eqref{eq:IEA} when the prescribed physical metric is interpreted as a
spatial field.

For comparison, the reference-frozen realization samples the
field at the beginning of the stage and attaches those labels to reference
points throughout the relaxation.  Specifically, define
$\overline\M_K(\xiB)=\M^{(n)}(F_K(\xiB;X^{(n)}))$ and
\begin{equation}\label{eq:Astage}
 \A_K^{\rm ref}(\xiB;X\mid X^{(n)})
 =\B_K(\xiB;X)^T\overline\M_K(\xiB)\B_K(\xiB;X).
\end{equation}
Its stage energy is
\begin{equation}\label{eq:IEAstage}
 I_{\EA,\rm ref}^{(n)}(X)
 =\sum_{K}\int_{\Khat}
 G\left((\A_K^{\rm ref})^{-1},
    \det((\A_K^{\rm ref})^{-1})\right)\dd\xiB.
\end{equation}
The two energies coincide at $X=X^{(n)}$ if sampled consistently,
but their gradients and subsequent stage energies generally differ.
Unless explicitly stated otherwise, the direct
spatial-field analysis below refers to \eqref{eq:IEAspatial}, while
\eqref{eq:IEAstage} remains a controlled alternative.

Treat $\J$ and $D=\det\J$ as independent variables and define
\begin{align}
\HH &:=\frac{\partial G}{\partial\J} = \theta\frac{mp}{2} D^{-1/2}[\tr(\J)]^{mp/2-1}\I_m, 
\label{eq:H}\\
g &:=\frac{\partial G}{\partial D} = -\frac{\theta}{2}D^{-3/2}[\tr(\J)]^{mp/2} + \frac{p-1}{2}(1-2\theta)m^{mp/2}D^{(p-3)/2}.
\label{eq:g}
\end{align}
Set
\begin{equation}\label{eq:S}
\bm S=\J\HH\J+Dg\J.
\end{equation}
Throughout the following matrix differentiation, \(X:Y=\operatorname{tr}(X^T Y)\) denotes the Frobenius inner product. 
Thus, for a scalar-valued function \(F(X)\) of a matrix argument, we can write
\[ \dd F=\frac{\partial F}{\partial X}:\dd X. \]
We first differentiate the integrand with respect to the element Jacobian while holding the metric value fixed.

\begin{proposition}[Derivative with respect to the isoparametric Jacobian]\label{prop:dGdB}
The partial derivative with respect to $\B$, holding the matrix value $\M$ fixed, is
\begin{equation}\label{eq:dGdB}
\frac{\partial G}{\partial\B} =-2\M\B\bm S.
\end{equation}
\end{proposition}
\begin{proof}
At a fixed reference point, suppress the element index and put
$\A=\B^T\M\B$, $\J=\A^{-1}$, and $D=\det\J$.
In this partial derivative only, $\dd\M=0$. 
Differentiating $\A\J=\I_m$ and applying Jacobi's determinant identity, respectively, we obtain
\begin{equation}\label{eq:reviewJacobi}
 \dd\J=-\J(\dd\A)\J,
 \qquad \dd D=D\tr(\J^{-1}\dd\J)
             =-D\tr(\J\dd\A).
\end{equation}
Here $\J$ is symmetric because $\A$ is symmetric positive definite.
Using the independent-variable partial derivatives defined in
\eqref{eq:H}--\eqref{eq:g}, the chain rule yields
\begin{align*}
 \dd G
 &=\HH:\dd\J+g\,\dd D\\
 &=-\tr\!\left(\HH^T\J(\dd\A)\J\right)
   -Dg\tr(\J\dd\A)\\
 &=-\bigl(\J\HH\J+Dg\J\bigr):\dd\A
 =-\bm S:\dd\A.
\end{align*}
The third equality uses the symmetry of $\HH$, $\J$, and the
Frobenius trace identity.  In particular, $\bm S$ is symmetric.
Next, differentiating $\A=\B^T\M\B$ gives
\[
 \dd\A=(\dd\B)^T\M\B+\B^T\M\dd\B.
\]
Taking the Frobenius product with $\bm S$ and using the symmetry
of $\M$ and $\bm S$ shows that the two terms contribute equally:
\[
 \bm S:\dd\A=2(\M\B\bm S):\dd\B.
\]
Consequently $\dd G=-2(\M\B\bm S):\dd\B$, which proves
\eqref{eq:dGdB}.  Notice that this derivation does not require
$\B$ to be square, and it applies to embedded curves and surfaces as
long as $\B$ has full column rank.
\end{proof}

Since
\[ \B_K =\sum_a\x_{K,a}(\nabla_{\xiB}N_a)^T, \]
the element contribution to the gradient with respect to local vertex $a$ is
\begin{equation}\label{eq:nodalEAforce}
\nabla_{\x_{K,a}}I_{\EA,K} =\int_{\Khat} \left(-2\M\B_K\bm S_K\right) \nabla_{\xiB}N_a\dd\xiB.
\end{equation}
This formula applies without change to $\B\in\R^{2\times2}$, $\R^{3\times2}$, and $\R^{3\times3}$.

For a spatial metric evaluated at moving quadrature points, we must also account for the variation of its argument. 
This gives the complete physical-coordinate derivative.

\begin{proposition}[Complete gradient for a fixed spatial metric field]
\label{prop:spatialgradient}
Assume $\M^{(n)}$ is $C^1$ in a neighborhood of the current element image,
with uniform positive-definiteness, and $\B_K$ has full column rank.
At every reference point define $\bm S_K$ by \eqref{eq:S} using
$\M=\M^{(n)}(F_K)$ and let
\begin{equation}\label{eq:metricresponse}
 \bm R_K:=\frac{\partial G}{\partial\M}
            =-\B_K\bm S_K\B_K^T\in\R^{d\times d},
 \qquad
 (\bm f_K)_k=\bm R_K:\partial_{x_k}\M^{(n)}(F_K),
 \quad k=1,\ldots,d.
\end{equation}
Then the complete elementwise physical-coordinate nodal gradient of
\eqref{eq:IEAspatial} is given by
\begin{equation}\label{eq:spatialforce}
 \nabla_{\x_{K,a}} I_{\EA,\rm sp,K}^{(n)}
 =\int_{\Khat}\left[
   -2\M\B_K\bm S_K\nabla_{\xiB}N_a
   +N_a\bm f_K\right]\dd\xiB.
\end{equation}
For the reference-frozen realization \eqref{eq:IEAstage}, only the
first term occurs.  The formula is valid for bulk and embedded elements
when the metric is an ambient $d\times d$ SPD field.
\end{proposition}
\begin{proof}
Fix a reference point and suppress the element index. Write
$\x=F_K(\xiB;X)$, $\B=D_{\xiB}F_K$, and
$\A=\B^T\M(\x)\B$. For arbitrary perturbations $\delta\x_a$ of
local vertices, the isoparametric definitions give
\[
 \delta\x=\sum_a N_a\,\delta\x_a,\qquad
 \delta\B=\sum_a\delta\x_a(\nabla_{\xiB}N_a)^T,
 \qquad
 \delta\M=\sum_{k=1}^d (\partial_{x_k}\M)\,\delta x_k.
\]
The last variation accounts for the change in the prescribed spatial
field evaluated at the moving physical point. It does not differentiate
the procedure by which that field was constructed at the beginning of
the stage. By Proposition~\ref{prop:dGdB}, for any symmetric variation
of $\A$ we have $\delta G=-\bm S:\delta\A$. The product rule gives
\[
\delta\A=(\delta\B)^T\M\B+\B^T\M\delta\B+
                 \B^T(\delta\M)\B.
\]
Since $\M$ and $\bm S$ are symmetric, cyclic invariance of the trace
implies
\begin{align*}
\bm S:\bigl((\delta\B)^T\M\B+\B^T\M\delta\B\bigr)
   &=2(\M\B\bm S):\delta\B,\\
\bm S:\bigl(\B^T(\delta\M)\B\bigr)
   &=(\B\bm S\B^T):\delta\M.
\end{align*}
It follows that
\[
\delta G=-2(\M\B\bm S):\delta\B
             -(\B\bm S\B^T):\delta\M.
\]
The coefficient of $\delta\M$ identifies the matrix derivative
$\bm R=-\B\bm S\B^T$. Substituting the nodal variations yields
\begin{align*}
\delta G
&=-2\sum_a (\M\B\bm S):
         \bigl[\delta\x_a(\nabla_{\xiB}N_a)^T\bigr]
   +\bm R:\sum_{k=1}^d(\partial_{x_k}\M)
             \sum_a N_a(\delta\x_a)_k\\
&=\sum_a\delta\x_a\cdot
  \left[-2\M\B\bm S\nabla_{\xiB}N_a
       +N_a\bm f\right],
\end{align*}
where $(\bm f)_k=\bm R:\partial_{x_k}\M$. Integration over the
fixed reference cell and identification of the coefficient of each
$\delta\x_a$ proves \eqref{eq:spatialforce}. For the reference-frozen
realization, $\delta\M=0$, so the second term vanishes. In either
case the metric field or its reference labels are kept fixed while
the mesh relaxes.
\end{proof}

A positive convex multilinear interpolation of SPD nodal matrices
remains SPD and is continuous, but its derivative may jump when a moving
quadrature point crosses a Cartesian interpolation-cell interface.
Away from those interfaces \eqref{eq:spatialforce} is the ordinary
gradient of the quadrature energy.  At an interface it gives a branchwise
or one-sided derivative, not necessarily a classical gradient.  Energy
dissipation still follows along any absolutely continuous trajectory for
which the energy chain rule holds almost everywhere and the implemented
velocity uses a compatible derivative almost everywhere. This is a
conditional statement, not an assertion of classical smooth ODE
existence or finite termination of sufficient-decrease backtracking at every interpolation knot.  To apply the
smooth versions of the theorems below without that qualification, we use
a $C^1$ SPD-preserving field representation, or restrict the stage to
regions in which the interpolated field is differentiable.  Regardless
of regularity, each accepted discrete trial is tested against the
actual stage energy and the full physical-element certificate.

When the prescribed metric field is $C^1$ on a neighborhood of the element images, 
both stage energies are continuously differentiable on the full-rank set. 
A continuously differentiable SPD-preserving interpolant is constructed next. 
In contrast, ordinary multilinear interpolation is only piecewise differentiable, as discussed above. 
When numerical quadrature is used, the gradient of the quadrature-defined energy 
must be obtained by differentiating the same quadrature sum.

We now construct a $C^1$ SPD-preserving representation to support the differentiability assumption in the spatial-field gradient.
The regularity of the direct spatial-field gradient can be obtained without
assuming that the application generating the metric supplies an analytic
formula. Consider a closed Cartesian sampling box $\overline U\subset\R^d$
and SPD tensors $\M_{\bm j}$ at its grid nodes $y_{\bm j}$. Define the
symmetric logarithmic data
\begin{equation}\label{eq:logcubicnodes}
 \bm L_{\bm j}=\log\M_{\bm j}.
\end{equation}
Interpolate the components of $\bm L_{\bm j}$ by tensor-product cubic
Hermite interpolation, with each nodal partial derivative shared by adjacent
cells. For equally spaced samples, one may use centered differences for
interior slopes and one-sided differences at the box boundary, with the
mixed nodal slopes formed by commuting these one-dimensional difference
operators. Denote the resulting interpolant by $\bm L_h(x)$ and set
\begin{equation}\label{eq:logcubicM}
 \M_h(x)=\exp(\bm L_h(x)), \qquad x\in\overline U.
\end{equation}
The following result gives the regularity needed for the smooth gradient-flow analysis.

\begin{proposition}[Regularity and positivity of the interpolated metric]
\label{prop:logcubicSPD}
For a finite Cartesian grid with strictly positive spacings and SPD nodal
tensors, the field \eqref{eq:logcubicM} is $C^1$ on $\overline U$, is SPD at
every point, and interpolates the nodal tensors. Its eigenvalues are bounded
above and away from zero on $\overline U$. On compact subsets of the
interior, its first derivatives are locally Lipschitz. If the field is
required outside $\overline U$, these claims apply there only after a
corresponding $C^1$ SPD extension has been specified.
\end{proposition}
\begin{proof}
Each matrix $\bm L_{\bm j}$ is real symmetric. Componentwise cubic Hermite
interpolation preserves symmetry and reproduces the nodal data. Because
neighboring cells use the same nodal values and first-derivative data,
the function and its first partial derivatives agree on all common cell
faces. In the tensor-product construction, the tangential data on a common
face are likewise identical since the one-dimensional difference and
interpolation operations commute. Thus $\bm L_h\in C^1(\overline U)$.
The matrix exponential maps any symmetric matrix to an SPD matrix, and
$\exp(\log\M_{\bm j})=\M_{\bm j}$. Continuity on the compact box gives
finite upper and strictly positive lower eigenvalue bounds. On each grid
cell, $\bm L_h$ is polynomial and has bounded second derivatives. The
shared continuous first derivatives imply their local Lipschitz
continuity across the finitely many cell interfaces. Smoothness of the
matrix exponential gives the corresponding property for $\M_h$.
\end{proof}

The derivatives required by the spatial-field force in
\eqref{eq:spatialforce} are evaluated by the Fr\'echet derivative
\begin{equation}\label{eq:logcubicderiv}
 \frac{\partial\M_h}{\partial x_k}
 =D\exp_{\bm L_h}\!\left[
  \frac{\partial\bm L_h}{\partial x_k}\right].
\end{equation}
Writing $\bm L_h=Q\diag(\lambda_1,\ldots,\lambda_d)Q^T$, this derivative is
\begin{equation}\label{eq:logcubicloewner}
 \frac{\partial\M_h}{\partial x_k}
 =Q\left(E\circ\left(Q^T\frac{\partial\bm L_h}{\partial x_k}Q\right)\right)Q^T,
 \qquad
 E_{ab}=\begin{cases}
 e^{\lambda_a},&\lambda_a=\lambda_b,\\[3pt]
 \dfrac{e^{\lambda_a}-e^{\lambda_b}}{\lambda_a-\lambda_b},
   &\lambda_a\ne\lambda_b.
 \end{cases}
\end{equation}
Here $\circ$ denotes the Hadamard (entrywise) matrix product, $(A\circ B)_{ij}=A_{ij}B_{ij}$.
The repeated-eigenvalue expression is the continuous divided-difference
limit, not an eigenvector derivative. A stable implementation uses a
small-gap expansion or an equivalent cancellation-free evaluation.

Cubic interpolation of log-metric data remains SPD after exponentiation but may overshoot the nodal eigenvalue range. 
Uniform eigenvalue bounds follow on the compact interpolation box, not from simple convexity of the cubic interpolant. 
All trial element images must remain within its domain or a specified smooth SPD extension.

\begin{remark}[Less expensive multilinear interpolation]
\label{rem:multilinearPractical}
An ordinary tensor-product multilinear interpolant is a convex combination
of the SPD nodal tensors and therefore remains SPD. It has lower evaluation
cost but, in general, only piecewise continuous first derivatives. Its
spatial gradient formula is valid wherever its derivatives exist, and
physical candidates can still be accepted using the actual computed energy
and whole-element validity checks. The classical differentiable-flow and
finite-backtracking arguments in this paper apply directly to the $C^1$
realization under their remaining hypotheses. They must not be asserted
unconditionally at multilinear interpolation knots. The two interpolants
agree at the sampled nodes but define different continuous metric fields
between them, so identical trajectories are not expected. The choice is an
implementation tradeoff, not a change to the intended purpose of the
prescribed metric.
\end{remark}

A one-dimensional calculation illustrates the different equidistribution 
behavior of the spatial-field and reference-frozen realizations.  
Let $m=d=1$, $\M=\rho(x)^2$
with $\rho>0$, and $h=\partial_\xi x>0$.  The density
\eqref{eq:G} simplifies to
\begin{equation}\label{eq:1dG}
  G=(1-\theta)(\rho(x)h)^{1-p}.
\end{equation}
If $\rho$ is attached to fixed reference-cell labels $\rho_i$, a
piecewise affine mesh with fixed endpoints minimizes a sum proportional
to $\sum_i\rho_i^{1-p}h_i^{1-p}$ under the constraint
$\sum_i h_i=|\Omega|$.  The Lagrange multiplier equation is
$\rho_i^{1-p}h_i^{-p}=\mathrm{constant}$, and therefore
\begin{equation}\label{eq:frozen1dresponse}
 h_i\propto\rho_i^{-(p-1)/p}.
\end{equation}
In general this differs from $h_i\propto\rho_i^{-1}$.

For the spatial-field energy, define the strictly increasing
cumulative metric coordinate
$y(x)=\int_{x_L}^{x}\rho(s)\dd s$.  Since
$\partial_\xi y=\rho(x)h$, the continuous energy becomes
$(1-\theta)\int(\partial_\xi y)^{1-p}\dd\xi$.
The function $s\mapsto s^{1-p}$ is strictly convex for $p>1$.
For fixed boundary values of $y$, Jensen's inequality implies that
the unique positive-derivative minimizer has constant
$\partial_\xi y$.  Thus the continuous spatial-field energy favors
\begin{equation}\label{eq:1dequidist}
 \rho(x)\partial_\xi x=\mathrm{constant}.
\end{equation}
This is the standard one-dimensional metric equidistribution condition, 
in contrast to the reference-frozen response \eqref{eq:frozen1dresponse}.

If the reference integral in \eqref{eq:IEA} is approximated by a
finite quadrature with fixed positive weights, the
implemented adaptation energy is the finite sum of its positive
integrands.  On a fixed metric stage, differentiating that same
quadrature sum yields the nodal force obtained by quadrature of
\eqref{eq:nodalEAforce}.  The energy-decay and barrier-preservation
proofs below apply to this discrete quadrature energy when the
complete gradient of the selected reference-frozen or spatial-field
realization is used, provided it is the energy tested for decrease.  

\subsection{Unified moving mesh PDE and energy decay}\label{sec:mmpde}

Following the MMPDE philosophy of \cite{HuangRenRussell1994,HuangRussell2011,HuangKamenski2018,ZhangHuang2026}, 
the negative variational gradient of the meshing energy defines the mesh motion.  
At the continuous level one may view the mesh as a time-dependent map $\x=\x(\xiB,t)$ and write schematically
\begin{equation}\label{eq:continuousMMPDEschematic}
\frac{\partial\x}{\partial t}
=-\frac{1}{\tau}\,\mathcal M_{\rm mob}(\x)\,
\frac{\delta I}{\delta\x},
\end{equation}
where $\mathcal M_{\rm mob}$ is a positive mobility operator and $t$ is a mesh or pseudo-time variable.  
After spatial discretization by the mesh vertices, \eqref{eq:continuousMMPDEschematic} becomes a finite-dimensional system of ODEs.  
Thus the nodal equations below are the spatially discrete MMPDE.

Let
\[
X=(\x_1^T,\ldots,\x_{N_v}^T)^T
\]
collect the mesh coordinates and let $I(X)$ be a differentiable discrete mesh energy.  
To accommodate bulk, surface, and boundary constraints in a unified formulation, 
we project the negative gradient of the meshing energy onto the admissible nodal-motion directions. 
The resulting projected nodal MMPDE is given by 
\begin{equation}\label{eq:mmpde}
\frac{\dd\x_i}{\dd t}
=-\frac{P_i}{\tau}
\PiB_i\nabla_{\x_i}I,
\qquad i=1,\ldots,N_v,
\end{equation}
where $\tau>0$ controls the pseudo-time scale, $P_i>0$ is a scalar mobility/scaling factor, 
and $\PiB_i$ is the admissible-motion projector from \eqref{eq:Pi}.  In block form,
\begin{equation}\label{eq:mmpdeBlock}
\dot X
=-\frac{1}{\tau}\mathcal D(X)\mathcal P(X)\nabla_X I(X),
\end{equation}
where $\mathcal D$ is block diagonal with positive factors $P_i$ and $\mathcal P$ is the block-diagonal projector.

Any positive $P_i$ preserves the gradient-flow sign property.  
To make the motion invariant under the scaling $\M\mapsto\chi\M$, 
a convenient intrinsic choice is
\begin{equation}\label{eq:PiScale}
P_i=
\det\bigl(\Tframe_i^T\M_i\Tframe_i\bigr)^{(p-1)/2}
\end{equation}
for $m<d$ and
\[
P_i=\det(\M_i)^{(p-1)/2}
\]
for $m=d$.  Under $\M\mapsto\chi\M$ with $\chi>0$, this scales as $\chi^{m(p-1)/2}$ 
and cancels the corresponding constant scaling of the energy gradient.

For stiff planar direct flows, the implementation also permits a
positive, stage-frozen Jacobi-type nodal balancing,
\begin{equation}\label{eq:jacobimobility}
 a_i=\sum_{K\ni i}\sum_q w_q G_{K,q}
       \|\nabla_xN_i\|_{K,q}^2,
 \qquad P_i=\frac{c}{a_i},\quad c>0,
\end{equation}
provided $0<a_i<\infty$. For embedded elements,
$\nabla_xN_i=\B_K(\B_K^T\B_K)^{-1}\nabla_{\xiB}N_i$
means the tangential physical gradient.  The normalization $c$ is
chosen independently of the nodal variation during a stage.  This
quantity is a positive stiffness estimate, not the exact diagonal
of the energy Hessian.  Replacing one positive scalar $P_i$ by another
changes the trajectory and its pseudo-time scaling, but not the
projected stationary equation $\PiB_i\nabla_{x_i}I=0$. The
energy-dissipation proof only requires $0<P_i<\infty$.

The projected negative-gradient structure gives an energy-dissipation identity whenever 
the stage energy and projectors satisfy the stated regularity assumptions.

\begin{theorem}[Energy decay]\label{thm:decay}
Assume $P_i>0$ and that $\PiB_i$ is symmetric positive semidefinite.  Along any smooth solution of \eqref{eq:mmpde},
\begin{equation}\label{eq:decay}
\frac{\dd I}{\dd t}
=-\sum_{i=1}^{N_v}\frac{P_i}{\tau}
\nabla_{\x_i}I^T\PiB_i\nabla_{\x_i}I
\le0.
\end{equation}
Hence $I(t)\le I(0)$.
\end{theorem}
\begin{proof}
Because $X(t)$ and $I$ are differentiable, the ordinary chain rule
in the finite-dimensional nodal-coordinate space gives
\[
 \frac{\dd}{\dd t}I(X(t))
 =\nabla_X I(X(t))^T\dot X(t)
 =\sum_{i=1}^{N_v}\nabla_{\x_i}I^T\dot\x_i.
\]
Substituting the nodal MMPDE \eqref{eq:mmpde} yields
\[
 \frac{\dd I}{\dd t}
 =-\frac1\tau\sum_{i=1}^{N_v}P_i
 \nabla_{\x_i}I^T\PiB_i\nabla_{\x_i}I.
\]
For a symmetric positive-semidefinite projector, every quadratic
form in this sum is nonnegative. For an orthogonal projector it equals
$\|\PiB_i\nabla_{\x_i}I\|_2^2$.  Since $P_i/\tau>0$,
\eqref{eq:decay} follows.  Integrating the inequality from $0$ to
$t$ proves $I(t)\le I(0)$.  In particular, when all projectors
are orthogonal, equality in the derivative holds precisely when
$\PiB_i\nabla_{\x_i}I=0$ at every mobile node.  No derivative of
$P_i$ or $\PiB_i$ enters this first-order chain-rule calculation.
\end{proof}

If the energy is only piecewise $C^1$ due to metric interpolation, the
identity may instead be integrated over an absolutely continuous
trajectory on which the appropriate chain rule and projected descent
identity hold almost everywhere.  Classical differentiability at every
time, unique ODE solvability at interpolation knots, 
and finite termination of backtracking do not follow from mere continuity of the metric.

Theorem~\ref{thm:decay} is the unified projector form of the separate bulk, surface, 
and curve decay arguments used in \cite{ZhangHuang2026}.  The identity is a stagewise result: 
it applies to the same differentiable energy whose complete gradient defines the velocity, 
and does not compare energies across a metric update, certificate refresh, or PDE remapping.  
It is the main dynamical ingredient in the certificate-preservation theorem of Section~\ref{sec:barrier}.

The MMPDE parameter $\tau$ and the numerical time step are distinct.  The parameter $\tau$ sets the response time of 
the mesh in \eqref{eq:mmpde}. Decreasing $\tau$ accelerates the continuous pseudo-time dynamics.  
A time integrator introduces an independent step $\Delta t_n=t_{n+1}-t_n$.  For example, forward Euler gives
\begin{equation}\label{eq:eulerMMPDEpreview}
\x_i^{n+1}
=\x_i^n
-\frac{\Delta t_n}{\tau}
P_i^n\PiB_i^n\nabla_{\x_i}I(X^n).
\end{equation}
The ratio $\Delta t_n/\tau$ therefore determines the displacement produced by one explicit step, 
whereas $\tau$ remains part of the underlying MMPDE model.  Section~\ref{sec:discrete} constructs $\Delta t_n$ 
so that the discrete step preserves the validity certificate and decreases the mesh energy.  
Higher-order time integrators are possible, provided every internal stage satisfies the same validity requirements,
 but temporal order is secondary when the main objective is a steady adapted mesh.

Besides dissipating the stage energy, the projected velocity preserves the admissible nodal constraint manifolds, 
as follows from the standard ODE invariance argument.

\begin{proposition}[Preservation of the geometric constraints]
\label{prop:constraintpreserve}
Suppose each constrained node belongs initially to a smooth embedded
constraint stratum $S_i\subset\R^d$ without boundary in its own
relative topology.  Assume the nodal vector field in
\eqref{eq:mmpde} is locally Lipschitz in a neighborhood of the
admissible nodal configurations and satisfies
$\dot\x_i\in T_{\x_i}S_i$ whenever $\x_i\in S_i$.
Then each such node stays on $S_i$ throughout the interval of
existence of the unique local ODE solution.  Fixed nodes remain fixed.
\end{proposition}
\begin{proof}
On a coordinate neighborhood of $S_i$, the restriction of the vector
field to its tangent bundle defines a locally Lipschitz ODE in local
manifold coordinates.  Starting from the prescribed point of $S_i$,
this equation has a unique local solution in $S_i$.  Viewed through
the embedding into $\R^d$, that solution also satisfies the original
ambient nodal ODE because the ambient vector field is tangent to
$S_i$.  Local uniqueness of the ambient initial-value problem makes
the two solutions coincide.  Repeating this argument in overlapping
charts proves invariance for the full interval on which the solution
exists.  A projector set to zero gives $\dot\x_i=0$ for a fixed node.
For nodes on geometric boundaries or corners, one uses the
corresponding smooth boundary or corner stratum as $S_i$. Projecting
merely onto the ambient surface tangent plane would not in general
keep a boundary node on that boundary.
\end{proof}

The geometric constraint statement concerns the nodal positions. 
It does not imply that a non-affine element is nonsingular between its vertices. 
We therefore turn to whole-element validity.

\section{Element validity and certified nonsingularity}\label{sec:validitymain}
The preceding energy and gradient formulas determine how the mesh moves.
We now establish geometric tests for accepting or certifying the resulting non-affine elements.

For an affine simplex, the Jacobian is constant.  Loss of element validity is therefore a whole-element event 
that can be related directly to element altitude and volume as in \cite{HuangKamenski2018,ZhangHuang2026}.  
For a non-affine tensor-product element, however, the Jacobian varies spatially.  
A small interior region can approach singularity while corner or quadrature samples remain benign.  
Consequently, a proof based only on elementwise metric volume or finitely sampled Jacobians is not sufficient for a trilinear hexahedron.

We introduce a finite certificate whose positivity implies regularity everywhere.
Let $R\subset\Khat$ be a closed reference subcell.  Choose a fixed oriented orthonormal frame
\[
\Tframe_{K,R}\in\R^{d\times m},
\qquad
\Tframe_{K,R}^T\Tframe_{K,R}=\I_m.
\]
Define the projected square Jacobian
\begin{equation}\label{eq:Cproj}
\bm C_{K,R}(\xiB)
=\Tframe_{K,R}^T\B_K(\xiB)
\in\R^{m\times m}
\end{equation}
and the projected orientation determinant
\begin{equation}\label{eq:qproj}
q_{K,R}(\xiB)=\det\bm C_{K,R}(\xiB).
\end{equation}

For a bulk domain, the local certificate frame is $\Tframe_{K,R}=\I_d$, and hence $q_{K,R}=\det\B_K$. 
Unlike the nodal tangent frame $\Tframe_i$ in \eqref{eq:Pi}, $\Tframe_{K,R}$ is fixed over the reference subcell.  
For a surface, $\Tframe_{K,R}$ consists of an oriented orthonormal basis of the selected projection plane. 
For a curve, it consists of a single unit tangent direction.

Since $\operatorname{rank}(\Tframe_{K,R}^T\B_K)\le\operatorname{rank}\B_K$, positivity of $q_{K,R}$ implies full column rank of $\B_K$. 
The projected condition is stronger than mere rank preservation, 
but this strength is useful because it gives a signed scalar polynomial for tensor-product elements.

\begin{theorem}[Finite projection certificate]\label{thm:finiteprojection}
Let $F_K\in C^1(\Khat;\R^d)$ satisfy
\[
\rank\B_K(\xiB)=m
\qquad\text{for all }\xiB\in\Khat.
\]
Then there exist finitely many closed tensor-product subcells $R_1,\ldots,R_L$ covering $\Khat$ 
and fixed oriented orthonormal frames $\Tframe_1,\ldots,\Tframe_L\in\R^{d\times m}$ such that
\begin{equation}\label{eq:positiveprojection}
\det(\Tframe_\ell^T\B_K(\xiB))>0,
\qquad \xiB\in R_\ell,
\quad \ell=1,\ldots,L.
\end{equation}
\end{theorem}
\begin{proof}
Fix $\xiB_0\in\Khat$.  Full column rank implies that the columns of
$\B_K(\xiB_0)$ span an $m$-dimensional subspace of $\R^d$.
Choose an orthonormal basis of this subspace, arranged as the
columns of $\Tframe_{\xiB_0}$.  Its orientation can be selected so
that the nonsingular square matrix
$\Tframe_{\xiB_0}^T\B_K(\xiB_0)$ has positive determinant.
With the frame held fixed, the scalar function
$\xiB\mapsto\det(\Tframe_{\xiB_0}^T\B_K(\xiB))$
is continuous.  It therefore remains strictly positive on some
relative open neighborhood $U_{\xiB_0}$ of $\xiB_0$ in $\Khat$.
These neighborhoods cover the compact reference cube.
Choose a finite subcover $U_1,\ldots,U_s$ with associated frames.
Then, by the Lebesgue number lemma, there is an $\ell>0$ such that every
subset of $\Khat$ of diameter below $\ell$ is contained in one
$U_j$.  Partition $\Khat$ into finitely many closed axis-aligned
cubes of diameter below $\ell$ and assign to each cube a frame
from a containing $U_j$.  The determinant is then positive on
that entire closed cube, including its boundary.
The frames need not agree across neighboring subcells, because the
certificate tests full rank separately on each closed subcell.
\end{proof}

For bulk elements the frame is fixed to $\I_d$, and initial positive orientation must be checked. 
For surfaces, the independently chosen local projection frames certify full rank 
but not global orientation or self-intersection freedom.

\subsection{Bernstein certification for bilinear quadrilaterals and trilinear hexahedra}\label{sec:bernstein}

On $[0,1]$, the degree-$r$ Bernstein basis is
\[
B_i^r(t)=\binom{r}{i}t^i(1-t)^{r-i},
\qquad i=0,\ldots,r.
\]
The tensor-product basis on $[0,1]^m$ is
\[
B_{\bm\alpha}^{\bm r}(\xiB)
=\prod_{j=1}^m B_{\alpha_j}^{r_j}(\xi_j).
\]
Every polynomial $q$ of tensor degree at most $\bm r=(r_1,\ldots,r_m)$ has a unique expansion
\begin{equation}\label{eq:bernexp}
q(\xiB)=\sum_{\bm\alpha}b_{\bm\alpha}
B_{\bm\alpha}^{\bm r}(\xiB).
\end{equation}
Since the Bernstein basis functions are nonnegative and form a partition of unity, we have
\begin{equation}\label{eq:bernbound}
\min_{\bm\alpha}b_{\bm\alpha}
\le q(\xiB)
\le \max_{\bm\alpha}b_{\bm\alpha},
\qquad \xiB\in[0,1]^m.
\end{equation}
If all the coefficients are positive, then $q(\xiB)\ge\min_{\bm\alpha}b_{\bm\alpha}>0$ on the entire closed reference cell. 
This convex-hull property underlies certified finite-element validity tests \cite{JohnenRemacleGeuzaine2013,JohnenWeillRemacle2017}.

If the lower bound is nonpositive while the true polynomial is positive, 
reference-space subdivision produces new Bernstein representations on smaller cells.  
The control values converge to the polynomial values as the cells shrink, 
so a strictly positive polynomial is eventually certified by sufficiently fine subdivision; 
this is the basis of adaptive validity algorithms in \cite{JohnenRemacleGeuzaine2013,JohnenWeillRemacle2017}.

Together with Theorem~\ref{thm:finiteprojection}, these bounds imply a finite positive-certificate result for polynomial isoparametric maps.

\begin{proposition}[Existence of a strictly positive Bernstein certificate]
\label{prop:certificateexists}
For a full-rank polynomial $Q_1$ map satisfying the hypotheses of
Theorem~\ref{thm:finiteprojection}, there is a finite partition of its
reference cube into tensor-product subcells with fixed projection frames
such that every Bernstein coefficient of each projected determinant is
strictly positive.  For an affine full-rank $P_1$ map on a reference
simplex, a single appropriately oriented projection frame and its
constant positive determinant provide the certificate without subdivision.
For bulk elements the identity frame is prescribed, and the corresponding
conclusion requires $\det\B_K>0$ everywhere.
\end{proposition}
\begin{proof}
First consider a $Q_1$ map.  By Theorem~\ref{thm:finiteprojection},
there is a finite tensor-product partition with one fixed frame on each
closed subcell $R$.  On $R$, the projected determinant $q$ is a
polynomial and, by continuity and compactness,
$\min_{\xiB\in R}q(\xiB)=\varepsilon_R>0$.
We show that its Bernstein coefficients become positive after a finite
further subdivision.

For a subcell $R'=\bm a+\operatorname{diag}(\bm h)[0,1]^m$
contained in $R$, express $q$ in local coordinates as
$\widehat q(\bm t)=q(\bm a+\operatorname{diag}(\bm h)\bm t)$.
Taylor expansion of the fixed polynomial about $\bm a$ is exact, and
its nonconstant power-basis coefficients are $O(\|\bm h\|)$ as
$\|\bm h\|\to0$, uniformly in $\bm a\in R$.
At the fixed tensor degree of $q$, the change from power-basis
coefficients to Bernstein coefficients is a fixed linear map.  It
maps the constant polynomial $q(\bm a)$ to the coefficient vector
whose entries all equal $q(\bm a)$.  Consequently every Bernstein
coefficient $b_\alpha(R')$ satisfies
\[
  b_\alpha(R')=q(\bm a)+O(\|\bm h\|),
\]
with a bound uniform over $\bm a\in R$ and the finitely many indices
$\alpha$.  Choose the subcell diameters small enough that the absolute
error is less than $\varepsilon_R/2$.  Every new coefficient is then
at least $\varepsilon_R/2>0$.  There are finitely many initial
subcells, so one finite common refinement suffices.

For a $P_1$ simplex, $\B_K$ is constant.  Choosing a frame that makes
$\det(\Tframe^T\B_K)>0$ proves positivity directly. Its degree-zero
Bernstein representation has one positive coefficient.  No
rectangular partition of a simplex is needed.  In the bulk case the
frame is fixed to $\I_d$, so the argument applies to positively
oriented elements and does not permit reversal of an inverted element.
\end{proof}

A nonpositive coefficient on an unsubdivided cell is inconclusive. Positivity of the polynomial can still be established on a sufficiently refined cover.

For a planar or projected surface bilinear map,
\[
F(\xi,\eta)=\bm a+\bm b\xi+\bm c\eta+\bm d\xi\eta.
\]
After projection to a fixed two-dimensional frame, we have
\[
\bm C(\xi,\eta)
=[\bm b_*+\bm d_*\eta,\,\bm c_*+\bm d_*\xi],
\]
where $\bm b_*=\Tframe^T\bm b$, etc.  Hence
\begin{align}
q(\xi,\eta)
&=\det(\bm b_*+\bm d_*\eta,
          \bm c_*+\bm d_*\xi)\notag\\
&=\det(\bm b_*,\bm c_*)
 +\xi\det(\bm b_*,\bm d_*)
 +\eta\det(\bm d_*,\bm c_*),
\label{eq:quadq}
\end{align}
with no $\xi\eta$ term.

The affine projected determinant in \eqref{eq:quadq} yields an exact test using only the four corners.

\begin{proposition}[Exact corner certificate for a $Q_1$ quadrilateral]\label{prop:quadcorners}
For a fixed projection frame, $q(\xi,\eta)>0$ on $[0,1]^2$ if and only if it is positive at the four corners.
\end{proposition}
\begin{proof}
Equation~\eqref{eq:quadq} has the form
$q(\xi,\eta)=a+b\xi+c\eta$.
It admits the exact bilinear convex-combination representation
\begin{align*}
 q(\xi,\eta)={}&(1-\xi)(1-\eta)q(0,0)
 +\xi(1-\eta)q(1,0)\\
 &+(1-\xi)\eta q(0,1)+\xi\eta q(1,1).
\end{align*}
All four weights are nonnegative and sum to one.  Thus strictly
positive corner values imply $q>0$ everywhere, including the edges.
The converse follows because the corners belong to the closed cell.
The same argument applies after an affine rescaling to a reference
subcell.
\end{proof}

The same result applies to a surface quadrilateral after projection to a fixed oriented local plane on a certificate subcell.

In contrast, the determinant of a trilinear hexahedron varies inside the reference cell, 
so a whole-element test is required. For $m=d=3$, $F_K$ is trilinear.  
Each column of $\B_K$ is bilinear in the other two variables, and
\begin{equation}\label{eq:hexdet}
q_K(\xi,\eta,\zeta)=\det\B_K(\xi,\eta,\zeta)
\end{equation}
is triquadratic.  Therefore it has the unique expansion
\begin{equation}\label{eq:hexbern}
q_K(\xi,\eta,\zeta)
=\sum_{i,j,k=0}^2
b_{ijk}
B_i^2(\xi)B_j^2(\eta)B_k^2(\zeta).
\end{equation}
There are 27 coefficients.  Johnen, Weill, and Remacle \cite{JohnenWeillRemacle2017} show that only 20 independent Jacobian evaluations 
are needed in an optimized implementation.  For a transparent reference implementation, however,
all 27 coefficients can be recovered from values on the tensor grid $\{0,1/2,1\}^3$.

In one dimension, if $c=(q(0),q(1/2),q(1))^T$ and $b=(b_0,b_1,b_2)^T$, then
\begin{equation}\label{eq:T2}
b=\bm T_2 c,
\qquad
\bm T_2=
\begin{bmatrix}
1&0&0\\
-1/2&2&-1/2\\
0&0&1
\end{bmatrix}.
\end{equation}
Thus, with a consistent lexicographic ordering,
\begin{equation}\label{eq:kronecker}
\bm b
=(\bm T_2\otimes\bm T_2\otimes\bm T_2)\bm c.
\end{equation}
A positive minimum coefficient gives the rigorous lower bound
\begin{equation}\label{eq:hexlower}
\det\B_K(\xi,\eta,\zeta)
\ge \min_{i,j,k}b_{ijk}>0
\qquad\text{throughout }\Khat.
\end{equation}
If the lower bound is inconclusive, de Casteljau subdivision (or an equivalent coefficient transformation) is 
applied recursively on the affected subcell. The boundary-based characterization of \cite{Zhang2026HexBoundaries} is 
another rigorous validity route for trilinear hexahedra. 
The method studied here retains Bernstein bounds for whole-element acceptance.
The preceding inequalities are exact mathematical statements.  For a
validated floating-point certificate, coefficient evaluation and
subdivision must use outward-rounded interval arithmetic or verified error
bounds sufficient to prove that each true coefficient remains strictly positive.  

\subsection{Certificate barrier and nonsingularity theorem}\label{sec:barrier}

We first define the fixed coefficients on which a certificate barrier acts.
A certificate stage consists of a finite collection $\mathcal C$ of certificate coefficients.  
Each $c\in\mathcal C$ is associated with an element, a reference subcell, a fixed projection frame, 
and one Bernstein coefficient of the projected determinant on that subcell.  
During a certificate stage, the subcell partition and projection frames are held fixed while the mesh vertices move.

Assume initially
\begin{equation}\label{eq:b0pos}
b_c^0:=b_c(\x(0))>0,
\qquad c\in\mathcal C.
\end{equation}
Define the normalized coefficient
\begin{equation}\label{eq:z}
z_c(\x)=\frac{b_c(\x)}{b_c^0}
\end{equation}
and the inverse-power validity barrier
\begin{equation}\label{eq:Ival}
I_{\val}(\x)
=\mu\sum_{c\in\mathcal C}
w_c z_c(\x)^{-r},
\qquad
\mu>0,
\quad r>0,
\quad w_c>0.
\end{equation}
The normalized coefficients are relative to the current
certificate-stage entry. The inverse-power barrier generally has a
nonzero gradient even at $z_c=1$, so its stationary mesh may depend
on the barrier parameters and the stage-refresh history. Under a
constant metric scaling $\M\mapsto\chi\M$, the adaptation energy
scales as $\chi^{-m(p-1)/2}$, while a fixed barrier weight does not.
A scale-consistent stage choice is
\begin{equation}\label{eq:muscale}
 \mu=\kappa\frac{I_{\EA}(X^{\rm entry})}
 {\sum_{c\in\mathcal C}w_c},\qquad\kappa>0,
\end{equation}
provided the initial adaptation energy is positive.
The full energy is
\begin{equation}\label{eq:Itotal}
\mathcal I_h(\x)
=I_{\EA}(\x)+I_{\val}(\x).
\end{equation}
The admissible configuration set for this stage is
\begin{equation}\label{eq:admissible}
\mathcal A_{\mathcal C}
=\{\x:\ b_c(\x)>0\text{ for every }c\in\mathcal C\}.
\end{equation}

To incorporate the validity barrier into the MMPDE, we next compute its nodal gradient.
At a fixed certificate evaluation point, let
\[
\bm C=\Tframe^T\B,
\qquad q=\det\bm C.
\]
For local node $a$,
\begin{equation}\label{eq:qgrad}
\nabla_{\x_a}q
=\Tframe\,\cof(\bm C)\,
\nabla_{\xiB}N_a.
\end{equation}
Indeed, $\delta q=\cof(\bm C):\delta\bm C$ and
$\delta\bm C=\Tframe^T\delta\x_a(\nabla N_a)^T$.

If a Bernstein coefficient is obtained from a fixed linear transformation of determinant values,
\begin{equation}\label{eq:blinear}
b_c=\sum_j Q_{cj}q(\xiB_j),
\end{equation}
then
\begin{equation}\label{eq:bgrad}
\nabla_{\x_a}b_c
=\sum_jQ_{cj}\,
\Tframe\,\cof(\bm C_j)\,
\nabla_{\xiB}N_a(\xiB_j).
\end{equation}
Consequently,
\begin{equation}\label{eq:Ivalgrad}
\nabla_{\x_a}I_{\val}
=-\mu r\sum_{c\in\mathcal C}
\frac{w_c}{b_c^0}
z_c^{-r-1}
\nabla_{\x_a}b_c.
\end{equation}
No finite-difference differentiation is required.

To avoid disturbing the minimizer when the validity margin is safely
above a prescribed threshold, an alternative is
\begin{equation}\label{eq:localbarrier}
 \phi_{r,z_0}(z)=
 \begin{cases}
 (z^{-1}-z_0^{-1})^r,&0<z<z_0,\\
 0,&z\ge z_0,
 \end{cases}
 \qquad 0<z_0\le 1,\quad r\in\{3,4,\ldots\},
\end{equation}
and $I_{\val}^{\rm loc}=\mu\sum_c w_c\phi_{r,z_0}(z_c)$.
This function is $C^2$ on $(0,\infty)$, is identically zero for
$z\ge z_0$, and diverges as $z\downarrow0$. In its active region,
\begin{equation}\label{eq:localbarriergrad}
 \nabla_X I_{\val}^{\rm loc}
 =-\mu r\sum_{c:z_c<z_0}\frac{w_c}{b_c^0}
 z_c^{-2}(z_c^{-1}-z_0^{-1})^{r-1}\nabla_Xb_c.
\end{equation}
The same energy-decay reasoning gives, on any smooth fixed stage
with nonnegative $I_{\EA}$ and initial total energy $E_0$,
\begin{equation}\label{eq:localbarrierbound}
 z_c(t)\ge
 \left[z_0^{-1}+(E_0/(\mu w_c))^{1/r}\right]^{-1}>0.
\end{equation}
Indeed, energy decay and nonnegativity imply
$\mu w_c\phi_{r,z_0}(z_c(t))\le E_0$; this is automatic for
$z_c\ge z_0$, while for $z_c<z_0$ rearrangement yields the bound.
Thus the alternative retains the fixed-stage positivity mechanism,
but its numerical influence and conditioning must be assessed rather
than assumed superior. Numerical comparisons of barrier choices under otherwise identical
validity tests are left to future work.

A fixed certificate stage has polynomial dependence on the nodal coordinates as follows.

\begin{lemma}[Smooth dependence of fixed-stage certificate coefficients]
\label{lem:coeffsmooth}
For fixed reference subcells, fixed projection frames, and fixed
Bernstein transformations, each coefficient $b_c(X)$ of a $P_1$ or
$Q_1$ element is a polynomial function of the local nodal coordinates
of total degree at most $m$.  In particular, $b_c$ is continuously
differentiable, the finite-coefficient admissible set
$\mathcal A_{\mathcal C}$ in \eqref{eq:admissible} is open in the
ambient nodal-coordinate space, and the barrier \eqref{eq:Ival} is
continuously differentiable on this set.
\end{lemma}
\begin{proof}
At each fixed reference point, the Jacobian $\B_K$ is linear in
nodal positions by \eqref{eq:Bdef}.  The frame is fixed, so the
entries of $\Tframe^T\B_K$ are also linear in those positions.
Its $m\times m$ determinant is a polynomial of total degree at most
$m$.  Each Bernstein coefficient is a fixed linear combination of
projected determinant evaluations, as in \eqref{eq:blinear}, and has
the same degree bound.  Polynomial functions are smooth. The inverse
image of $(0,\infty)$ under a continuous function is open, and the
finite intersection of these sets is open.  Finally, the inverse
powers in the barrier are smooth wherever all $b_c$ are positive.
\end{proof}

The preceding gradient and the energy-decay identity give a lower bound on the coefficients during each fixed stage.
We now use the projected MMPDE \eqref{eq:mmpde} with $I=\mathcal I_h$.

\begin{theorem}[Certificate preservation and quantitative lower bound]\label{thm:certificate}
Assume $0<\theta\le1/2$, so that $I_{\EA}\ge0$, and let
$\mathcal I_h=I_{\EA}+I_{\val}$ be the same fixed-stage
energy used to generate the direct motion.  Suppose its finite-dimensional
trajectory has nonincreasing energy on $[0,T)$. This holds for a smooth
projected negative gradient flow using the complete gradient, or for an
absolutely continuous trajectory satisfying the corresponding chain-rule
and dissipation identity almost everywhere.  If \eqref{eq:b0pos} holds,
then for every $c\in\mathcal C$ and $0\le t<T$,
\begin{equation}\label{eq:blower}
b_c(t)
\ge b_c^0
\left(\frac{\mu w_c}{\mathcal I_h(0)}\right)^{1/r}
>0.
\end{equation}
Hence no certificate coefficient can reach zero in finite time along such a trajectory.
\end{theorem}
\begin{proof}
We work within one certificate stage: the prescribed spatial metric field (or, for the alternative
realization, its reference samples), the certificate partition,
the projection frames, the positive normalizing values $b_c^0$,
and all weights are fixed.
Lemma~\ref{lem:coeffsmooth} ensures that the barrier is differentiable
where all coefficients remain positive.  Since the adaptation density
\eqref{eq:G} is nonnegative for $0<\theta\le1/2$, the total stage
energy satisfies $\mathcal I_h=I_{\EA}+I_{\val}\ge I_{\val}>0$.
The assumed fixed-stage energy dissipation (in the smooth case supplied
by Theorem~\ref{thm:decay} using the complete projected gradient) gives
$\mathcal I_h(X(t))\le\mathcal I_h(X(0))<\infty$.

Fix any $c\in\mathcal C$.  All terms of \eqref{eq:Ival} are positive,
so their sum is at least its $c$th term.  Consequently,
\[
 \mu w_c\left(\frac{b_c(t)}{b_c^0}\right)^{-r}
 \le I_{\val}(X(t))\le\mathcal I_h(X(t))
 \le\mathcal I_h(X(0)).
\]
All quantities in this inequality are positive, and $r>0$.
Rearranging therefore gives exactly \eqref{eq:blower}.
The bound holds for every coefficient individually and uniformly over
$0\le t<T$. It does not depend on the duration of the smooth stage.
Were one of the coefficients to tend to zero along the trajectory,
its barrier term would diverge, contradicting the finite initial
energy bound.  This excludes loss of the fixed certificate during the
trajectory's interval of existence.  It does not by itself establish
that the ODE has a global solution, nor compare energies across a
certificate or metric update.
\end{proof}

For every covered cell, the coefficient lower bound and the Bernstein convex-hull property imply $q_{K,R}(\xiB,t)>0$ 
on each certificate subcell. Thus all certified elements retain full column rank. In the bulk case, 
the fixed identity frame also preserves their chosen positive orientation. For planar $Q_1$ elements, 
the four corner values give this conclusion without subdivision. If the Jacobian norms are bounded, 
the argument yields a quantitative singular-value estimate.

\begin{corollary}[Quantitative local nonsingularity]\label{cor:sigma}
Suppose, in addition to the hypotheses of Theorem~\ref{thm:certificate},
that $\|\B_K(\xiB,t)\|_2\le L$ for all elements, reference points, and
$t\in[0,T)$, with $0<L<\infty$.  Set
\[
 \delta=\min_{c\in\mathcal C}
 b_c^0\left(\frac{\mu w_c}{\mathcal I_h(0)}\right)^{1/r}>0.
\]
Then for every certified element and reference point,
\begin{equation}\label{eq:sigmaLower}
 \sigma_{\min}(\B_K(\xiB,t))\ge \frac{\delta}{L^{m-1}}>0.
\end{equation}
In particular, when the metric eigenvalues are bounded below by
$\underline\rho$, the minimum singular value of
$\M^{1/2}\B_K$ is at least
$\sqrt{\underline\rho}\,\delta/L^{m-1}$.
\end{corollary}
\begin{proof}
Let $\bm C=\Tframe^T\B_K$ on a certified subcell.  By the convex-hull
bound and Theorem~\ref{thm:certificate}, $\det\bm C\ge\delta$.
Since $\|\bm C\|_2\le\|\B_K\|_2\le L$, the product formula for
singular values gives
\[
 \sigma_{\min}(\bm C)
 \ge\frac{\det\bm C}{L^{m-1}}
 \ge\frac{\delta}{L^{m-1}}.
\]
Orthogonal projection cannot increase norm.  Thus, for all $\bm v$,
$\|\B_K\bm v\|_2\ge\|\bm C\bm v\|_2$ and consequently
$\sigma_{\min}(\B_K)\ge\sigma_{\min}(\bm C)$.
The metric bound gives the final assertion.
\end{proof}

The preceding coefficient bound excludes approach to the certificate boundary, 
but existence for all pseudo-times requires additional assumptions. 
To conclude global-in-time existence of the finite-dimensional ODE, 
one also needs the nodal configuration to remain in a bounded set and 
the remaining coefficients of the vector field to stay regular.  
These properties do not follow merely from positivity of the certificate.
A global-in-time conclusion holds, for example, if the trajectory stays in a
compact subset of the permitted nodal-configuration manifold, the metric
samples and projectors remain sufficiently smooth with bounded derivatives,
and the fixed finite certificate stage is used without discontinuous updates.
Indeed, the coefficient bound then keeps the trajectory inside a compact
subset of the admissible region. The ODE right-hand side remains locally
Lipschitz there, so the standard continuation theorem applies.  Without
these additional hypotheses, Theorem~\ref{thm:certificate} proves
preservation only on the interval of existence of the smooth trajectory.

In a moving surface patch, a fixed projection frame may become restrictive after rotation. 
Similarly, a coarse Bernstein lower bound may become small even while the projected determinant remains positive. 
We therefore allow certificate refreshes as follows.
\begin{enumerate}[label=(\alph*)]
\item stop the current certificate stage at a valid mesh;
\item recompute local projection frames and/or subdivide reference cells;
\item certify strictly positive new coefficients;
\item set the current coefficients as the new $b_c^0$ and restart the barrier flow.
\end{enumerate}
Each fixed stage is covered by Theorem~\ref{thm:certificate}.  If
refresh times do not accumulate within the time interval of interest,
the endpoint of each certified stage is a valid initial state for the
next. An induction over stages then preserves regularity at all such
times.  Energy monotonicity remains stagewise, since refreshing
changes the barrier and possibly the metric.  For infinitely many
refreshes accumulating at a finite time, local stagewise positivity
alone is insufficient to certify the limiting mesh, and a uniform
lower bound or an independent continuation argument is additionally
required.  This distinction is particularly important on moving
surfaces, where projection frames may need repeated updates.

The explicit barrier bound can be conservative when the number of coefficients is large. 
Its dependence on coefficient weights and barrier parameters should not be interpreted as 
a measured numerical margin or as a substitute for trial validity checks.

\subsection{Validity-constrained discrete integration of the MMPDE}\label{sec:discrete}

Theorem~\ref{thm:Q1exactcoercivity} prevents local $Q_1$ rank loss on bounded sublevels of 
the exactly integrated energy, and Theorem~\ref{thm:certificate} provides a separate coefficient bound 
for a barrier-augmented stage. Neither conclusion automatically holds for the finite-quadrature time integrator, 
and an explicit step could cross the certificate boundary.  We therefore discretize the nodal MMPDE 
with a validity-preserving, energy-stable step-selection procedure.

At pseudo-time $t_n$, define the MMPDE velocity
\begin{equation}\label{eq:mmpdevelocity}
\bm v_i^n
=-\frac{P_i^n}{\tau}
\PiB_i^n\nabla_{\x_i}\mathcal I_h(X^n).
\end{equation}
The default time integrator is forward Euler,
\begin{equation}\label{eq:explicitEulerCertified}
\x_i^{n+1}
=\x_i^n+\Delta t_n\bm v_i^n.
\end{equation}
Thus the trial path is
\begin{equation}\label{eq:line}
\x_i(\Delta t)
=\x_i^n+\Delta t\,\bm v_i^n.
\end{equation}
For unconstrained interior nodes this is a time discretization of the MMPDE.
For motion on a curved surface, the Euler trial must additionally be retracted
to the surface, and the actual retracted candidate must be checked for validity
and energy decrease.

Because $\B_K$ is linear in the nodal coordinates,
\[
\bm C(\Delta t)=\bm C_0+\Delta t\,\bm C_1.
\]
Therefore
\begin{equation}\label{eq:qalpha}
q(\Delta t)
=\det(\bm C_0+\Delta t\,\bm C_1)
\end{equation}
is a polynomial of degree at most $m$ in the trial time step.  A Bernstein certificate coefficient is 
a fixed linear combination of such determinants, so
\begin{equation}\label{eq:balpha}
b_c(\Delta t)
=a_{c,0}+a_{c,1}\Delta t+\cdots+a_{c,m}(\Delta t)^m.
\end{equation}
Thus the validity-limiting step polynomial is at most linear for curves, quadratic for quadrilateral/surface certificates, 
and cubic for hexahedral certificates.

Define
\begin{equation}\label{eq:alphac}
\Delta t_c^*
=\inf\{\Delta t>0:\ b_c(\Delta t)=0\},
\end{equation}
with $\Delta t_c^*=+\infty$ if no positive root exists.  Let
\begin{equation}\label{eq:alphaval}
\Delta t_{\val}
=\eta\min_{c\in\mathcal C}\Delta t_c^*,
\qquad 0<\eta<1.
\end{equation}
If all $\Delta t_c^*=+\infty$, set $\Delta t_{\val}$ to a prescribed maximum trial step.  
Coefficient polynomials and their first positive roots are valid for a straight-line nodal trial, 
not for a nonlinear surface retraction.

The first positive root of each coefficient polynomial provides a conservative step cap that preserves the current certificate.

\begin{proposition}[Discrete certificate preservation]\label{prop:discrete}
For the unretracted Euler trial path, if $b_c(0)>0$ for every $c\in\mathcal C$ and the step satisfies
\[
0\le\Delta t_n\le\Delta t_{\val},
\]
then every certificate coefficient remains positive at the trial state.
\end{proposition}
\begin{proof}
Fix $c\in\mathcal C$.  Along the unretracted nodal line, $b_c(h)$ is
continuous (polynomial), and $b_c(0)>0$.  If it became
nonpositive at some $h>0$, the intermediate value theorem would
produce a zero in $(0,h]$.  Hence $b_c(h)>0$ for every
$0\le h<\Delta t_c^*$, with the same conclusion for all finite
$h$ when $\Delta t_c^*=+\infty$.
If any positive root exists among the finitely many certificate
polynomials, $0<\eta<1$ implies
$\Delta t_{\val}<\Delta t_c^*$ for every coefficient with a
finite positive root.  Coefficients without positive roots remain
positive for the whole positive ray.  

If there are no positive roots, the separately prescribed maximum trial step is admissible.
Therefore all coefficients remain strictly positive along the chosen
step.  Since the unchanged reference partition covers every cell,
the Bernstein bounds imply local element regularity at the trial
configuration.  A nonlinear retraction changes these polynomial
paths and requires an independent coefficient check at the actual
retracted state.
\end{proof}

After imposing the validity cap, we reduce the numerical time step until the following Armijo sufficient-decrease condition holds:
\begin{equation}\label{eq:armijo}
\mathcal I_h(X^n+\Delta t_n V^n)
\le
\mathcal I_h(X^n)
+c_A\Delta t_n\nabla\mathcal I_h(X^n)^T V^n,
\qquad 0<c_A<1,
\end{equation}
holds, where $V^n$ denotes the assembled nodal MMPDE velocity.  Since $\nabla\mathcal I_h^T V^n\le0$, 
this is an energy-stable adaptive time-step condition.  It should be interpreted as a numerical stabilization 
of the explicit MMPDE integration, rather than as changing the underlying method into a generic line-search optimizer.  
The accepted bulk update is then \eqref{eq:explicitEulerCertified} with
$t_{n+1}=t_n+\Delta t_n$. A retracted surface update is accepted only after a
fresh certificate and an energy test at the retracted nodal coordinates. The
polynomial step bound alone does not certify a nonlinear retraction.

Under differentiability assumptions, a sufficiently small admissible step also satisfies the energy-decrease test.

\begin{proposition}[Existence of an accepted sufficiently small step]
\label{prop:armijotermination}
Let the fixed-stage energy be continuously differentiable in a neighborhood
of the current valid configuration $X^n$, and let all its certificate
coefficients be strictly positive there.  Suppose $P_i>0$, $\tau>0$, and
the projectors are orthogonal.  For a surface update, suppose additionally
that the retraction $\mathcal R_{X^n}$ is $C^1$ near the origin, satisfies
$\mathcal R_{X^n}(0)=X^n$, and has
$D\mathcal R_{X^n}(0)V=V$ for every admissible tangent velocity $V$.
If the projected MMPDE velocity $V^n$ is nonzero, then a sufficiently small
positive trial step has positive certificate coefficients and satisfies the
Armijo inequality, with the energy evaluated at the actual retracted state
when applicable.  Thus geometric backtracking with a reduction factor in
$(0,1)$ terminates in exact arithmetic.
\end{proposition}
\begin{proof}
Write $g=\nabla\mathcal I_h(X^n)$ and $V=V^n$.  By the nodal
projected-gradient formula and the symmetry and idempotence of every
orthogonal projector,
\begin{equation}\label{eq:reviewSlope}
 d:=g^TV=-\frac1\tau\sum_{i=1}^{N_v}P_i
 \|\PiB_i^n\nabla_{\x_i}\mathcal I_h(X^n)\|_2^2.
\end{equation}
If $V\ne0$ then at least one projected gradient is nonzero, so
$d<0$.  For a bulk straight-line trial put $X(h)=X^n+hV$.
For a surface trial put $X(h)=\mathcal R_{X^n}(hV)$.
The asserted $C^1$ retraction and its first-order consistency yield
$X(0)=X^n$ and $X'(0)=V$ in either case.  Differentiability of the
stage energy therefore gives
\[
 \mathcal I_h(X(h))
 =\mathcal I_h(X^n)+h d+o(h)
 \quad\text{as }h\downarrow0.
\]
Since $d<0$ and $0<c_A<1$, choose $h_1>0$ so that
$|o(h)|\le(1-c_A)(-d)h$ for $0<h\le h_1$.  Then
\[
 \mathcal I_h(X(h))
 \le\mathcal I_h(X^n)+c_A h d,
\]
which is the Armijo condition with the energy evaluated on the
actual candidate path.
Moreover, Lemma~\ref{lem:coeffsmooth} and $b_c(X^n)>0$ give,
for each of the finitely many coefficients, an $h_c>0$ such that
$b_c(X(h))>0$ for $0\le h\le h_c$.  Take $h_2=\min_{c\in\mathcal C}h_c>0$.
Then all $0<h\le\min(h_1,h_2)$ pass both tests. In particular, a trial
sequence $h_0,\beta h_0,\beta^2h_0,\ldots$, with $0<\beta<1$,
reaches such a step after finitely many reductions.  

The argument uses a differentiable stage energy.  At a knot of the
piecewise-multilinear metric interpolant, the selected branch gradient
need not be a descent direction, so finite backtracking is not guaranteed
by this proposition.  An actual-energy reduction test can still reject
such proposals, or a differentiable SPD interpolation can restore its
hypotheses.  Exact energy and coefficient tests are assumed here.
validated floating-point acceptance needs the error enclosures
described above.
\end{proof}

If $V^n=0$, the scheme has reached a projected stationary state and no
positive step is required.  For a nonlinear retraction, the certificate
and energy checks must be performed on the actual retracted candidate,
not merely on its unretracted linear predictor.

The step-root calculation must be certified in finite precision: an approximate root alone does not provide a guarantee. 
Validated interval bounds or a positivity test on a one-dimensional Bernstein subdivision of the candidate step interval may be used.

Hard whole-element acceptance preserves local validity independently of the barrier. 
The barrier provides a continuous fixed-stage coefficient margin, while Theorem~\ref{thm:Q1exactcoercivity} concerns 
exact-energy sublevel sets. These guarantees are distinct, and none establishes global injectivity.

\subsection{Local validity versus global injectivity}\label{sec:global}

The barrier theory proves local element regularity. This is necessary but conceptually distinct from global one-to-one-ness 
of the entire mesh map.  Positive element-interior Jacobians do not exclude overlap between images of nonadjacent regions of a mesh.  
This distinction is classical in nonlinear elasticity and global inversion theory; see, for example, \cite{CiarletNecas1987,Traver2026}.

For a conforming bulk mesh, let
\[
F_h:\overline{\widehat\Omega}\to\R^d
\]
be the global piecewise-isoparametric map obtained by gluing the element maps.
A sufficient global-injectivity argument additionally requires a local homeomorphism across element interfaces, 
a one-to-one boundary map, and suitable topological-degree assumptions. 
These properties do not follow from positive Jacobians within individual cells and are not certified by the present method.

For a mesh constrained to a fixed embedded surface, local regularity plus preservation of the surface topology 
and degree can similarly be used to rule out local folds.  If the surface itself is evolving freely in $\R^3$, 
however, distant surface patches may collide without any element becoming singular.  
Such self-contact requires a nonlocal collision/contact barrier and is outside the scope of the local certificate theory.

\section{Prescribed metrics}\label{sec:adaptationmain}
Having established the movement and validity mechanisms, we describe how the prescribed metrics are constructed. 
The discrete procedure and numerical examples are presented in Section~\ref{sec:numerics}.

\subsection{Metric construction}\label{sec:adaptive}
The method accepts any prescribed, sufficiently regular SPD metric.
For solution-dependent adaptation, a stage consists of recovering a field
from the current solution, constructing a spatial metric, and relaxing the
physical mesh while this metric field is fixed. The metric may be updated
after an accepted stage. Its dependence on the preceding PDE solve is not
differentiated within the stage. The spatial derivative at the moving
quadrature locations is retained in \eqref{eq:spatialforce}.

For an $m$-dimensional surface in $\R^d$, an intrinsic tangent metric
$\M_T$ can be extended to an ambient SPD field by
\begin{equation}\label{eq:lift}
\M=\Tframe\M_T\Tframe^T+\mu_n(\I_d-\Tframe\Tframe^T),
\qquad \mu_n>0,
\end{equation}
where $\Tframe$ is an orthonormal tangent frame. The ambient metric
must be defined in a neighborhood of the surface to evaluate its spatial
derivatives on curved isoparametric patches. Since such patches need not
lie on the exact surface between nodes, their pullback metric can depend
on the normal extension $\mu_n$.

\subsection{Nodal Hessian recovery}\label{sec:hessianrecovery}
For Hessian-based monitors, we use one nodal quadratic least-squares
recovery for both element families. This avoids choosing a diagonal
subdivision of $Q_1$ cells solely for derivative recovery. Let
$\mathcal P_i$ denote neighboring nodes of $\x_i$ and
$\widehat{\bm r}_{ij}=S_i^{-1}(\x_j-\x_i)$, with a nonsingular
local scaling matrix $S_i$. The coefficients are determined by
\begin{equation}\label{eq:patchlsq}
\min_{c_i,\widehat{\bm g}_i,\widehat H_i=\widehat H_i^T}
\sum_{j\in\mathcal P_i}w_{ij}
\left[u_j-c_i-\widehat{\bm g}_i^T\widehat{\bm r}_{ij}
-\tfrac12\widehat{\bm r}_{ij}^T\widehat H_i\widehat{\bm r}_{ij}\right]^2,
\qquad w_{ij}>0.
\end{equation}
The recovered physical Hessian is
\begin{equation}\label{eq:hessianunscale}
H_i=S_i^{-T}\widehat H_i S_i^{-1}.
\end{equation}
Full column rank of the weighted quadratic design matrix ensures
uniqueness and exact reproduction of quadratic polynomials in the fitted
coordinates. Patches are enlarged when necessary to satisfy this rank
condition. For surface fields, the fit is performed in the local
$m$-dimensional tangent coordinates before the metric is lifted using
\eqref{eq:lift}.

\subsection{Metrics used in the examples}\label{sec:metrics_choices}
We illustrate the framework with accuracy-oriented and diffusion-aligned metric tensors. 
The scalar-monitor examples use the SPD tensor
\begin{equation}\label{eq:Madaptimpl}
\M_{\rm adapt}=\|\bm R_H\|_F^{1/m}\det(\bm R_H)^{-1/(2m)}\bm R_H,
\qquad \bm R_H=\I_m+|H|/\varepsilon_{\rm adapt},
\end{equation}
where $|H|$ is the spectral absolute value of the recovered Hessian
and $\varepsilon_{\rm adapt}>0$. The regularization and scaling in
\eqref{eq:Madaptimpl} specify the monitor used here. 

For diffusion-aligned movement, let $\bm D_K$ denote a representative SPD
diffusion tensor. We use the DMP-oriented metric of Li and Huang
\cite{LiHuang2010DMP}
\begin{equation}\label{eq:MDMPimpl}
\M_{{\rm DMP},K}=\bm D_K^{-1},
\end{equation}
and the combined diffusion-alignment/adaptation metric
\begin{equation}\label{eq:MDMPadaptimpl}
\M_{{\rm DMP+adapt},K}
=\left(1+\frac{\mathcal B_K}{\alpha_h}\right)^{\frac{2}{d+2}}
\det(\bm D_K)^{1/d}\bm D_K^{-1},
\end{equation}
where
\begin{align}
\mathcal B_K&=\det(\bm D_K)^{-1/d}\|\bm D_K^{-1}\|_2
\frac{1}{|K|}\int_K\|\bm D_K|H(u)|\|_2^2\dd x,
\label{eq:BKdmp}\\
\alpha_h&=\left[\frac{1}{|\Omega|}
\sum_{K\in\T_h}|K|\mathcal B_K^{d/(d+2)}\right]^{(d+2)/d}.
\label{eq:alphahdmp}
\end{align}
Here $|H(u)|$ is the spectral absolute value. The anisotropic
non-obtuse-angle theory underlying these metrics applies to simplicial
$P_1$ finite elements. Using the tensor as a monitor for $Q_1$ meshes
does not imply a $Q_1$ discrete maximum principle. In the examples,
metric construction may include nodal smoothing, complexity scaling,
and eigenvalue bounds. These operations alter the prescribed metric and
are distinct from the MMPDE itself.

\section{Numerical results}\label{sec:numerics}
We present numerical examples demonstrating metric-driven mesh redistribution for $P_1$ and $Q_1$ elements 
in one, two, and three dimensions, including a curved surface. 

\phantomsection\label{sec:implementation}
The direct method evaluates the isoparametric adaptation energy and
its complete spatial-field gradient, projects the velocity onto the
admissible geometric constraints, and uses backtracking steps that
decrease the discrete energy while passing the whole-element validity
test of Section~\ref{sec:discrete}. For curved surfaces, retraction
and certification are performed on the actual trial geometry. The
Bernstein orientation test is distinct from the quadrature used to
integrate the adaptation energy. After a metric refresh, the next
stage has a new energy reference.

For comparison, the indirect computational-coordinate MMPDE
\cite{HuangMMPDELab2019} varies computational coordinates with the
physical integration weights and metric samples fixed, and then
reconstructs physical nodes. It is a different gradient flow: decrease
of its computational-coordinate energy does not imply the direct
barrier theorem or validity of the reconstructed $Q_1$ mesh. 

The quantitative planar, graph-surface, and bulk experiments below
use the direct physical-coordinate method, with the metric fixed
during each relaxation stage and refreshed between stages. The
one-dimensional and two-dimensional figure pairs are 
qualitative illustrations obtained with different numerical settings.
They are not used to establish properties of the direct gradient flow
or to support the quantitative planar comparisons in
Table~\ref{tab:planarquant}.
Unless stated otherwise, $p=3/2$, $\theta=1/3$, and the spatial SPD metric is 
represented by a $C^1$ log-Hermite interpolant. The localized
barrier uses $\kappa=10^{-3}$, $r=3$, and $z_0=0.5$, together with
independent whole-element trial validity tests. Stagewise decrease
refers to the quadrature-discretized energy. Successive stages need
not have the same objective. Saved meshes were also checked for local
orientation using exact arithmetic on their stored coordinates.
These checks concern saved elements, not every intermediate
floating-point trial or global injectivity.

The alignment measure $Q_{\rm ali}\ge1$ has its ideal value at one.
To quantify equidistribution, let $\overline V_{\M}=N^{-1}\sum_{K\in\T_h}|K|_{\M}$,
where $N$ is the number of elements and $|K|_{\M}$ is defined in
\eqref{eq:metricvolume}. Then we use the coefficient of variation
\begin{equation}\label{eq:CVM}
\mathrm{CV}_M=\frac{1}{\overline V_{\M}}
\left[\frac{1}{N}\sum_{K\in\T_h}
\bigl(|K|_{\M}-\overline V_{\M}\bigr)^2\right]^{1/2}
\end{equation}
as an equidistribution measure. Its ideal value is zero. The numerical experiments have finite step
budgets and are not asserted to attain stationary or exactly
$\M$-uniform meshes.

\subsection{Representative one-dimensional adaptations}\label{sec:1Dmesh}
\begin{figure}[!htbp]\centering
\begin{subfigure}[t]{.485\linewidth}\centering
\includegraphics[width=\linewidth]{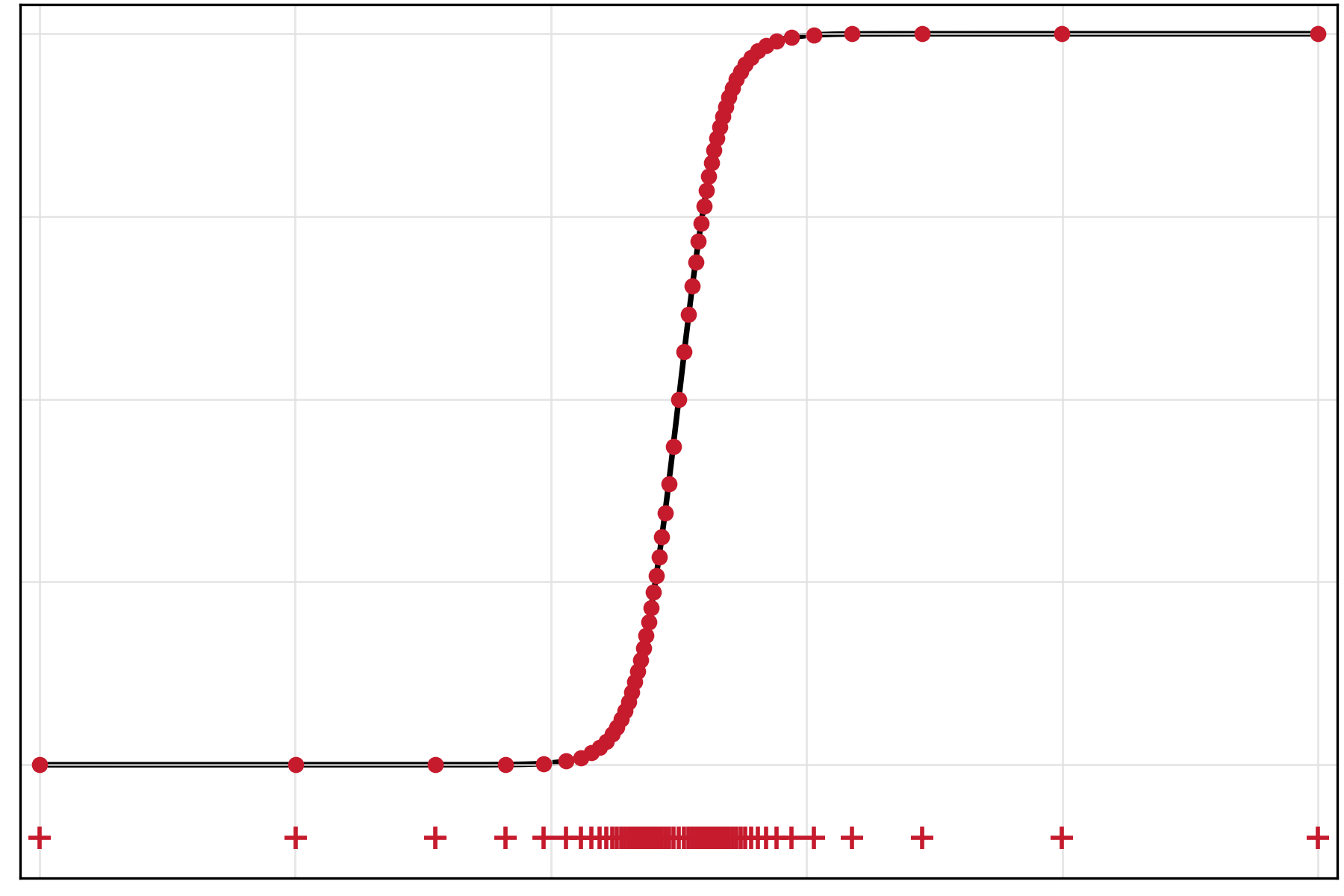}
\caption{$P_1$ intervals.}\label{fig:legacy1dP1}\end{subfigure}\hfill
\begin{subfigure}[t]{.485\linewidth}\centering
\includegraphics[width=\linewidth]{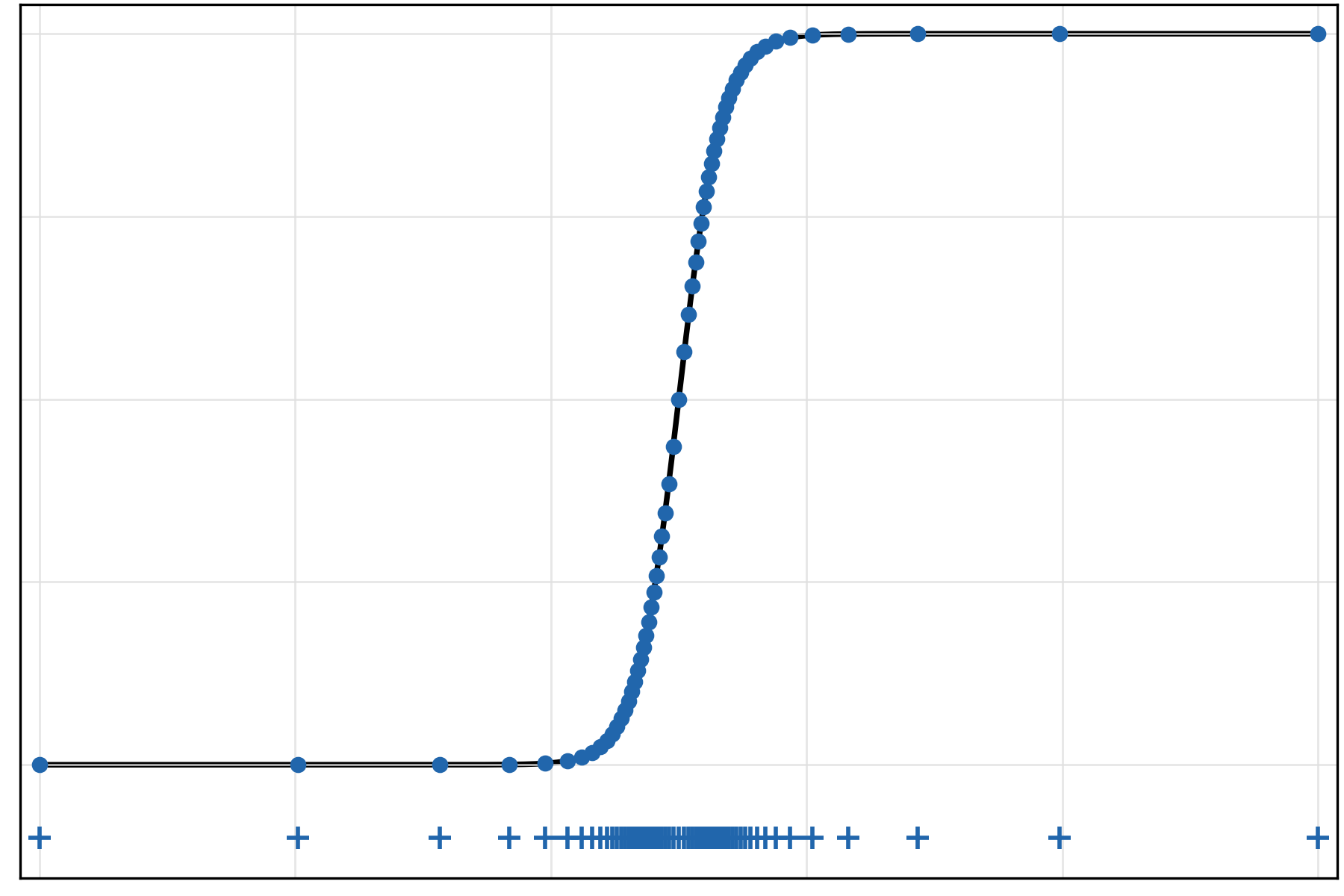}
\caption{$Q_1$ intervals.}\label{fig:legacy1dQ1}\end{subfigure}
\caption{One-dimensional Hessian-based adaptation using $\M_{\rm adapt}$ for solution $u(x)=\tanh[30(x-\tfrac12)]$. 
(a) $P_1$ mesh. (b) $Q_1$ mesh. The nodal positions become concentrated around the steep solution transition, 
as shown by the coordinate markers below the curves. The $P_1$ and $Q_1$ element maps are identical in one dimension.}
\label{fig:legacy1d}\end{figure}

We begin with a Hessian-based moving-mesh example to illustrate redistribution in one dimension. 
Figure \ref{fig:legacy1d} presents the resulting $P_1$ and $Q_1$ meshes under the same framework using $\M_{\rm adapt}$. 
The prescribed field is $u(x)=\tanh[30(x-\tfrac12)]$.
In one dimension, $P_1$ and $Q_1$ maps both describe affine intervals. The
figure illustrates the resulting concentration of mesh nodes near the
transition.

\FloatBarrier

\subsection{Prescribed scalar monitors on planar meshes}\label{sec:example_protocol}
We next demonstrate two-dimensional mesh redistribution using the three
smooth prescribed monitors from \cite{HuangMMPDELab2019} on $[0,1]^2$:
\begin{align}
 u_{\rm front}(x,y)&=\tanh[-30(y-1/2-\tfrac14\sin(2\pi x))],\notag\\
 u_{\rm cross}(x,y)&=\tanh(60y)-\tanh[60(x-y-1/2)],\notag\\
 u_{\rm rings}(x,y)&=\sum_{(a,b)\in S}\tanh\{30[(4x-2-a)^2+(4y-2-b)^2-1/8]\},\label{eq:analyticmonitors}\\
 S&=\{(0,0),(\tfrac12,\tfrac12),(\tfrac12,-\tfrac12),
       (-\tfrac12,\tfrac12),(-\tfrac12,-\tfrac12)\}.\notag
\end{align}

\begin{figure}[!htbp]\centering
\begin{subfigure}[t]{.485\linewidth}\centering
\includegraphics[width=\linewidth]{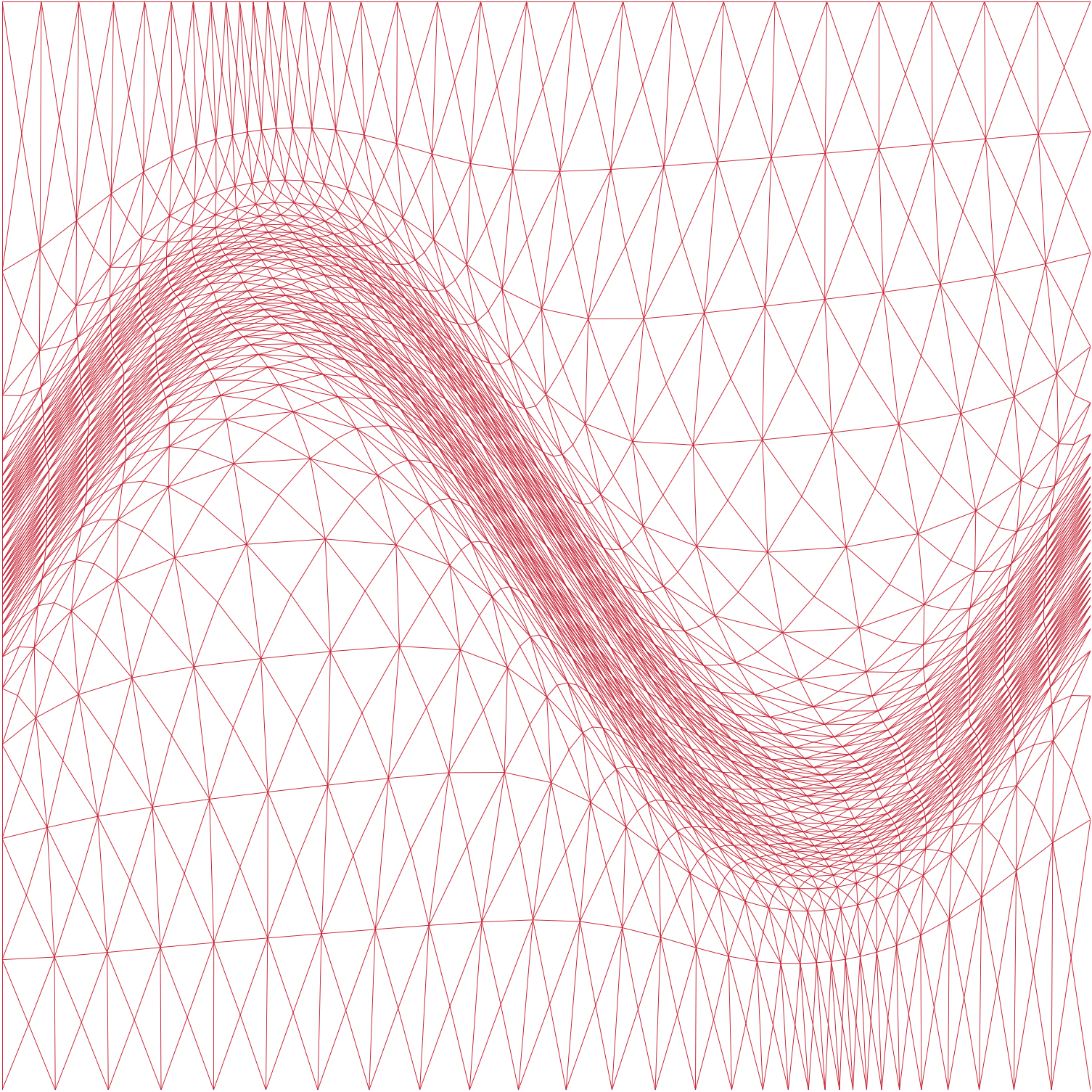}
\caption{$P_1$ triangles.}\label{fig:legacy2d3P1}\end{subfigure}\hfill
\begin{subfigure}[t]{.485\linewidth}\centering
\includegraphics[width=\linewidth]{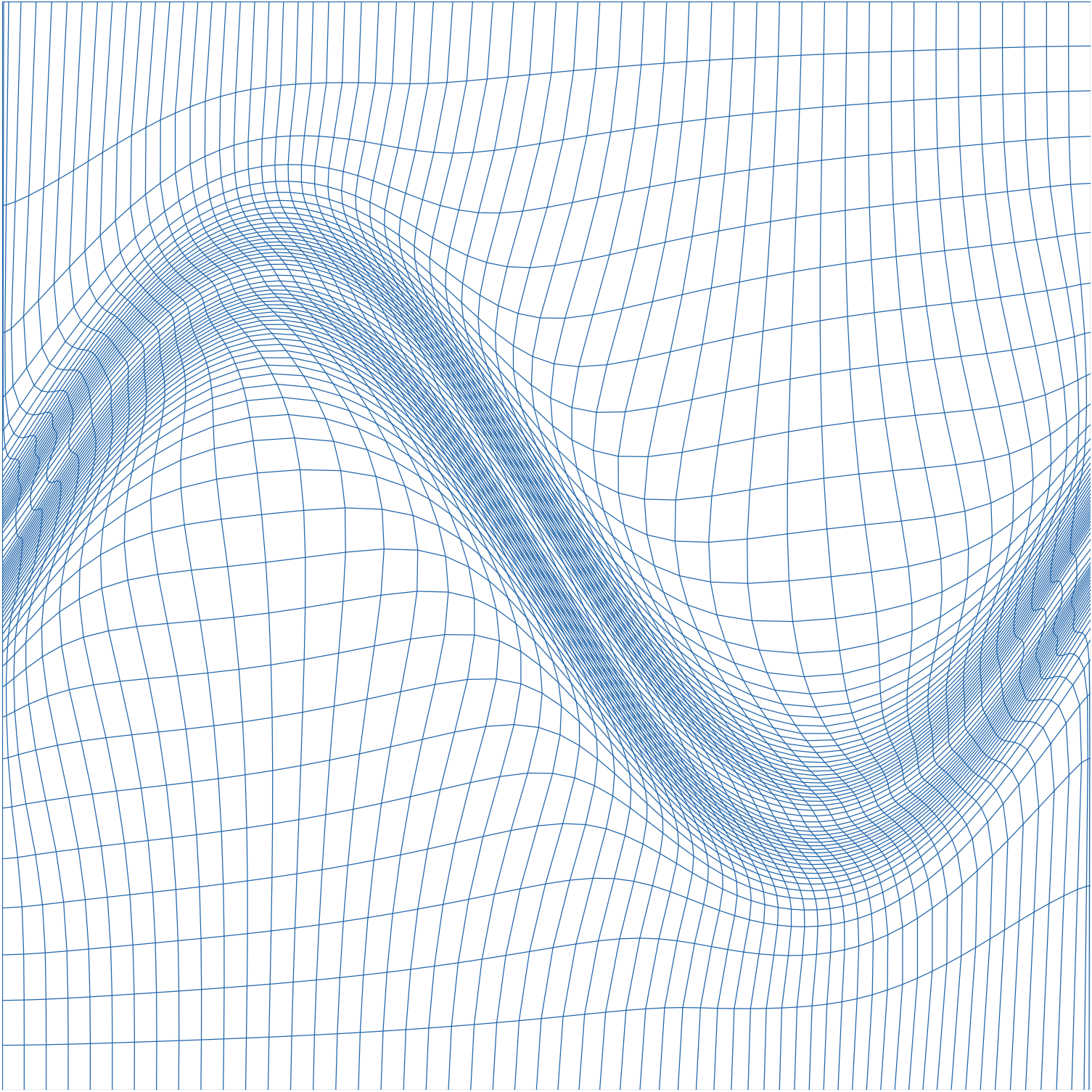}
\caption{$Q_1$ quadrilaterals.}\label{fig:legacy2d3Q1}\end{subfigure}
\caption{Hessian-based adaptations to the curved-front monitor
$u_{\rm front}$ in \eqref{eq:analyticmonitors} using $\M_{\rm adapt}$. 
(a) $P_1$ mesh. (b) $Q_1$ mesh. }
\label{fig:legacy2d3}\end{figure}

\begin{figure}[!htbp]\centering
\begin{subfigure}[t]{.485\linewidth}\centering
\includegraphics[width=\linewidth]{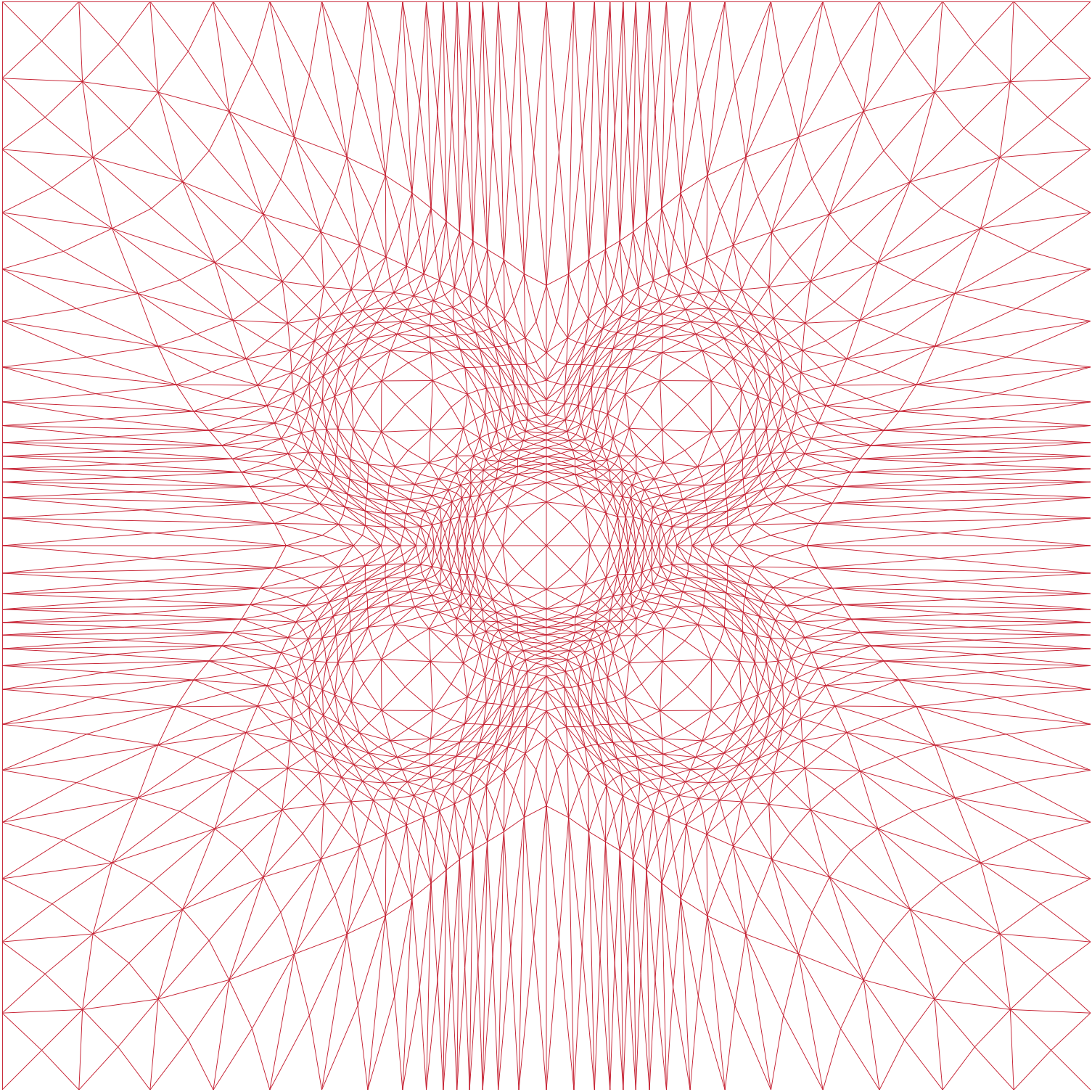}
\caption{$P_1$ triangles.}\label{fig:legacy2d5P1}\end{subfigure}\hfill
\begin{subfigure}[t]{.485\linewidth}\centering
\includegraphics[width=\linewidth]{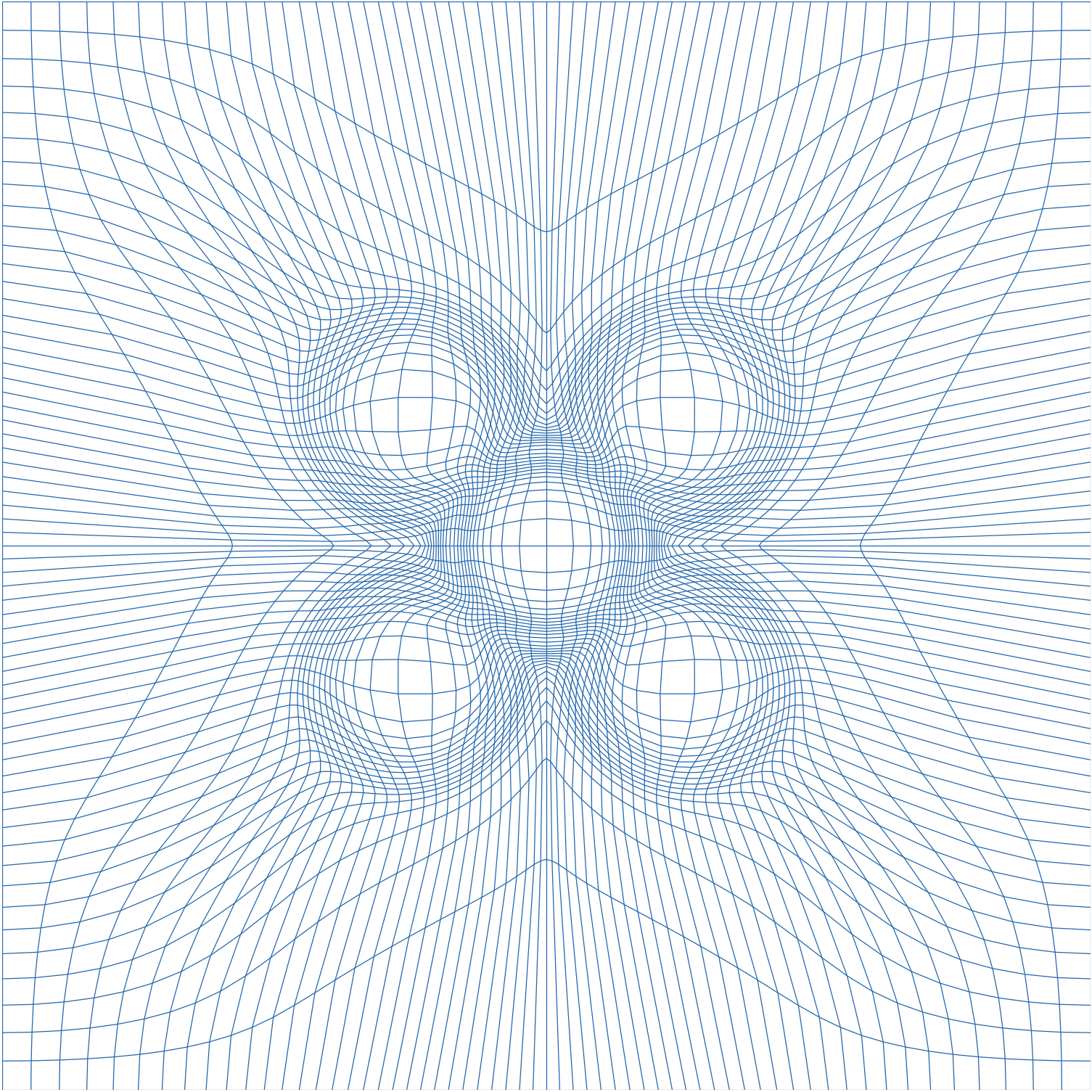}
\caption{$Q_1$ quadrilaterals.}\label{fig:legacy2d5Q1}\end{subfigure}
\caption{Hessian-based adaptations to the localized-rings
monitor $u_{\rm rings}$ in \eqref{eq:analyticmonitors} using $\M_{\rm adapt}$. 
(a) $P_1$ mesh. (b) $Q_1$ mesh. }
\label{fig:legacy2d5}\end{figure}

\begin{figure}[!htbp]\centering
\begin{subfigure}{.485\linewidth}\centering
\resultpanel{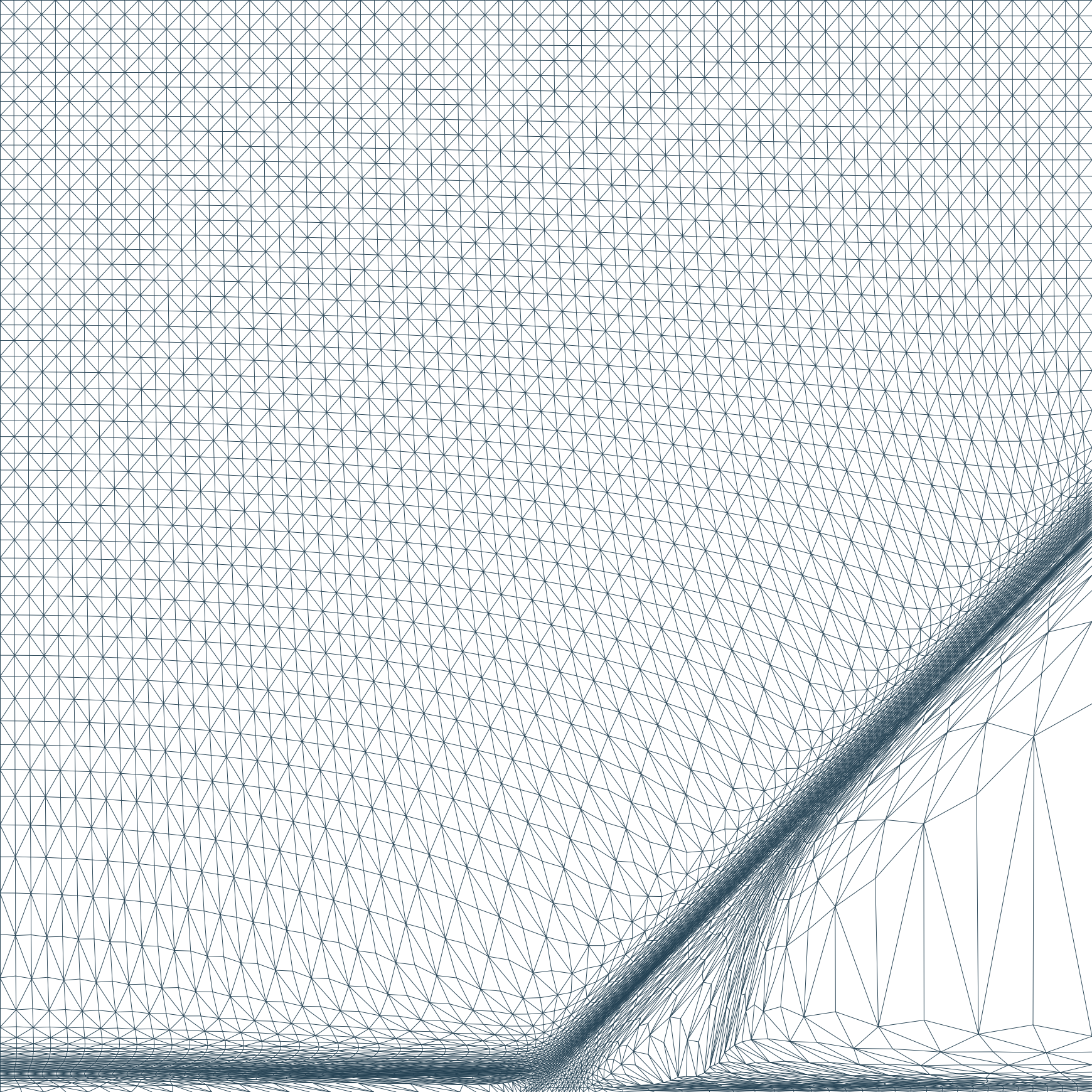}
\caption{$P_1$ triangles.}\end{subfigure}\hfill
\begin{subfigure}{.485\linewidth}\centering
\resultpanel{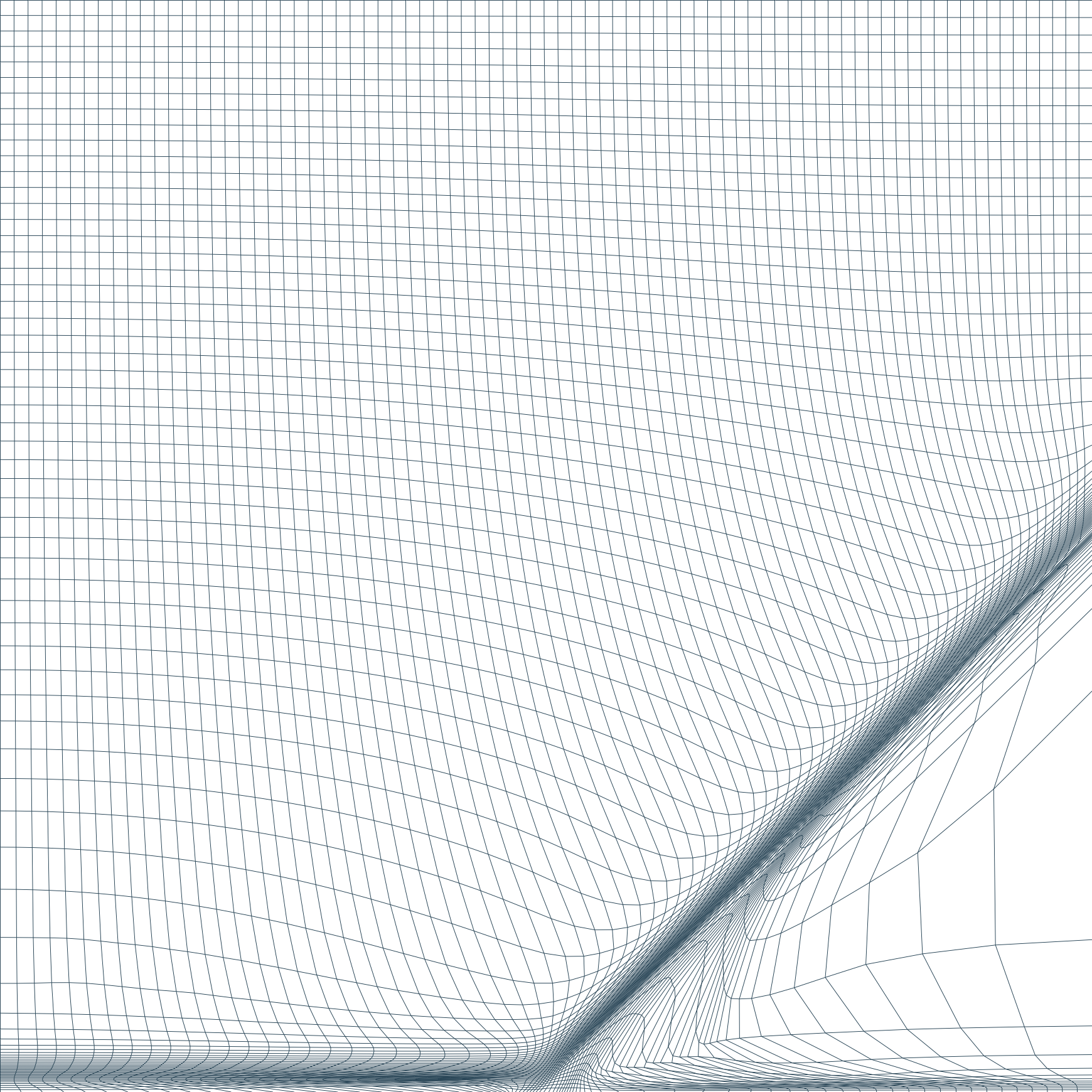}
\caption{$Q_1$ quadrilaterals.}\end{subfigure}
\caption{Direct mesh movement for the crossing-layer monitor
$u_{\rm cross}$ in \eqref{eq:analyticmonitors} using $\M_{\rm adapt}$. 
(a) $P_1$ mesh. (b) $Q_1$ mesh.}
\label{fig:ex2d4pair}\end{figure}

The $P_1$ and $Q_1$ adaptations to $u_{\rm front}$ and
$u_{\rm rings}$ are shown in Figures~\ref{fig:legacy2d3}
and~\ref{fig:legacy2d5}, respectively. Figure~\ref{fig:ex2d4pair}
shows direct movement for $u_{\rm cross}$. 

\FloatBarrier

For the quantitative direct-MMPDE experiments in Table~\ref{tab:planarquant},
we form $\M_{\rm adapt}$ from nodal values of each exact monitor,
using quadratic least-squares recovery (at most five patch rings,
condition limit 100, $\varepsilon_{\rm adapt}=1$),
three metric-smoothing passes and a $161\times161$ Cartesian grid
for interpolating the metric field. This interpolation grid is not the
moving-mesh resolution. In these direct experiments, both element
families have $81\times81=6561$ mesh vertices but different numbers
of elements. The metric field is evaluated at the moving physical
quadrature locations.
The diagnostic transition regions used in Table~\ref{tab:planarquant}
are, with
 $\eta_{0.99}=\operatorname{arctanh}(0.99)$,
\begin{align*}
 \mathcal R_{\rm front}&=\left\{(x,y):\left|y-\tfrac12-\tfrac14\sin(2\pi x)\right|<\eta_{0.99}/30\right\},\\
 \mathcal R_{\rm cross}&=\left\{(x,y):\min\{|y|,|x-y-\tfrac12|\}<\eta_{0.99}/60\right\},\\
 \mathcal R_{\rm rings}&=\left\{(x,y):\min_{(\alpha,\beta)\in S}
 \left|(4x-2-\alpha)^2+(4y-2-\beta)^2-\tfrac18\right|<\eta_{0.99}/30\right\}.
\end{align*}

The occupancy is the proportion of element centers in the indicated
region. The metric-mass fraction in Table~\ref{tab:planarquant} is the
integral of $\sqrt{\det\M(x)}$ over this same region divided
by its integral over the full domain, using each run's final prescribed
spatial metric field. The implementation
uses $\tau=0.01$, a positive scalar Jacobi mobility,
the three-node degree-two triangle rule for $P_1$ and a $3\times3$
tensor-product Gauss rule for $Q_1$ (nine quadrature nodes per cell),
and a maximum proposed coordinate displacement of $0.2/80$ per step. 
Recovered metric data are first transferred to Cartesian control nodes by 
piecewise-linear interpolation on the fixed recovery simplices. 
A $C^1$ log-Hermite SPD field is then constructed from these control values 
for moving-point metric evaluation. Thus the intermediate piecewise-linear transfer 
does not describe the final field differentiated in the direct spatial-field gradient. 
Both recovery and transfer introduce approximation beyond the prescribed analytic scalar.
No numerical PDE solve is involved in these tests.
For the six saved final meshes, a fixed-state re-evaluation of the
adaptation energy with independent Gauss--Duffy ($P_1$) or tensor-product
($Q_1$) quadrature yields differences of at most $0.115\%$ between
the production rule and the order-12 rule. The corresponding difference
between independent order-6 and order-12 rules is at most $0.00921\%$.
These comparisons use the same frozen metric and saved coordinates.
They are quadrature-sensitivity checks, not estimates of the error in
the continuous energy along the mesh trajectory. 
The results are presented in Table \ref{tab:planarquant}.

\begin{table}[!htbp]\centering\small
\caption{Direct planar response under $\M_{\rm adapt}$. Values are
initial/final for metric-volume CV, region occupancy (percent of element
centers), and analytic-monitor interpolation $L^2$ error. The last
region number is the metric-mass fraction on the final prescribed
field, not an imposed target or a convergence test. $Q_{\rm ali}^{95\%}$
is evaluated on the final mesh. $L^2$ values use independent order-six
quadrature, not errors of a PDE solution. Diagnostic regions are defined
in Section~\ref{sec:example_protocol}.}
\label{tab:planarquant}
\setlength{\tabcolsep}{3.1pt}
\begin{tabular}{@{}llrrrr@{}}\toprule
Monitor & Family & $\mathrm{CV}_M$ & $Q_{\rm ali}^{95\%}$ & Region: initial/final/metric & $L^2$: initial/final \\\midrule
Curved front & $P_1$ & 2.54 / 0.78 & 3.20 & 17.7 / 72.4 / 95.3 & 0.006184 / 0.0008882 \\
Curved front & $Q_1$ & 2.53 / 0.76 & 2.95 & 17.6 / 73.4 / 95.2 & 0.005423 / 0.0009713 \\
Crossing layers & $P_1$ & 3.79 / 3.54 & 3.20 & 8.4 / 40.0 / 77.5 & 0.01245 / 0.005535 \\
Crossing layers & $Q_1$ & 3.76 / 3.35 & 2.47 & 8.9 / 45.8 / 80.3 & 0.009606 / 0.0026 \\
Localized rings & $P_1$ & 2.69 / 0.94 & 3.61 & 14.8 / 56.9 / 95.2 & 0.01837 / 0.003076 \\
Localized rings & $Q_1$ & 2.66 / 0.85 & 3.39 & 14.6 / 61.4 / 95.0 & 0.01695 / 0.002693 \\
\bottomrule\end{tabular}
\end{table}
Each run completed its allowance of 2,000 accepted steps.
The differences between region occupancy and metric-mass fractions
show that visible concentration does not imply full equidistribution.
The interpolation errors decrease on these prescribed fields but are
not evidence of an application-independent accuracy gain. At the final
crossing-layer endpoint, the minimum certificate-to-stage-baseline
ratios are $0.1908$ for $P_1$ and $0.1680$ for $Q_1$. The localized
barrier is active on $4.57\%$ and $5.27\%$ of the respective final
certificate coefficients.
In the crossing-layer example, occupancy increases throughout the ten
stages, from 8.4\% to 40.0\% for $P_1$ and from 8.9\% to 45.8\%
for $Q_1$. Metric-volume variation is nonmonotone, and the calculation
does not establish metric equidistribution.

\FloatBarrier

\subsection{Further adaptation of a graded surface mesh}\label{sec:surface_mesh_example}
We next examine whether the unified MMPDE can further redistribute an
already graded mesh on a curved surface. 
We use the graph $\x(s,t)=(s,t,f(s,t))$ on $[-1,1]^2$, with
\begin{equation}\label{eq:surface_example_geometry}
 f(s,t)=0.08(s^2-t^2)+
 [0.10+0.25\exp(-6(s-0.15)^2)]
 \exp\{-160[t-0.35\sin(\pi s)]^2\}.
\end{equation}
The prescribed surface monitor is
\begin{equation}\label{eq:surface_example_field}
 u(s,t)=\exp\!\left[-\frac{(s-0.15)^2}{2(0.40)^2}
 -\frac{[t-0.35\sin(\pi s)]^2}{2(0.075)^2}\right].
\end{equation}
Both $P_1$ and $Q_1$ computations begin with an $81\times81$
chart grid that is already concentrated near the curved ridge. This
initial grading is prescribed analytically, rather than produced by a
prior MMPDE relaxation. Specifically, from uniformly spaced $(s,r)$
parameters we set
$\eta=\sinh(2.5r)/\sinh(2.5)$ and
$t=\eta+0.35\sin(\pi s)(1-\eta^2)$, and then lift the nodes onto
the graph. The experiment tests the additional redistribution produced
by moving this graded initial mesh according to the prescribed metric.
The reported mesh images and diagnostics compare this initial state
with the state after further MMPDE relaxation.

Both element-center occupancy
and chart metric-mass fraction are evaluated over the same full-width
ridge band
$\mathcal R_{\rm surf}=\{(s,t)\in[-1,1]^2:
|t-0.35\sin(\pi s)|\le0.15\}$.
Occupancy counts element centers, while metric mass integrates
$\sqrt{\det C(s,t)}$ over that region, divided by its integral
over the entire chart.
The mesh moves in chart coordinates and is lifted to the graph.
Whole-cell chart-Jacobian checks preserve orientation. A localized
barrier with $\kappa_{\rm chart}=10^{-3}$ is applied above a permanent
relative certificate floor of $10^{-3}$. The surface metric has a
positive normal extension. The results are obtained after three
stages of at most 512 accepted steps each and do not represent
stationary configurations. This graph-specific constraint is distinct
from a general retraction method for arbitrary embedded surfaces.

\begin{figure}[!htbp]\centering
\begin{subfigure}[t]{.485\linewidth}\centering
\includegraphics[width=\linewidth]{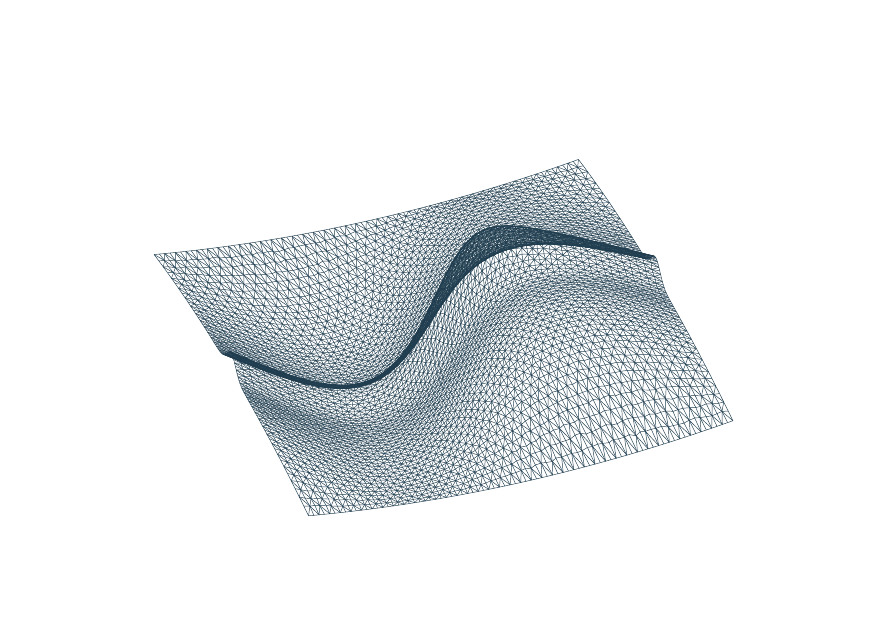}
\caption{Initially graded $P_1$ mesh.}\label{fig:surface-P1-initial}\end{subfigure}\hfill
\begin{subfigure}[t]{.485\linewidth}\centering
\includegraphics[width=\linewidth]{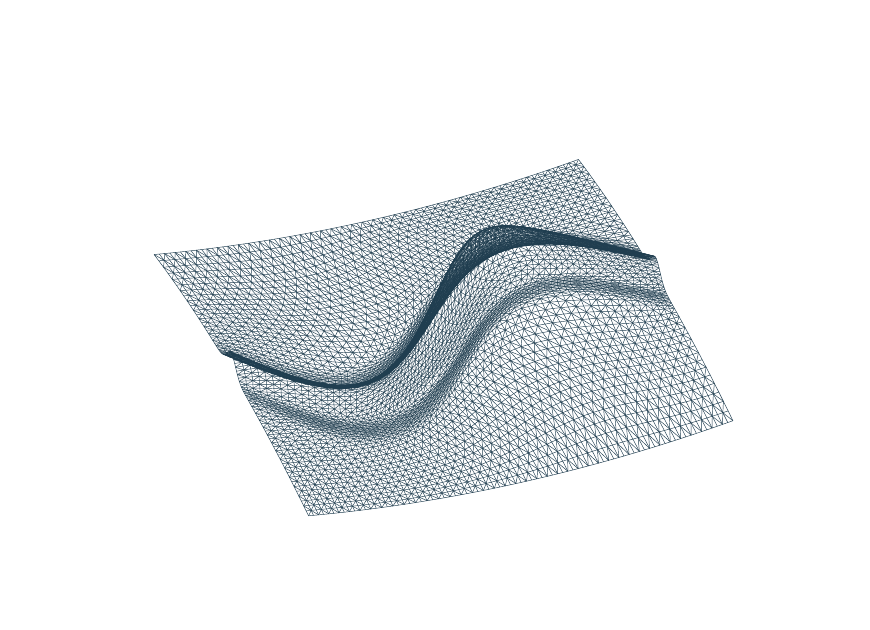}
\caption{$P_1$ after further movement.}\label{fig:surface-P1-moved}\end{subfigure}
\medskip
\begin{subfigure}[t]{.485\linewidth}\centering
\includegraphics[width=\linewidth]{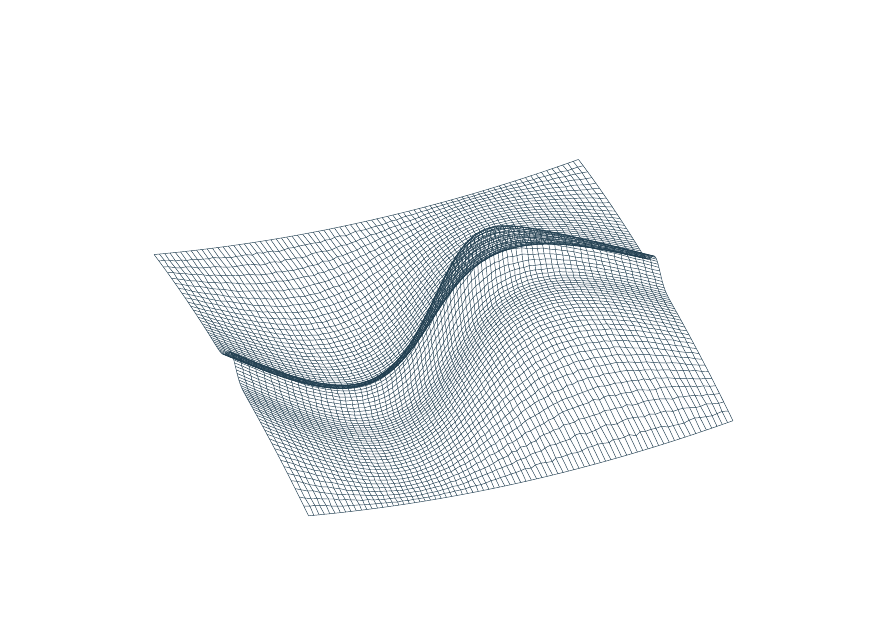}
\caption{Initially graded $Q_1$ mesh.}\label{fig:surface-Q1-initial}\end{subfigure}\hfill
\begin{subfigure}[t]{.485\linewidth}\centering
\includegraphics[width=\linewidth]{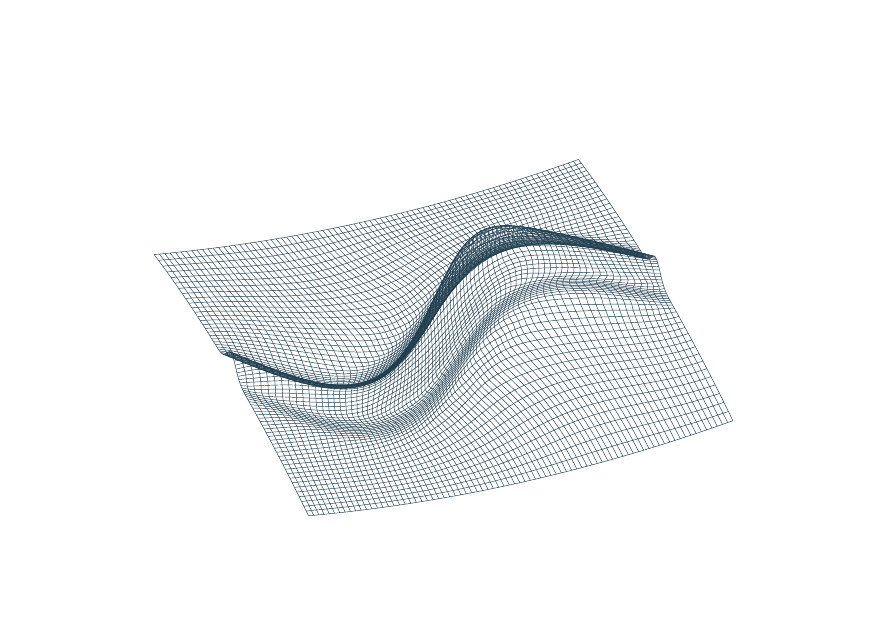}
\caption{$Q_1$ after further movement.}\label{fig:surface-Q1-moved}\end{subfigure}
\caption{Further metric-driven redistribution of initially graded
graph-surface meshes, viewed obliquely (elevation $38^\circ$, azimuth $-110^\circ$),
without height exaggeration.
(a) Initially graded $P_1$ mesh. (b) $P_1$ mesh after further movement.
(c) Initially graded $Q_1$ mesh. (d) $Q_1$ mesh after further movement.}
\label{fig:surface-theta-third-full}\end{figure}

\begin{figure}[!htbp]\centering
\begin{subfigure}[t]{.485\linewidth}\centering
\includegraphics[width=\linewidth]{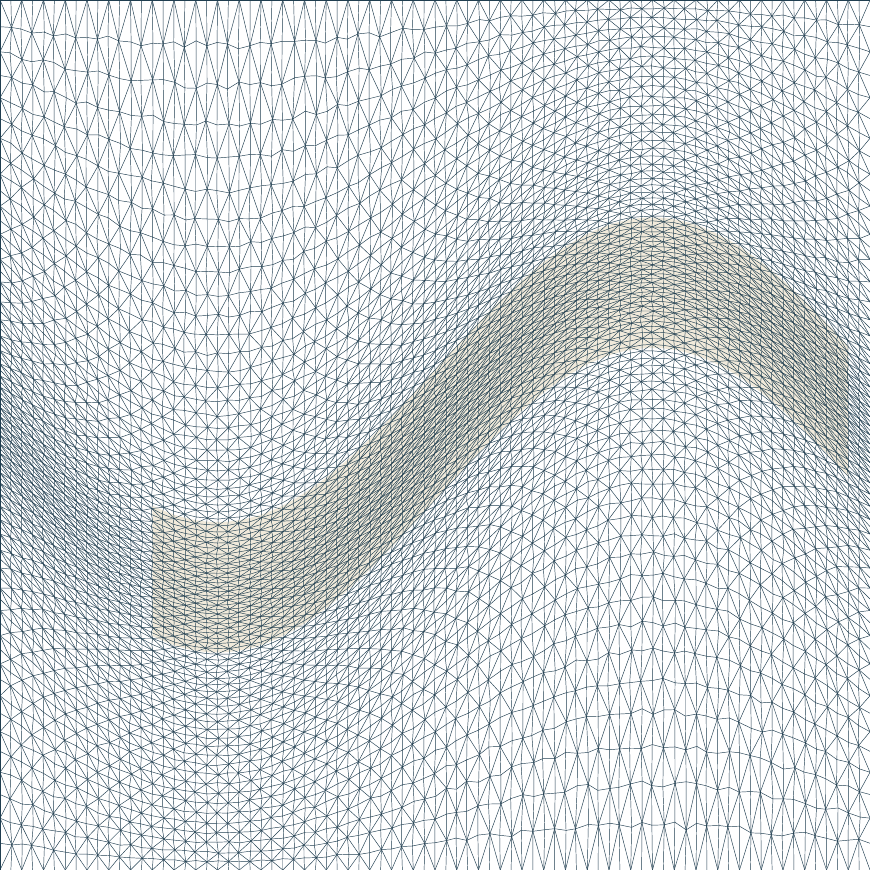}
\caption{Initially graded $P_1$ mesh.}\label{fig:surface-top-P1-initial}\end{subfigure}\hfill
\begin{subfigure}[t]{.485\linewidth}\centering
\includegraphics[width=\linewidth]{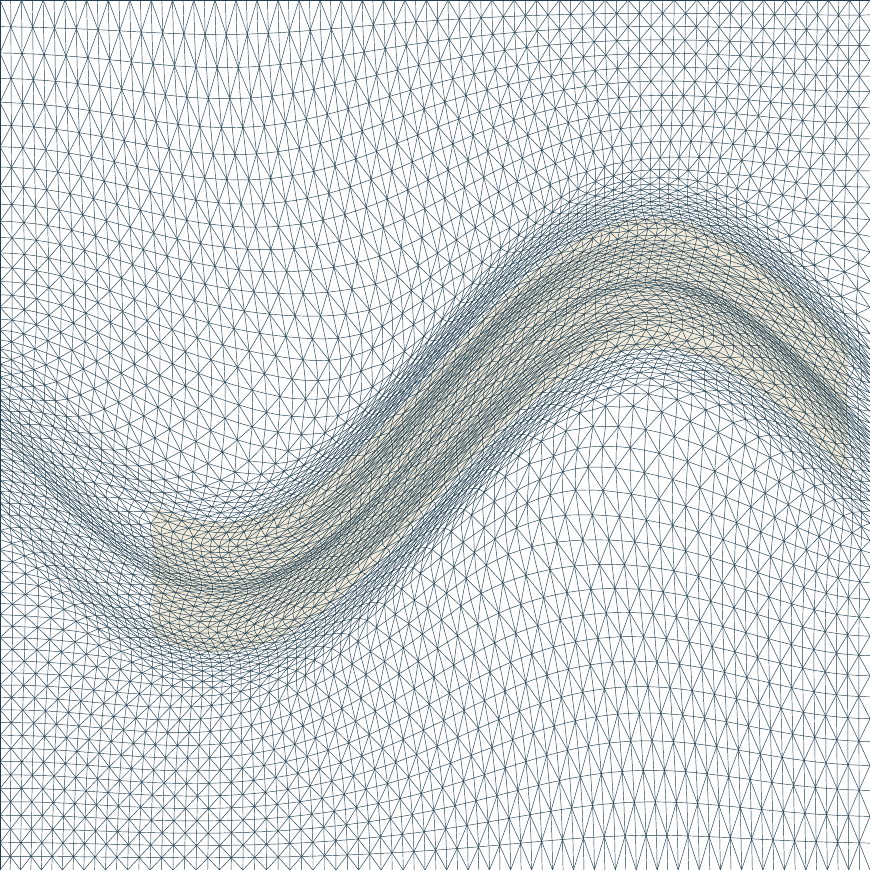}
\caption{$P_1$ after further movement.}\label{fig:surface-top-P1-moved}\end{subfigure}
\medskip
\begin{subfigure}[t]{.485\linewidth}\centering
\includegraphics[width=\linewidth]{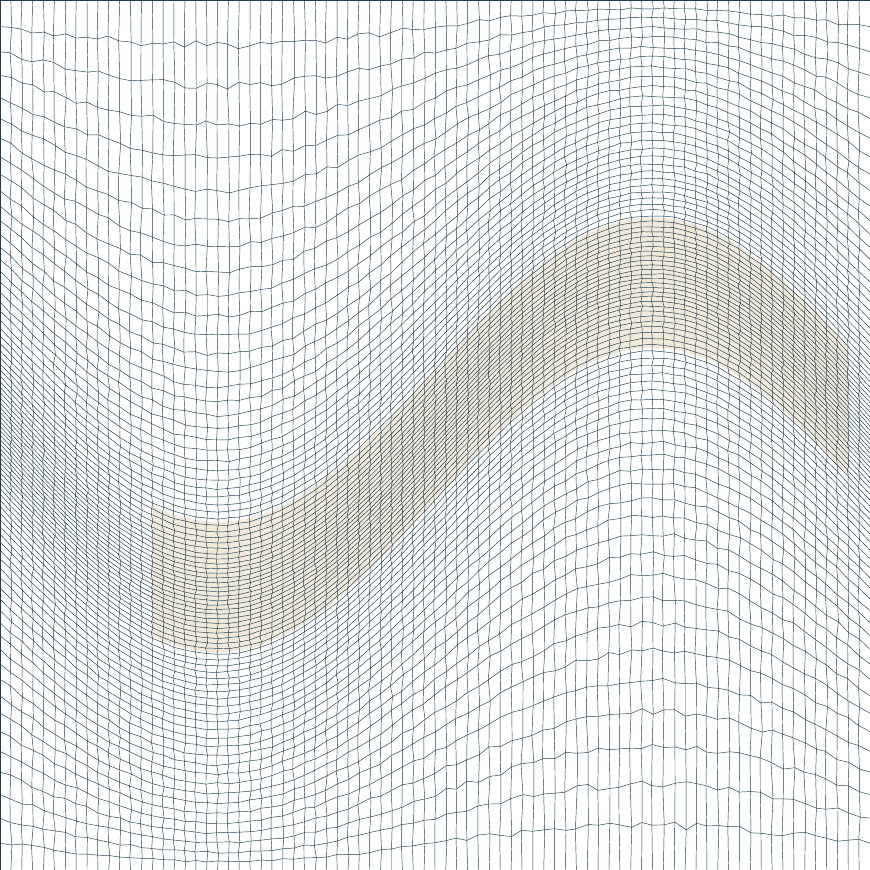}
\caption{Initially graded $Q_1$ mesh.}\label{fig:surface-top-Q1-initial}\end{subfigure}\hfill
\begin{subfigure}[t]{.485\linewidth}\centering
\includegraphics[width=\linewidth]{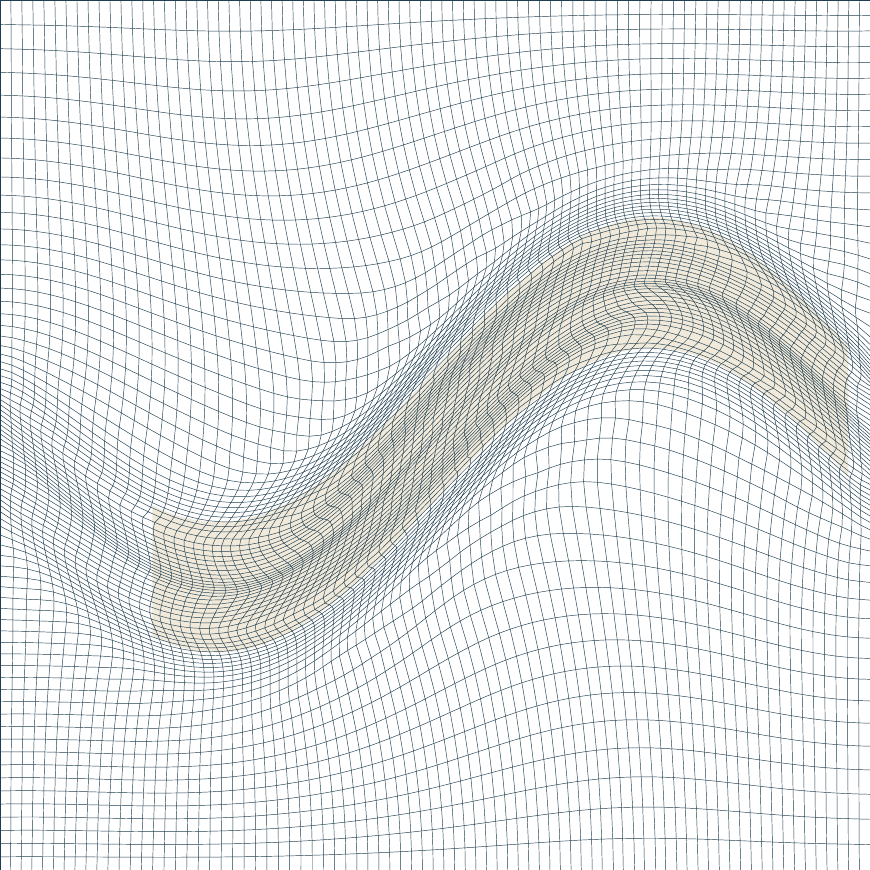}
\caption{$Q_1$ after further movement.}\label{fig:surface-top-Q1-moved}\end{subfigure}
\caption{Top-down chart-coordinate views of the initial and further
adapted surface meshes in Figure~\ref{fig:surface-theta-third-full}.
The curved band is centered on $t=0.35\sin(\pi s)$. The initial
meshes already contain ridge concentration, and the moved meshes
show the additional redistribution obtained by relaxation.}
\label{fig:surface-topdown}\end{figure}

The final occupancy of the full-width ridge band rises from
32.66\% to 35.90\% for $P_1$ and from 32.94\% to 37.03\% for
$Q_1$, compared with final metric-mass fractions of 50.98\% and
51.03\%, respectively. Metric-volume CV decreases from
$0.918$ to $0.626$ for $P_1$ and from $0.907$ to $0.591$ for $Q_1$.
The maximum sampled surface-normal deviations are $44.26^\circ$
and $37.98^\circ$, respectively. The results show visible constrained
surface redistribution, but not improved geometric approximation or
validation of the general closed-surface case.
Figure~\ref{fig:surface-theta-third-full} displays the 3D view of the meshes,
while Figure~\ref{fig:surface-topdown} shows their corresponding 2D top-down chart
views.

\FloatBarrier

\subsection{Bulk three-dimensional test}\label{sec:bulk3d}
We next demonstrate mesh movement in three dimensions.
On the unit cube we prescribe the scalar field
\begin{equation}\label{eq:bulkscalar}
 u(x,y,z)=\tanh\{-7[1/2-\tanh(-7(y-1/2-\tfrac14\sin2\pi x))]\}.
\end{equation}
This monitor is independent of $z$ and therefore tests bulk
isoparametric validity and extrusion of anisotropic adaptation,
not a fully three-directional metric variation. The direct runs use $17^3=4913$ vertices, five metric-update stages,
and up to 2000 accepted steps per stage. The Hessian regularization is
$\varepsilon_{\rm adapt}=1$. The $P_1$ adaptation integral uses a
four-point tetrahedral rule and the $Q_1$ integral uses two Gauss
points per reference direction. Every saved hexahedral endpoint has
positive whole-element Bernstein coefficients, checked by exact
rational arithmetic on the stored coordinates. The smallest final
coefficient lower bound is $4.01043\times10^{-7}$. The normalized
certificate margin is $\min_c(b_c/b_c^0)$, with $b_c^0$ fixed at the
beginning of the corresponding stage. The exact-arithmetic checks
certify the saved local geometry, not the complete floating-point
trajectory or global injectivity. The $Q_1$ energy and spatial-field gradient were integrated using
two Gauss points per reference direction. For the final saved mesh,
this gives an adaptation energy of $403316.0824$, whereas nine points
per direction give $438997.6693$, a relative difference of $8.8471\%$
with respect to the lower-order value. On the initial uniform mesh,
the corresponding discrepancy is $0.2545\%$. Table \ref{tab:bulkquant} summarizes the results.
The projected adaptation forces on the final mesh differ by $84.21\%$ relative to the
higher-order force norm. Independent differentiation confirms that
the discrete gradient is consistent with the energy under its own
quadrature rule, but the final $Q_1$ mesh is quadrature-limited.
The example illustrates locally valid hexahedral movement, not a
quadrature-converged minimizer. No corresponding quadrature-sensitivity
audit was performed for the four-point tetrahedral $P_1$ rule.

\begin{table}[!htbp]\centering\small
\caption{Bulk direct results at 4913 vertices and 2,000 accepted
inner steps. $Q_{\rm ali}^{95\%}$ and normalized certificate margins
refer to the final saved configuration. All endpoints pass the
independent whole-element check.}
\label{tab:bulkquant}
\begin{tabular}{@{}lrrrrr@{}}\toprule
Family & Elements & Initial CV & Final CV & Final $Q_{\rm ali}^{95\%}$ & Min normalized certificate \\\midrule
$P_1$ & 24,576 & 9.48 & 2.96 & 4.73 & 0.063 \\
$Q_1$ & 4,096 & 13.69 & 4.69 & 4.47 & 0.019 \\
\bottomrule\end{tabular}\end{table}
\begin{figure}[!htbp]\centering
\begin{subfigure}{.485\linewidth}\centering
\resultpanel{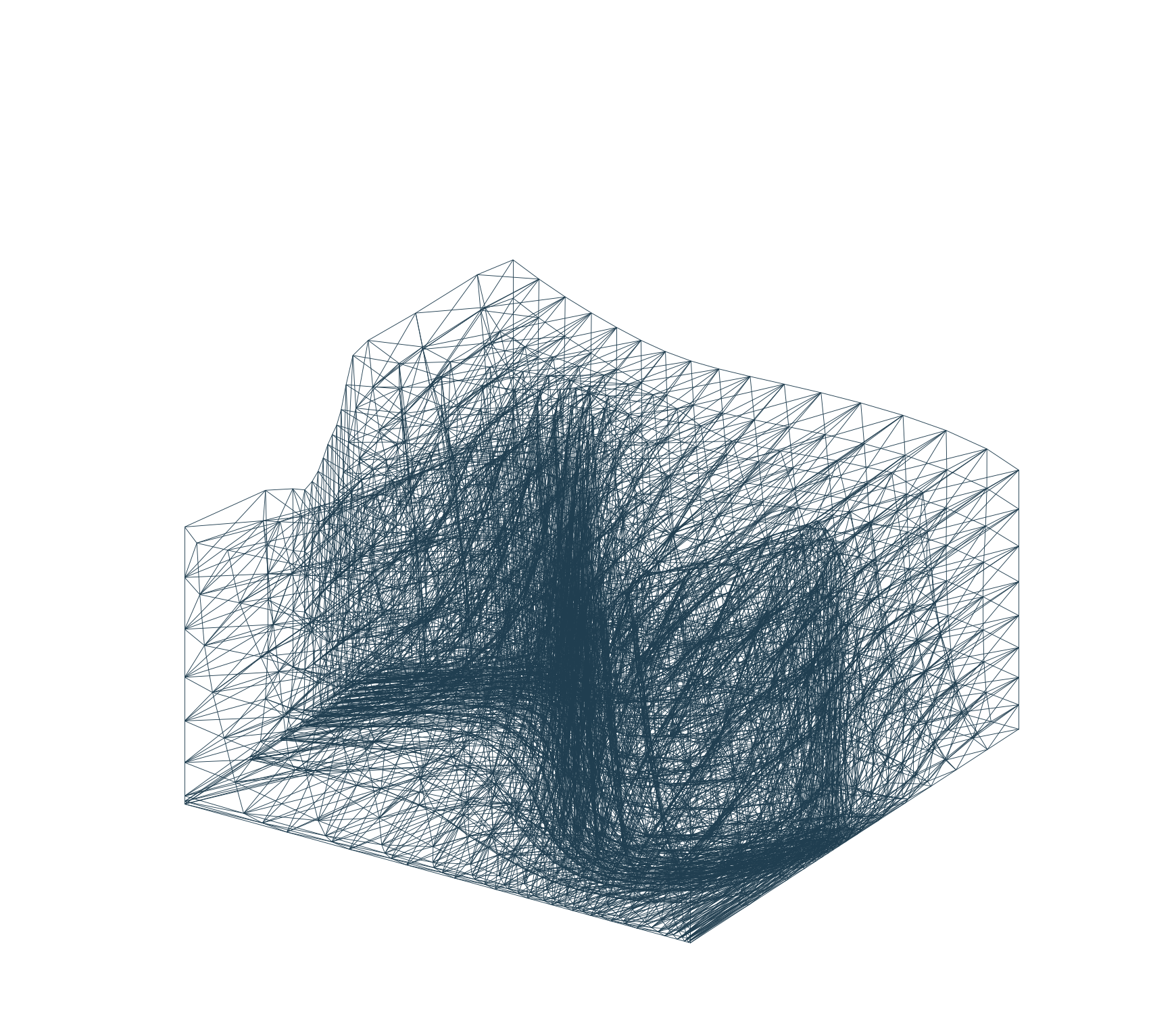}\caption{$P_1$ tetrahedra.}\end{subfigure}\hfill
\begin{subfigure}{.485\linewidth}\centering
\resultpanel{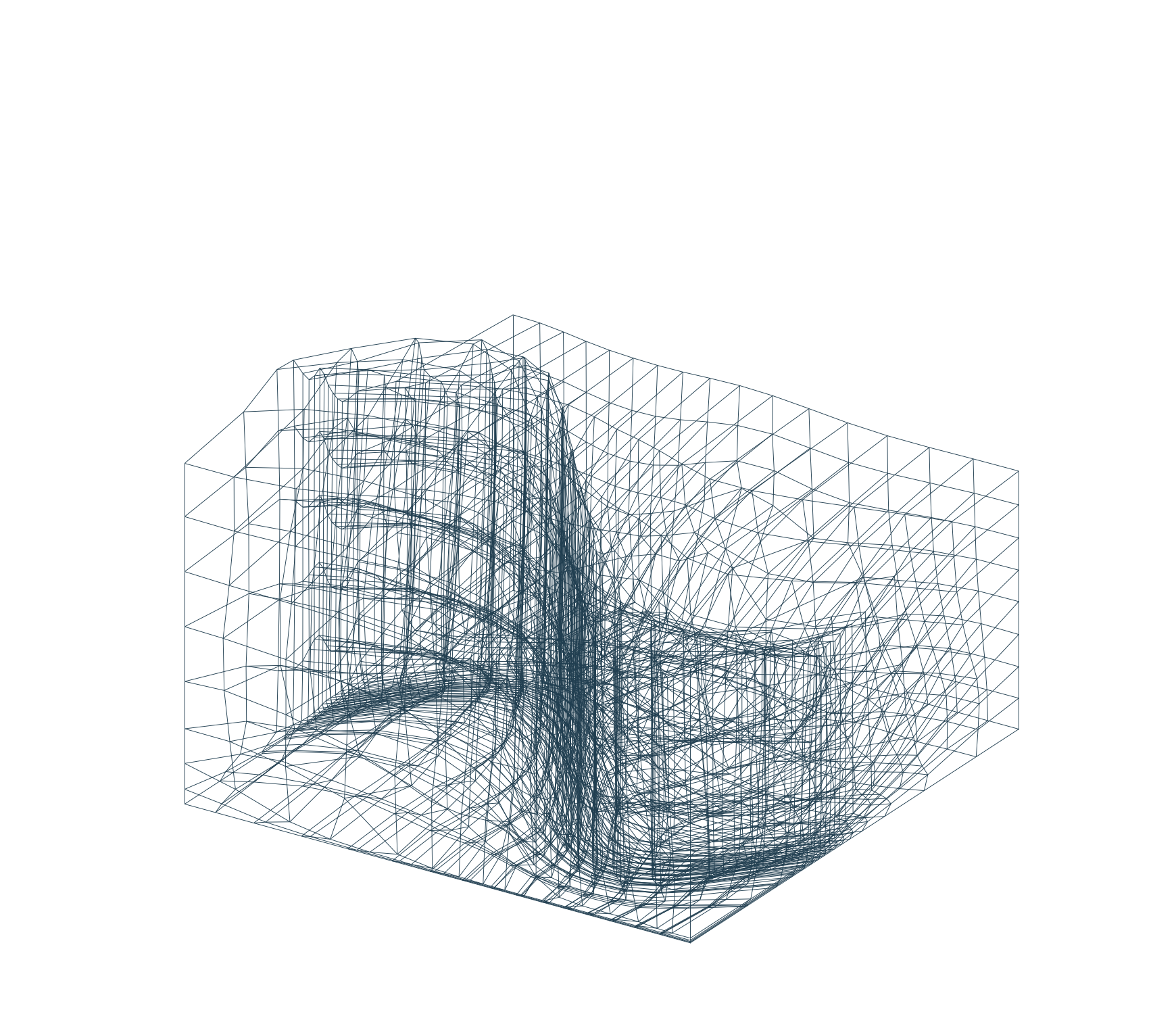}\caption{$Q_1$ hexahedra.}\end{subfigure}
\caption{Cutaway views of the moved unit-cube meshes, showing cells
whose initial centroids satisfy $z\le1/2$. The prescribed scalar
is independent of $z$. (a) $P_1$ mesh. (b) $Q_1$ mesh. 
The figure illustrates metric-driven deformation
of hexahedral and tetrahedral meshes with verified local validity.}
\label{fig:ex3d5pair}\end{figure}
\FloatBarrier

\subsection{Diffusion-aligned metric demonstration}\label{sec:dmp_examples}
Lastly, we present mesh movement based on the solution-dependent $\M_{\rm DMP+adapt}$ metric in
\eqref{eq:MDMPadaptimpl}, with $\theta=1/3$, for the stationary
anisotropic diffusion problem
\begin{equation}\label{eq:LHex51pde}
 -\nabla\cdot(D\nabla u)=0 \quad\text{in }\Omega,
 \qquad\Omega=[0,1]^2\setminus(4/9,5/9)^2,
\end{equation}
with $u=0$ on the outer boundary and $u=2$ on the hole boundary.
The diffusion tensor is
\begin{equation}\label{eq:LHdiffusion}
 D(x,y)=R(\beta)\operatorname{diag}(1000,1)R(\beta)^T,
 \quad R(\beta)=\begin{pmatrix}\cos\beta&-\sin\beta\\
 \sin\beta&\cos\beta\end{pmatrix}.
\end{equation}
We use both $\beta=\pi/4$ and the spatial orientation
\begin{equation}\label{eq:LHvariabletheta}
 \beta(x,y)=\frac{\pi}{2}-\pi\cos x\sin y.
\end{equation}
The Galerkin discretization uses the physical isoparametric
gradients and the spatial diffusion tensor at quadrature points,
with 12-point $P_1$ and 16-point $Q_1$ element rules. The
solution-dependent metric is refreshed after accepted mesh movement.
Each mesh contains 5120 elements. The $P_1$ and $Q_1$ grids have
2640 and 5280 vertices, respectively. The five metric-update stages
use the localized validity barrier and impose maximum whole-cell
width and height $0.2$ and boundary-node gap $0.15$.

\begin{figure}[!htbp]\centering
\begin{subfigure}[t]{.485\linewidth}\centering
\includegraphics[width=\linewidth]{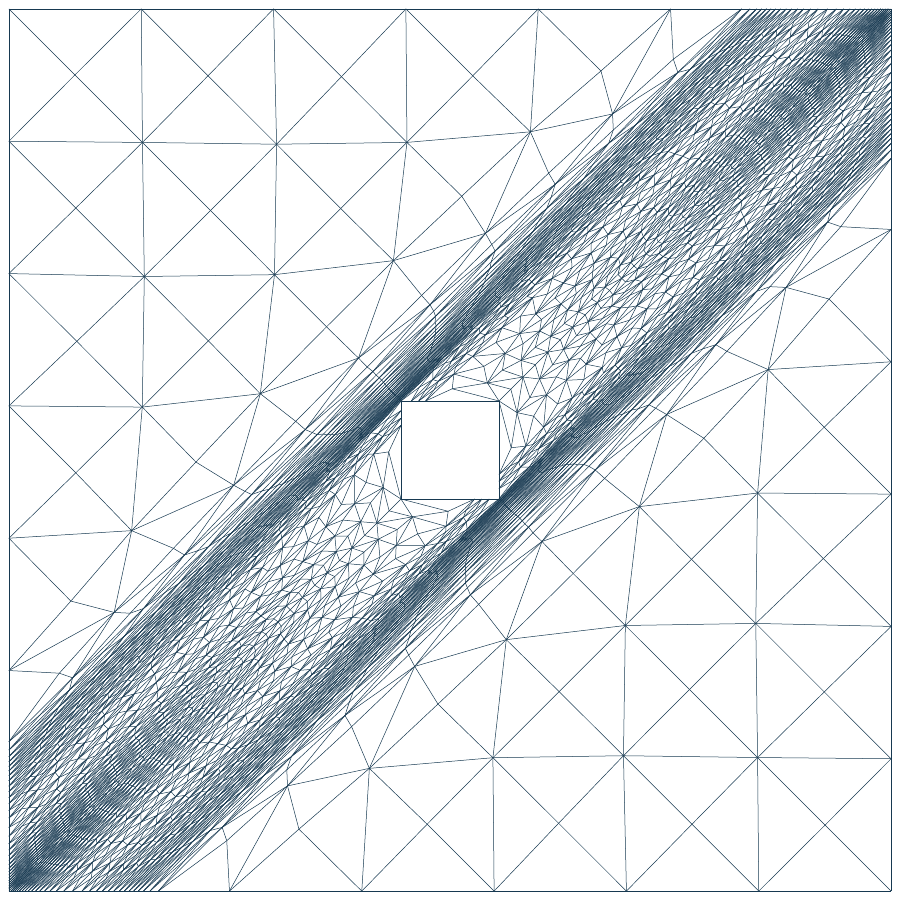}
\caption{$P_1$: $U_-=2.03164\times10^{-4}$.}\end{subfigure}\hfill
\begin{subfigure}[t]{.485\linewidth}\centering
\includegraphics[width=\linewidth]{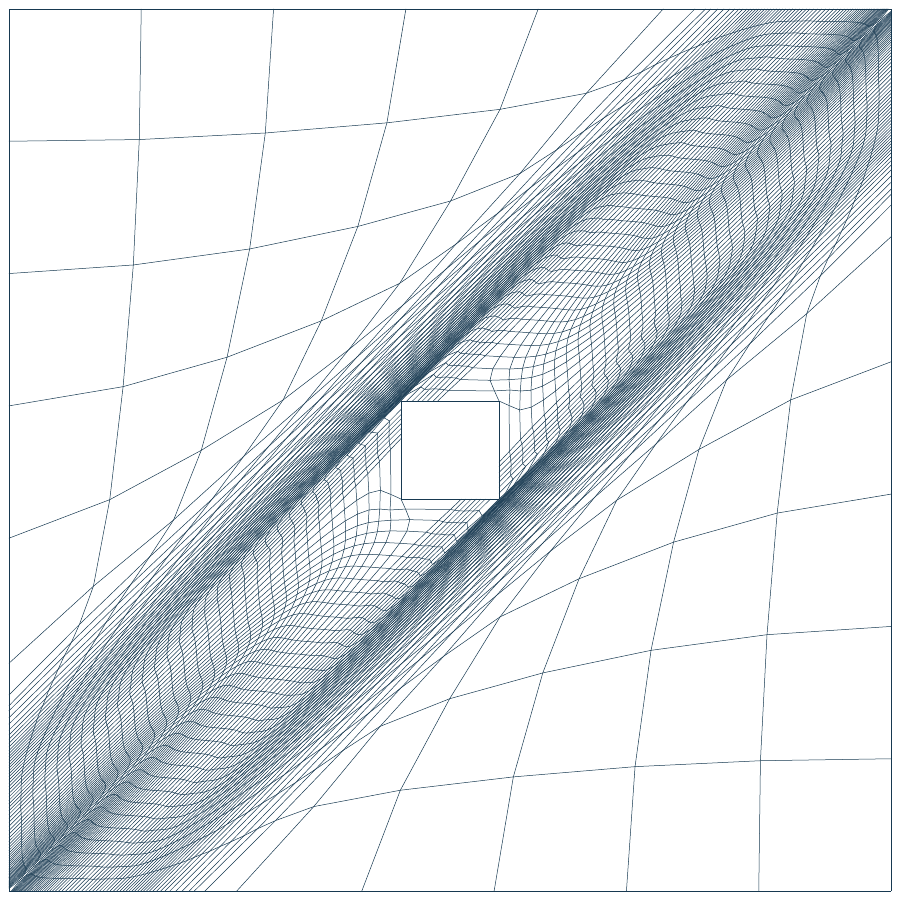}
\caption{$Q_1$: $U_-=9.21391\times10^{-6}$.}\end{subfigure}
\medskip
\begin{subfigure}[t]{.485\linewidth}\centering
\includegraphics[width=\linewidth]{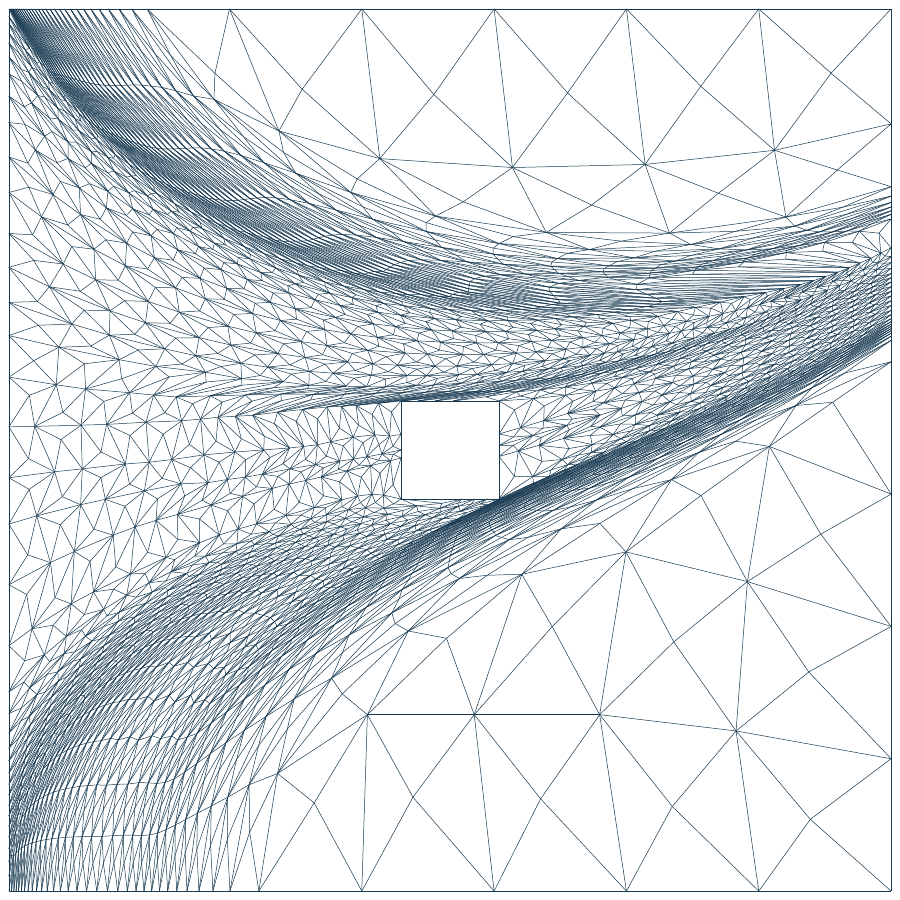}
\caption{$P_1$: $U_-=8.78471\times10^{-5}$.}\end{subfigure}\hfill
\begin{subfigure}[t]{.485\linewidth}\centering
\includegraphics[width=\linewidth]{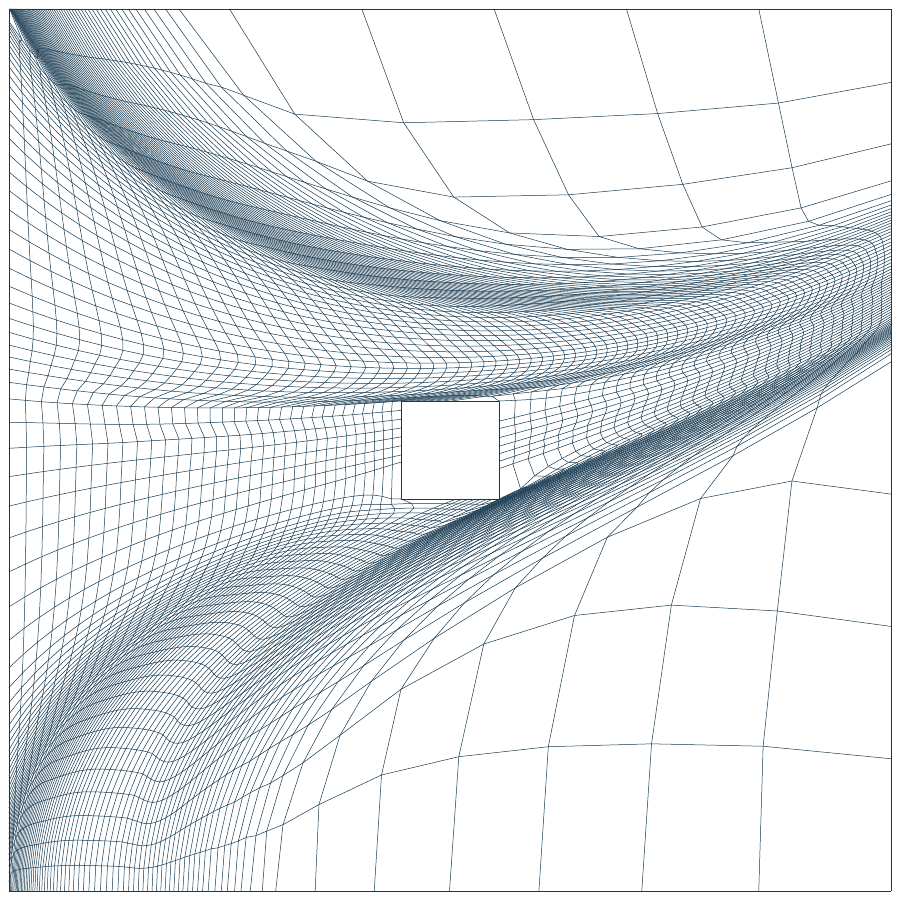}
\caption{$Q_1$: $U_-=1.50145\times10^{-6}$.}\end{subfigure}
\caption{Metric-driven meshes for $\M_{\rm DMP+adapt}$ with
$\theta=1/3$ and 5120 elements. Top row: constant diffusion
orientation with $\beta=\frac{\pi}{4}$; bottom: spatial orientation \eqref{eq:LHvariabletheta}.
Left: $P_1$ mesh; right: $Q_1$ mesh.}
\label{fig:dmpadapt-four-panels}\end{figure}

The raw undershoot is
$U_- = \max\{0,-\min_{x_i}u_h(x_i)\}$.
The four values in Figure~\ref{fig:dmpadapt-four-panels} describe the
individual solutions. The reported undershoots are smaller than those obtained on the
corresponding initial uniform meshes, illustrating the effect of
$\M_{\rm DMP+adapt}$. However, these results neither compare approximation
orders nor establish satisfaction of a discrete maximum principle.

\FloatBarrier

\subsection{A hexahedral validity counterexample}\label{sec:hexexample}
\begin{figure}[!htbp]\centering
\includegraphics[width=.88\linewidth]{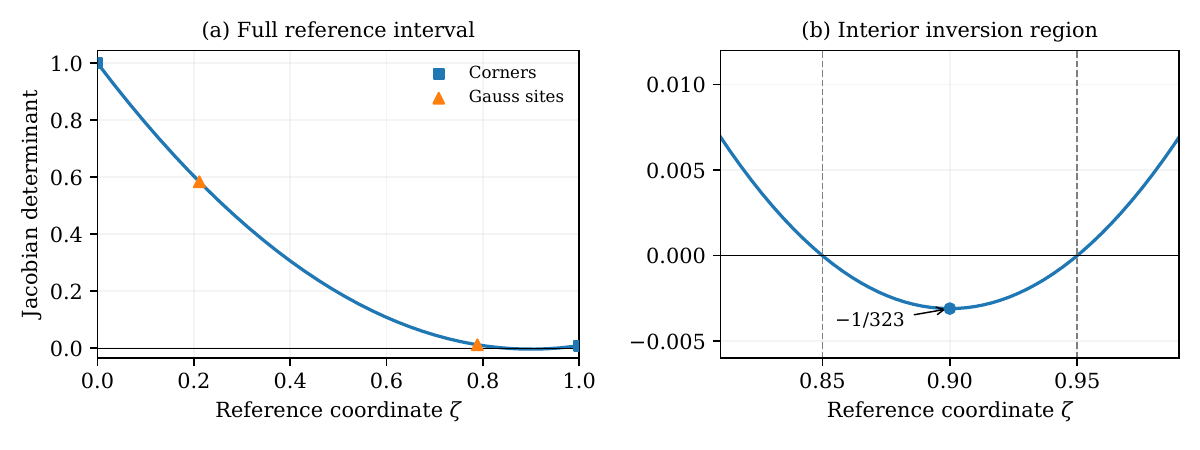}
\caption{Jacobian determinant of the explicit trilinear map
\eqref{eq:hex_sample_counterexample}. The determinant is positive at
the corner planes ($\zeta=0,1$) and the two Gauss planes but negative
for $17/20<\zeta<19/20$. Sampling does not certify whole-element orientation.}
\label{fig:hex-diagnostics}\end{figure}

Here, we illustrate why pointwise Jacobian sampling is insufficient to
establish validity of a trilinear hexahedron. Consider the explicit
$Q_1$ map on $[0,1]^3$,
\begin{equation}\label{eq:hex_sample_counterexample}
 F(\xi,\eta,\zeta)
 =\left(\left(1-\frac{20}{17}\zeta\right)\xi,
         \left(1-\frac{20}{19}\zeta\right)\eta,\zeta\right).
\end{equation}
Its Jacobian determinant is independent of $\xi$ and $\eta$:
\begin{equation}\label{eq:hex_sample_det}
 \det D_{\xiB}F
 =\left(1-\frac{20}{17}\zeta\right)
  \left(1-\frac{20}{19}\zeta\right).
\end{equation}
At every reference corner this determinant is either $1$ or
$3/323$, and at all eight two-point tensor-product Gauss nodes it
is positive. Nevertheless,
\[
 \det D_{\xiB}F(\xi,\eta,9/10)=-1/323<0.
\]
Indeed, the determinant vanishes on the planes $\zeta=17/20$
and $\zeta=19/20$ and is negative between them, as shown in Figure \ref{fig:hex-diagnostics}. 
Thus, even positive corner and Gauss-site determinants cannot certify orientation of the
whole element. Positivity of all Bernstein coefficients is sufficient,
whereas a negative coarse Bernstein coefficient is inconclusive and
can require subdivision. In this example, the quadratic Bernstein coefficients in the
$\zeta$ direction are $1$, $-37/323$, and $3/323$, repeated in
the other two reference directions. The negative coefficient is
consistent with, but by itself would not prove, the observed inversion.

The discrete acceptance test preserves local orientation of accepted
meshes independently of whether a barrier is included.  The barrier
instead supplies the separate continuous-time positive-margin result
of Theorem~\ref{thm:certificate}. No numerical example is used to
claim that it uniquely prevents inversion or step-size stagnation.
\FloatBarrier

\section{Conclusions}\label{sec:conclusion}

We have developed a unified isoparametric MMPDE framework for direct
motion of simplicial $P_1$ and tensor-product $Q_1$ meshes under a
prescribed symmetric positive-definite metric tensor $\M$. The
formulation uses the same intrinsic pullback metric and
reference-integrated equidistribution--alignment energy for planar,
bulk, and embedded elements. A complete physical-coordinate gradient
accounts for changes in a prescribed spatial metric field at moving
quadrature points, and the projected gradient flow respects admissible
nodal tangent directions under the stated smoothness assumptions.
The framework does not depend on a particular metric-construction
algorithm. Accuracy-oriented and DMP-oriented metrics, for example,
may be supplied without changing the mesh-motion formulation.

A new analytical contribution is the exact-energy nonsingularity
theorem for $Q_1$ maps. On a fixed finite mesh with bounded nodal
coordinates and a uniformly SPD metric, bounded exact adaptation
energy prevents loss of local Jacobian rank when $p\ge3/2$ for
bilinear elements or $p\ge4/3$ for trilinear elements. The thresholds
are sharp for general maps, and the bilinear result also applies to
embedded surface cells as preservation of local immersion. The proof
uses the singular-value structure of the alignment density and the
geometry of rank-deficient $Q_1$ Jacobians, including the singular
segments forced by rank-one trilinear Jacobians. This result supplies
a finite-mesh existence margin but no explicit certificate or
mesh-refinement-uniform bound.

The exact-energy guarantee and computational validity tests have
distinct roles. Since finite quadrature need not capture the energy
blowup near a small singular region, the accepted discrete meshes
still require whole-element validity tests. Positive corner
Jacobians give an exact test for planar bilinear quadrilaterals;
positive Bernstein determinant coefficients, with subdivision when
necessary, provide a sufficient test for trilinear hexahedra and
projected surface patches. A fixed-stage optional certificate barrier
yields a quantitative coefficient margin, whereas hard trial
acceptance preserves local validity even without that barrier. These
statements concern local element regularity and orientation, not
global injectivity or avoidance of surface self-contact.

The numerical results show that the prescribed metric can drive
visible redistribution for the supported element families in one dimension, 
on planar meshes, on three-dimensional volume meshes, and on a graph surface. The graded-start surface example demonstrates further redistribution
under chart constraints, while
an explicit hexahedral counterexample demonstrates that positivity at
corners and Gaussian quadrature nodes does not certify validity
throughout a trilinear element. The reported relaxations have finite
step budgets and do not establish exact metric equidistribution or
stationarity. In particular, the three-dimensional $Q_1$ energy and
gradient exhibit non-negligible quadrature sensitivity. An equivalent
sensitivity study has not been carried out for the four-point $P_1$
tetrahedral rule. Numerical accuracy and discrete maximum-principle
properties remain dependent on the prescribed metric and the
underlying PDE discretization, rather than following from the mesh
motion alone.

\appendix

\section{Proof details and sharpness of exact-energy \texorpdfstring{$Q_1$}{Q1} coercivity}
\label{app:Q1coercivity}

This appendix supplies the nonintegrability lemma used in
Theorem~\ref{thm:Q1exactcoercivity}. The argument is local in
reference coordinates and applies equally to bulk Jacobians and to
full-column-rank Jacobians of polynomial elements embedded in a higher
ambient dimension. Singularities on the boundary of the closed
reference cube require one-sided neighborhoods. They cannot be
excluded by imposing assumptions only at interior quadrature points.

\subsection*{Singular-value and integrability estimates}
For a full-column-rank matrix $\B\in\R^{d\times m}$, with
$\sigma_1\ge\cdots\ge\sigma_m>0$, recall
\[
H_p(\B)=\left(\prod_{j=1}^m\sigma_j\right)
\left(\sum_{j=1}^m\sigma_j^{-2}\right)^{mp/2}.
\]
At rank loss, set $H_p=+\infty$.

The first lemma relates the alignment contribution to the smallest singular value, 
given a lower bound on the nondegenerate singular directions.

\begin{lemma}[Rank-dependent singular-value bound]
\label{lem:Q1rankestimate}
Let $1\le r<m$, and suppose the matrix field $\B$ satisfies
$\sigma_r(\B(\xiB))\ge a>0$ on a set $V\subset\Khat$. Then, at every
full-rank point of $V$,
\begin{equation}\label{eq:Q1rankbound}
H_p(\B)\ge a^r\sigma_m(\B)^{-(mp-m+r)}.
\end{equation}
If $\B$ is $L$-Lipschitz and has deficient rank at $z\in\Khat$,
then $\sigma_m(\B(\xiB))\le L|\xiB-z|$. If $\B$ has deficient
rank on a compact segment $L_*\subset\Khat$, the corresponding estimate is
$\sigma_m(\B(\xiB))\le L\dist(\xiB,L_*)$.
\end{lemma}
\begin{proof}
For $j\le r$ we have $\sigma_j\ge a$, and for $j>r$ we have
$\sigma_j\ge\sigma_m$. Hence
\[
\prod_{j=1}^m\sigma_j\ge a^r\sigma_m^{m-r},\qquad
\sum_{j=1}^m\sigma_j^{-2}\ge\sigma_m^{-2},
\]
which gives \eqref{eq:Q1rankbound}. The singular-value perturbation
inequality
\[
|\sigma_m(\B(\xiB))-\sigma_m(\B(z))| \le\|\B(\xiB)-\B(z)\|_2
\]
implies the pointwise distance estimate. For a segment, choose the
nearest point of $L_*$, where the smallest singular value is zero.
Q1 Jacobians are polynomial and thus Lipschitz on the closed cube.
\end{proof}

The divergence estimates also require control of neighborhoods at
faces, edges, and corners. The following geometric argument avoids
assuming that a singular point or segment lies in the interior.

\begin{lemma}[One-sided neighborhoods of the reference cube]
\label{lem:Q1cubemeasure}
Let $\Khat=[0,1]^m$, $c=(1/2,\ldots,1/2)$, and
$\kappa_m=(\sqrt m+1)/2$. Denote by $\omega_m$ the volume of the
unit ball of $\R^m$. Then for any $z\in\Khat$ and
$0<\rho\le\kappa_m$,
\begin{equation}\label{eq:Q1ballmeasure}
\bigl|\Khat\cap B(z,\rho)\bigr|
\ge\omega_m\left(\frac{\rho}{2\kappa_m}\right)^m.
\end{equation}
If $m=3$ and $L_*$ is any segment of length $\ell>0$ in
$\Khat$, then
\begin{equation}\label{eq:Q1tubemeasure}
\left|\{\xiB\in\Khat:\dist(\xiB,L_*)<\rho\}\right|
\ge\frac{\pi\ell}{32\kappa_3^2}\rho^2.
\end{equation}
\end{lemma}
\begin{proof}
The ball $B(c,1/2)$ is contained in the reference cube.
For $z\in\Khat$ and $0<t\le1$, convexity of the cube implies
\[
 B\bigl((1-t)z+tc,t/2\bigr)\subset\Khat.
\]
Each point of this ball is within distance $t\kappa_m$ of $z$.
With $t=\rho/\kappa_m$, its volume proves
\eqref{eq:Q1ballmeasure}. To prove \eqref{eq:Q1tubemeasure},
take $t=\rho/(2\kappa_3)\le1/2$ and apply the same construction
to every $z\in L_*$. The resulting balls have radius $t/2$ and
lie within distance $\rho$ of $L_*$. Their centers sweep a segment
of length $(1-t)\ell\ge\ell/2$. Their union contains a cylinder
of radius $t/2$ and length $\ell/2$, giving the claimed lower bound.
The argument remains valid when the original segment lies on a face
or edge of the cube.
\end{proof}

The neighborhood estimates yield the following general nonintegrability criterion.

\begin{lemma}[Nonintegrability near a point or segment]
\label{lem:Q1divergencecriterion}
Let $f\ge0$ be measurable on $\Khat$, and let $Z$ be either a point
($k=m$) or, when $m=3$, a segment of positive length ($k=2$).
Suppose $f(\xiB)\ge C\dist(\xiB,Z)^{-\alpha}$ almost everywhere
in a sufficiently small neighborhood of $Z$, for constants
$C,\alpha>0$. If $\alpha\ge k$, then
$\int_{\Khat}f\dd\xiB=+\infty$. Equality $\alpha=k$
corresponds to logarithmic divergence.
\end{lemma}
\begin{proof}
By Lemma~\ref{lem:Q1cubemeasure}, the volume of the $\rho$-
neighborhood of $Z$ inside $\Khat$ is bounded below by $c_Z\rho^k$
for all sufficiently small $\rho>0$. For large $h$, the superlevel
set $\{f>h\}$ therefore has volume at least
$c_Z(C/h)^{k/\alpha}$. Removing the null set $Z$ does not change
this volume. The layer-cake identity gives
\[
\int_{\Khat}f\dd\xiB
 =\int_0^\infty|\{\xiB\in\Khat:f(\xiB)>h\}|\dd h
 \ge c_Z C^{k/\alpha}\int_{h_0}^\infty h^{-k/\alpha}\dd h
 =+\infty,
\]
since $k/\alpha\le1$.
\end{proof}

We can now combine these estimates with the algebraic structure of bilinear and trilinear maps.

\begin{lemma}[Divergence at a singular $Q_1$ Jacobian]
\label{lem:Q1singularenergy}
Let $F:[0,1]^m\to\R^d$, with $d\ge m$, be bilinear for $m=2$
or trilinear for $m=3$, and set $\B=D_{\xiB}F$. If
$\rank\B(z)<m$ at any $z$ of the closed reference cube, then
\[
\int_{[0,1]^m}H_p(\B(\xiB))\dd\xiB=+\infty
\]
when $p\ge3/2$ for $m=2$ or $p\ge4/3$ for $m=3$.
\end{lemma}
\begin{proof}
We distinguish the possible Jacobian ranks. If $\B$ is rank
deficient on a region of positive measure, the extended integral is
immediately infinite. The estimates below also handle isolated and
lower-dimensional singular sets.

Case (i): Bilinear maps, rank one.
Write the translated map in local coordinates around $z$ as
\[
F(s,t)=F(z)+bs+ct+dst.
\]
At $z$, the first singular value is positive and the second is zero.
Lemma~\ref{lem:Q1rankestimate} with $(m,r)=(2,1)$ therefore gives
\[
H_p(\B(s,t))\ge C\sigma_2(\B(s,t))^{-(2p-1)}
\ge C'|(s,t)|^{-(2p-1)}
\]
wherever $\B$ has full rank near $z$. Integrability cannot improve
when there is additional rank deficiency. Lemma~\ref{lem:Q1divergencecriterion} 
with $Z=\{z\}$ and $k=2$ gives an infinite integral when $2p-1\ge2$, or $p\ge3/2$. 
The same conclusion holds if additional rank loss occurs nearby.

Case (2): Bilinear maps, rank zero.
Both first derivatives vanish at $z$, so the same shifted bilinear
expansion reduces to $F(s,t)=F(z)+dst$. The two Jacobian columns are
$dt$ and $ds$ and are parallel wherever they are defined. Thus
$\rank\B\le1$ on the whole reference cell, and the extended
integral is infinite.

Case (3): Trilinear maps, rank two.
Here the first two singular values stay bounded below near $z$.
Lemma~\ref{lem:Q1rankestimate} with $(m,r)=(3,2)$ yields
\[
H_p(\B(\xiB))\ge C|\xiB-z|^{-(3p-1)}.
\]
Lemma~\ref{lem:Q1divergencecriterion} with $Z=\{z\}$ and $k=3$ gives divergence when $3p-1\ge3$, or $p\ge4/3$.

Case (4): Trilinear maps, rank one.
A single-point bound would not establish the required exponent in
this case. We prove that rank one forces a nontrivial singular
segment, which reduces the transverse integration dimension.
Write the three Jacobian columns at $z$ as
$\B_i(z)=a_i u$, where $u\in\R^d$ is nonzero and the scalars
$a_i$ are not all zero. If $a_i=0$ for one index, then the
$i$th partial derivative of the trilinear map is independent of
its own reference coordinate. It remains identically zero along the
coordinate-axis segment through $z$. Therefore the Jacobian has rank
at most two on that segment.

Suppose next that all three $a_i$ are nonzero. Choose two coordinates
$i,j$ and hold the third coordinate fixed at its value at $z$.
Because $F$ is separately affine in each reference variable, its
restriction to this coordinate plane has exactly the form
\[
F(s,t)=F(z)+u(a_i s+a_j t)+vst,
\]
where $s,t$ denote increments from $z$ and $v\in\R^d$ is fixed.
Differentiation gives
$\partial_s F=a_i u+vt$ and
$\partial_t F=a_j u+vs$, whence
\begin{equation}\label{eq:Q1wedge}
\partial_sF\wedge\partial_tF
=(a_i s-a_j t)u\wedge v.
\end{equation}
On the line $a_i s=a_j t$, the two selected Jacobian columns are
linearly dependent. With the third column unrestricted, the full
$3$-column Jacobian has rank at most two along that line.

It remains to justify that the line contains a nontrivial segment
inside the reference cube, even when $z$ lies on its
boundary. For every coordinate choose an admissible nonzero
increment sign: at a lower face use the positive sign; at an upper
face use the negative sign; at an interior coordinate either sign is
available. Among the three nonzero signed coefficients
$a_i\operatorname{sign}(s_i)$, two have the same sign. Choose this
pair. Then the equation $a_i s=a_j t$ admits sufficiently small
nonzero increments in their respective admissible directions.
Holding the remaining coordinate fixed defines the desired segment
in the closed cube, with $z$ as an endpoint. This argument handles
interior points, faces, edges, and corners without assuming that the
Jacobian singularity is isolated.

Shrink the segment if needed so that its first singular value stays
bounded below, since it is positive at $z$. On a sufficiently small
neighborhood of a shorter subsegment, Lemma~\ref{lem:Q1rankestimate}
with $(m,r)=(3,1)$ gives
\[
H_p(\B(\xiB))\ge
C\dist(\xiB,L_*)^{-(3p-2)}.
\]
The distance estimate follows from rank deficiency at each point
of $L_*$ and Lipschitz continuity of the Jacobian. Lemma~\ref{lem:Q1divergencecriterion} 
applies to the segment with $k=2$, even when the segment lies on a boundary face or edge. 
The integral therefore diverges when $3p-2\ge2$, or $p\ge4/3$.

Case (5): Trilinear maps, rank zero.
All three first derivatives vanish at $z$. In coordinates $s,t,r$
measured from $z$, the shifted trilinear expansion is
\[
F(s,t,r)=F(z)+Pst+Qsr+Rtr+Tstr,
\qquad P,Q,R,T\in\R^d.
\]
The Jacobian on the three coordinate axes is, respectively,
\[
\B(s,0,0)=[\,0,\ sP,\ sQ\,],\qquad
\B(0,t,0)=[\,tP,\ 0,\ tR\,],\qquad
\B(0,0,r)=[\,rQ,\ rR,\ 0\,].
\]
Each axis intersects the reference cube in a segment of positive
length, even if $z$ lies at its boundary. If one of $P,Q,R$ is
nonzero, at least one of these axis segments contains a nearby
point where the Jacobian has rank one or two. The previous cases
then imply the claimed divergence. If all three vectors vanish,
all Jacobian columns are scalar multiples of $T$, so the Jacobian
has rank at most one throughout the cube. Its extended energy is
therefore infinite everywhere. This exhausts the rank-zero case.
\end{proof}

\subsection*{Sharpness of the exponents}
The following two explicit families show that the critical exponents
cannot be reduced for the class of all $Q_1$ maps in the theorem.
Choose $\M=\I$ and let $0<\varepsilon\le1$. In two dimensions set
\begin{equation}\label{eq:Q1sharp2d}
F_\varepsilon^{(2)}(s,t)
=\bigl(s-t,\,st+\varepsilon(s+t)\bigr),
\qquad
\det DF_\varepsilon^{(2)}=s+t+2\varepsilon>0.
\end{equation}
In three dimensions set
\begin{equation}\label{eq:Q1sharp3d}
\begin{split}
F_\varepsilon^{(3)}(s,t,r)
&=\bigl(s-r,\,t-r,\,st+sr+tr+\varepsilon(s+t+r)\bigr),\\
\det DF_\varepsilon^{(3)}
&=2(s+t+r)+3\varepsilon>0.
\end{split}
\end{equation}
These are respectively bilinear and trilinear maps on the unit cube.
For completeness, the corresponding nodal positions at
$\varepsilon=0$ are as follows. In two dimensions, the images of
$(0,0),(1,0),(0,1),(1,1)$ are
$(0,0),(1,0),(-1,0),(0,1)$, respectively. In three dimensions,
using the ordering
\[
\begin{aligned}
&(0,0,0),(1,0,0),(0,1,0),(1,1,0),\\
&(0,0,1),(1,0,1),(0,1,1),(1,1,1),
\end{aligned}
\]
the corresponding images are
\[
\begin{aligned}
&(0,0,0),(1,0,0),(0,1,0),(1,1,1),\\
&(-1,-1,0),(0,-1,1),(-1,0,1),(0,0,3).
\end{aligned}
\]
The three-dimensional map is not asserted to be globally injective.
The limiting Jacobian at the origin has rank one for
\eqref{eq:Q1sharp2d} and rank two for \eqref{eq:Q1sharp3d}.
In both cases, the singular values not approaching zero are bounded
above and away from zero uniformly on the cube and in
$\varepsilon$: the fixed first output component $s-t$ has gradient
$(1,-1)$ in 2D, while the first two output components in 3D have
gradients $(1,0,-1)$ and $(0,1,-1)$, a fixed full-row-rank matrix.
The remaining singular value is therefore comparable to the positive
Jacobian determinant, with constants independent of $\varepsilon$.

It follows that, with $R=s+t$ in 2D and $R=s+t+r$ in 3D, the
alignment part of the energy density is bounded above by
$C(R+\varepsilon)^{1-mp}$, and its equidistribution part by
$C(R+\varepsilon)^{1-p}$. The constants are uniform in
$\varepsilon\in(0,1]$. To see the relevant integrability directly,
partition the positive coordinate cube by the value of $R$. Its
cross-sectional measure is bounded by $CR^{m-1}$ near the origin.
Thus the potentially singular alignment integral is controlled by
\[
C\int_0^{m}\!R^{m-1}(R+\varepsilon)^{1-mp}\dd R.
\]
For $0<\varepsilon\le1$ this stays uniformly finite when
$m-1+1-mp>-1$, equivalently $p<p_m:=1+1/m$, consistent with
\eqref{eq:Q1criticalp}. The equidistribution term is less
singular for $p>1$. Nevertheless,
$\sigma_{\min}(DF_\varepsilon^{(m)}(0))\to0$ as
$\varepsilon\downarrow0$. At the critical exponent the limiting
alignment integral diverges logarithmically. At subcritical exponents, setting $\varepsilon=0$ already
gives a rank-deficient map with finite exact energy for $1<p<p_m$, 
since assigning $+\infty$ at the isolated singular point does not alter the Lebesgue integral. 
The positive-$\varepsilon$ families show that the failure also occurs along regular, positively oriented maps. 
Therefore no uniform positive singular-value bound can hold on bounded-energy sublevels for $1<p<p_m$.

\section{Polynomial degree of the trilinear Jacobian determinant}\label{app:degree}

For completeness, consider the trilinear map
\[
F(\xi,\eta,\zeta)
=\sum_{i,j,k\in\{0,1\}}
N_{ijk}(\xi,\eta,\zeta)\x_{ijk}.
\]
The derivative $F_\xi$ is bilinear in $(\eta,\zeta)$ and independent of $\xi$; $F_\eta$ is 
bilinear in $(\xi,\zeta)$ and independent of $\eta$; and $F_\zeta$ is bilinear in $(\xi,\eta)$ 
and independent of $\zeta$.  In the determinant
\[
\det[F_\xi,F_\eta,F_\zeta],
\]
each reference variable can therefore appear with degree at most two.  
Thus the determinant lies in the tensor-product polynomial space $\mathbb Q_2$, which has dimension $3^3=27$.  
This is the algebraic basis for \eqref{eq:hexbern} and for the validity method of \cite{JohnenWeillRemacle2017}.

For a bilinear quadrilateral, the corresponding cancellation
\[
\det(\bm d,\bm d)=0
\]
removes the mixed $\xi\eta$ term from the determinant, yielding the affine expression \eqref{eq:quadq}.

\section{Scaling of the adaptation density}\label{app:scaling}

Under the metric scaling $\M\mapsto\chi\M$ with $\chi>0$, one has
\[
\A\mapsto\chi\A,
\qquad
\J\mapsto\chi^{-1}\J,
\qquad
D\mapsto\chi^{-m}D.
\]
Both terms in \eqref{eq:G} scale as
\[
\chi^{-m(p-1)/2}.
\]
Hence $I_{\EA}$ and its nodal gradient scale by the same factor.
The intrinsic prefactor \eqref{eq:PiScale} scales as $\chi^{m(p-1)/2}$ and
therefore makes the mesh velocity due to the adaptation energy alone
invariant under constant rescaling of the metric.  For the total
barrier-augmented energy, this invariance does not hold automatically:
$I_{\val}$ is unaffected by $\M\mapsto\chi\M$ when $\mu$, $w_c$, and the
certificate normalization are fixed.  To retain invariance, the barrier
weight $\mu$ must be rescaled as $\chi^{-m(p-1)/2}$, or an equivalent
normalization must be imposed.

\end{document}